\documentclass[11pt,a4paper]{article}
\usepackage[margin=28mm]{geometry}
\usepackage{amsmath,amssymb,amsthm,mathtools,booktabs,array,microtype,mathrsfs}
\usepackage[T1]{fontenc}
\usepackage{lmodern}
\usepackage{tikz,flafter}
\usetikzlibrary{arrows.meta,calc,positioning,decorations.pathreplacing}
\usepackage[colorlinks=true,linkcolor=blue,citecolor=blue,urlcolor=blue]{hyperref}
\hypersetup{pdftitle={Logarithmic corrections for four-dimensional loop percolation}}
\numberwithin{equation}{section}
\newtheorem{theorem}{Theorem}[section]
\newtheorem{lemma}[theorem]{Lemma}
\newtheorem{proposition}[theorem]{Proposition}
\newtheorem{corollary}[theorem]{Corollary}
\theoremstyle{definition}

\theoremstyle{remark}
\newtheorem{remark}[theorem]{Remark}
\newcommand{\Pp}{\mathbb P}
\newcommand{\E}{\mathbb E}
\newcommand{\Z}{\mathbb Z}
\newcommand{\R}{\mathbb R}

\newcommand{\capc}{\operatorname{cap}}
\newcommand{\diam}{\operatorname{diam}}
\newcommand{\rad}{\operatorname{rad}}
\newcommand{\Cl}{\operatorname{Cl}}
\newcommand{\one}{\mathbf 1}

\newcommand{\tr}{\operatorname{Tr}}

\title{Logarithmic corrections for four-dimensional loop percolation}
\author{Quirin Vogel$^1$}
\date{2 November 2026}
\begin{document}
\maketitle
\begin{abstract}
We study the one-arm probability of the discrete random-walk loop soup on
$\mathbb{Z}^4$ at small positive intensity. Using the small-loop crossing and
macroscopic attachment estimates of Gao, we prove that the expected capacity
of the origin cluster formed by loops contained in a box grows at most as
the square of the logarithm of the box radius. Consequently, there are
constants $\alpha_0,c,C,\beta>0$ such that, for $0<\alpha\le\alpha_0$,
\[
 c\alpha R^{-2}[\log(R+2)]^{\beta\alpha}
 \le \mathbb{P}_\alpha(0\leftrightarrow B_R^c)
 \le C\alpha R^{-2}[\log(R+2)]^2.
\]
The interval $(0,\alpha_0]$ is strictly subcritical. The upper bound rests on
a simultaneous capacity estimate for neighborhoods of loop ranges, followed
by an expansion over chains in which sufficiently small intermediate loops
have been removed. A deterministic bound on the variation of the surviving
scale sequences makes the expansion summable. The positive logarithmic-power lower bound is already known from work of
Chang and Sapozhnikov; we give a proof using macroscopic attachments and
a layer-cake argument, with constants uniform at small intensity. In particular,
the one-arm exponent equals two, and the logarithm of its correction has
order $\log\log R$, for every fixed intensity in this interval. We also derive polynomial tail exponent one for the
cardinality and capacity of the full cluster, determine their positive
moment thresholds, and give an application to random wires.
\end{abstract}
\centerline{\textit{$^1$ University of Klagenfurt,
Department of Statistics, Klagenfurt, Austria}}
\centerline{\textsc{quirin.vogel@aau.at}}
\setcounter{tocdepth}{1}
{\small\tableofcontents}
\clearpage

\section{Introduction and main result}\label{sec:introduction}

The loop soup on $\Z^4$ has two sources of long connections. A single long
loop may reach a distant box, while a sequence of intersecting loops may
produce the same connection without any one loop spanning the full distance.
The second mechanism makes a direct first-moment estimate on loop chains
ineffective: a long loop can meet many much smaller loops, and counting all
such contacts separately loses the screening caused by their overlap.

Capacity is a natural way to measure the set available for a further
connection. In dimension four, a random-walk segment of duration $n$ has
expected capacity of order $n/\log n$. This logarithmic saving is substantial,
but it must survive both the addition of nearby small-loop clusters and
the selection of loops through a connectivity event. The purpose of this
paper is to establish an upper bound that retains enough of that saving to
sum the resulting chains at sufficiently small intensity.

\subsection{The discrete loop soup}

We work throughout on $\Z^4$. For $x\in\Z^4$ and $r\ge0$, write
\[
 |x|=\max_{1\le a\le4}|x_a|,\qquad
 B(x,r)=\{y\in\Z^4:|y-x|\le r\},\qquad B_r=B(0,r).
\]
All diameters and radii use this norm. For finite sets $A,F\subset\Z^4$,
$A+F=\{a+f:a\in A,f\in F\}$. The internal vertex boundary of $B(x,r)$
is denoted by $S(x,r)$, and $S_r=S(0,r)$. Thus the notation for a sphere
also makes sense at a noninteger radius. An annular crossing means that a
connected set meets both boundary sets; equivalent formulations using
entrance and exit from boxes differ only by a bounded lattice rounding.

Let $S=(S_k)_{k\ge0}$ be nearest-neighbor simple random walk. Its law and
expectation when started at $x$ are $\Pp_x$ and $\E_x$, and
\[
 P(x,y)=\tfrac18\one_{\{|x-y|_1=1\}},\qquad
 p_n(x,y)=\Pp_x(S_n=y),\qquad
 G(x,y)=\sum_{n=0}^\infty p_n(x,y).
\]
Here $|\cdot|_1$ is used only to specify nearest neighbors. For every even
$n\ge2$, let $\Pp_n^{\mathrm{br}}$ be the probability law
$\Pp_0(\,\cdot\mid S_n=0)$ on paths of length $n$, with expectation
$\E_n^{\mathrm{br}}$. More explicitly, for an event $F$ determined by
the first $k\le n$ steps,
\begin{equation}\label{eq:model-bridge}
 \Pp_n^{\mathrm{br}}(F)
 =\frac{\E_0[\one_F p_{n-k}(S_k,0)]}{p_n(0,0)}.
\end{equation}
There is no bridge law at odd $n$, since the walk is bipartite.
For a path of length $n$ its vertex range is
$U_n=\{S_0,\ldots,S_n\}$.

A rooted loop is a nearest-neighbor path
$\dot\ell=(x_0,\ldots,x_n)$ with $n\ge2$ and $x_n=x_0$.
Its duration is $N(\dot\ell)=n$. Two rooted loops represent the same
unrooted loop if one is obtained from the other by a circular time shift.
Time reversal is not included in this equivalence relation. Write
$\pi\dot\ell$ for the resulting unrooted loop and
$U(\ell)=\{x_0,\ldots,x_{n-1}\}$ for its range. Both $N$ and $U$
are well defined on unrooted loops. Let $\mathscr L$ be the countable set
of unrooted loops, equipped with its power-set $\sigma$-field.

The measure $\mu$ on $\mathscr L$ is the pushforward of rooted loop
measure assigning mass $n^{-1}8^{-n}$ to each rooted $n$-step loop. In
particular, for every nonnegative function $F$ on $\mathscr L$,
\begin{equation}\label{eq:model-rooted}
 \int F(\ell)\,\mu(d\ell)
 =\sum_{n=2}^\infty\frac1n
   \sum_{x\in\Z^4}
   \E_x\!\left[F\bigl(\pi(S_0,\ldots,S_n)\bigr);S_n=x\right].
\end{equation}
The value inside the expectation is used only on the return event.
We write $\mu_n$ for the restriction of $\mu$ to duration $n$.
Formula~\eqref{eq:model-rooted}, rather than a choice of representative,
fixes the weight of periodically repeated loops. It is the normalization
used in Gao~\cite[Section~7.1]{Gao} and in the loop-measure construction
of Lawler and Limic~\cite[Chapter~9]{LL}.

For $\alpha\ge0$, the loop soup $\mathcal L_\alpha$ is the Poisson random
counting measure on $\mathscr L$ of intensity $\alpha\mu$.
Its law and expectation are $\Pp_\alpha$ and $\E_\alpha$.
Thus the numbers of occurrences of distinct loop shapes are independent
Poisson variables. Distinct occurrences of the same shape are treated as
distinct objects whenever a chain or a disjoint witness is specified.
For a set $D\subset\Z^4$, $\mathcal L_\alpha^D$ denotes the restriction
to loops with $U(\ell)\subset D$. A duration cutoff $L$ means the
restriction $N(\ell)<L$; a diameter cutoff $s$ means
$\diam U(\ell)\le s$. These are different restrictions.

Two vertices are connected if they coincide or if there is a finite chain
of soup occurrences whose successive ranges intersect, with the first
range containing the first vertex and the last range containing the
second. Denote the full origin cluster by $\mathcal C_\alpha(0)$.
The finite-box cluster $K_R$ is the origin cluster formed using only
$\mathcal L_\alpha^{B_R}$. In particular, $0\in K_R\subset B_R$ even
if no loop meets the origin. The set $K_R$ is generally smaller than
$\mathcal C_\alpha(0)\cap B_R$.

For a finite set $A$, ordinary capacity means the full-space capacity
\[
 \capc(A)=\sum_{x\in A}
   \Pp_x(S_k\notin A\text{ for every }k\ge1),
 \qquad \capc(\varnothing)=0.
\]
The potential-theoretic identities used below are proved in
Section~\ref{sec:potential}. We set
\begin{equation}\label{eq:model-observables}
 A_\alpha(R)=\E_\alpha[\capc(K_R)],\qquad
 p_\alpha(R)=\Pp_\alpha(\mathcal C_\alpha(0)\cap B_R^c\ne\varnothing).
\end{equation}
In particular, capacity is not computed for a walk killed on leaving
$B_R$. Write $c_0=\capc\{0\}=G(0,0)^{-1}$, and define the percolation
threshold by
\[
 \alpha_c=\inf\{\alpha\ge0:
       \Pp_\alpha(|\mathcal C_\alpha(0)|=\infty)>0\},
 \qquad\inf\varnothing=\infty.
\]

\subsection{The theorem}

\begin{theorem}\label{thm:main}
There are constants $\alpha_0>0$, $\beta>0$, and
$c_A,C_A,c_p,C_p\in(0,\infty)$, independent of $\alpha$ and $R$, such
that $\alpha_0<\alpha_c$ and, for every $R\ge1$,
\begin{equation}\label{eq:maincap}
 c_A[\log(R+2)]^{\beta\alpha}
 \le A_\alpha(R)
 \le c_0+C_A\alpha[\log(R+2)]^2,
 \qquad 0\le\alpha\le\alpha_0.
\end{equation}
Moreover,
\begin{equation}\label{eq:mainarm}
 c_p\alpha R^{-2}[\log(R+2)]^{\beta\alpha}
 \le p_\alpha(R)
 \le C_p\alpha R^{-2}[\log(R+2)]^2,
 \qquad 0<\alpha\le\alpha_0.
\end{equation}
Consequently, for every fixed $\alpha\in(0,\alpha_0]$,
\begin{equation}\label{eq:exponent-conclusion}
 p_\alpha(R)=R^{-2+o(1)},\qquad
 \log\bigl(R^2p_\alpha(R)\bigr)
     =\Theta_\alpha(\log\log R)\quad(R\longrightarrow\infty).
\end{equation}
\end{theorem}

The constants in this statement are existence constants, not numerical
estimates of the size of the subcritical interval. An explicit choice of
$\alpha_0$ in terms of the constant in the loop transition estimate is
given in Section~\ref{sec:conclusion}. The symbol $\Theta_\alpha$ in
\eqref{eq:exponent-conclusion} means two-sided bounds by positive multiples
whose values may depend on the fixed intensity. The theorem establishes
the order of the logarithmic correction; it does not establish a limit for
$\log(R^2p_\alpha(R))/\log\log R$, or extend the conclusion to the whole
interval $(0,\alpha_c)$.
The proof also yields a sharper, intensity-dependent upper bound for
the logarithmic correction, implying that the exponent of 2 is not sharp at small intensities; see Remark~\ref{end:refined-log}.

The same expansion also controls the number of vertices in a finite-box
cluster. Together with the macroscopic attachment estimate, it gives
tail and moment statements for the full cluster. For the following
corollary, its radius means
$\rho_\alpha:=\max\{|x|:x\in\mathcal C_\alpha(0)\}$.
The cluster is almost surely finite in the stated intensity interval,
so its full-space capacity is defined by the finite-set formula above.

\begin{corollary}[Cluster-size and capacity tails]\label{cor:cluster-tails}
There are constants $c,C>0$ and $t_{\mathrm{cap}}>1$, independent of
$0<\alpha\le\alpha_0$, such that, for $t\ge t_{\mathrm{cap}}$,
\begin{align}
 c\alpha t^{-1}[\log(t+2)]^{\beta\alpha-1}
 &\le\Pp_\alpha\bigl(\capc(\mathcal C_\alpha(0))>t\bigr)
 \le C\alpha t^{-1}[\log(t+2)]^2,
 \label{eq:full-capacity-tail}\\
 c\alpha t^{-1}[\log(t+2)]^{\beta\alpha-1}
 &\le\Pp_\alpha\bigl(|\mathcal C_\alpha(0)|>t\bigr)
 \le C\alpha t^{-1}[\log(t+2)]^3.
 \label{eq:full-volume-tail}
\end{align}
Consequently, for each fixed $\alpha\in(0,\alpha_0]$, both tail
probabilities in these displays are $t^{-1+o(1)}$ as $t\to\infty$.
For every $s>0$,
\begin{align}
 \E_\alpha[\rho_\alpha^s]<\infty
 &\quad\Longleftrightarrow\quad s<2,
 \label{eq:radius-moments}\\
 \E_\alpha[|\mathcal C_\alpha(0)|^s]<\infty
 &\quad\Longleftrightarrow\quad s<1,\notag\\
 \E_\alpha[\capc(\mathcal C_\alpha(0))^s]<\infty
 &\quad\Longleftrightarrow\quad s<1.
 \label{eq:size-capacity-moments}
\end{align}
\end{corollary}
The proof is given in Section~\ref{sec:cluster-tails}. The conclusion
concerns the cardinality of the entire cluster, not only its radius;
the infinite first moment of the cardinality was already proved by
Chang and Sapozhnikov~\cite[Theorem~1.4(a)]{CS16}.

\subsection{Relation to earlier work and the role of dimension four}

The continuum counterpart is the Brownian loop soup of Lawler and
Werner~\cite{LW04}, whose planar cluster boundaries underlie the
construction of conformal loop ensembles by Sheffield and
Werner~\cite{SW12}. Lawler and Trujillo Ferreras~\cite{LTF07}
introduced the random-walk loop soup and its approximation by Brownian
loops, while Le Jan and Lemaire~\cite{LJL13} developed loop-cluster
percolation on graphs.

Chang and Sapozhnikov~\cite{CS16} established a nontrivial phase transition
for the discrete loop soup on $\Z^d$, $d\ge3$. Their
Theorem~1.3 gives a polynomial upper bound at sufficiently small intensity,
with an exponent tending to $d-2$ as the intensity tends to zero.
In dimension four, they also proved
\[
 p_\alpha(R)\ge c(\alpha)R^{-2}(\log R)^{\varepsilon(\alpha)},
 \qquad \varepsilon(\alpha)>0,
\]
and a corresponding lower bound for the expected capacity of the
finite-box origin cluster; see their Theorem~1.6 and Lemma~6.2.
Thus a positive logarithmic correction to the one-arm probability was
already known. The upper bound in the present paper supplies the other
side, on a fixed interval of positive intensities, and identifies the
polynomial decay exponent there as exactly two. Our proof of the lower
bound records constants uniformly down to zero intensity.

The contrast with higher dimensions is useful. Put
\[
 \alpha_\#=\sup\{\alpha\ge0:
          \E_\alpha[|\mathcal C_\alpha(0)|]<\infty\}.
\]
For $d\ge5$, Chang and Sapozhnikov~\cite[Theorem~1.4]{CS16} proved
$\alpha_\#>0$ and one-arm bounds of order $R^{2-d}$ for
$\alpha<\alpha_\#$. Vogel~\cite[Theorems~1 and~2]{Vog26} obtained
asymptotic equivalents, with identified leading constants, for both the
one-arm probability and the two-point function in this interval. His
proof describes a long connection by attaching one large loop to localized
clusters. In dimensions three and four, however,
$\alpha_\#=0$~\cite[Theorem~1.4(a)]{CS16}; estimates based on the
expected cardinality of the unrestricted cluster cannot be transferred
directly. The finite-box expected capacity $A_\alpha(R)$ is the quantity
we control instead.

The supercritical regime has a different finite-cluster behavior:
Chang~\cite{Cha17} proved exponential decay of the diameter of the
finite origin cluster for every $\alpha>\alpha_c$. Recent work of
Betz, Drewitz, Klippel and M\"onch~\cite{BDKM26} also proves that the
discrete threshold is strictly larger than $1/2$ in $d\ge5$, further
separating discrete loop percolation from its metric-graph counterpart.

The logarithm in the four-dimensional range-capacity scale is central
to this distinction. Asselah, Schapira and Sousi~\cite[Theorem~1.1 and
(1.4)]{ASS} proved the almost-sure asymptotic and recorded its mean version,
\[
 \frac{\log n}{n}\capc(U_n)\longrightarrow\frac{\pi^2}{8}
 \quad\text{almost surely},
 \qquad
 \E_0[\capc(U_n)]\sim\frac{\pi^2}{8}\frac n{\log n}.
\]
The same dimension is the threshold for summability of the square of the
Green function: $G(0,x)\asymp(1+|x|)^{2-d}$ gives
$\sum_{|x|\le R}G(0,x)^2\asymp\log R$ in $d=4$, whereas
$\sum_xG(0,x)^2<\infty$ for $d>4$. This is one concrete sense in which
four-dimensional intersections are marginal. Our proof must preserve the
factor $1/\log n$ when small loops are added to a long loop, then control
the accumulation of errors over logarithmically many duration scales.

Logarithmic corrections also arise in other four-dimensional random-walk
models. For example, Bauerschmidt, Brydges and Slade~\cite[Theorem~1.1]{BBS}
proved that the susceptibility of continuous-time weakly self-avoiding
walk, at sufficiently small positive coupling, diverges as the inverse
distance to criticality times the one-quarter power of its logarithm.
This comparison places the phenomenon in a broader context; it does not
identify a logarithmic exponent for loop percolation, whose observable
and mechanism here are different.

Finally, the questions in Chang and Sapozhnikov~\cite[Section~8]{CS16}
should be distinguished. They ask both about a logarithmic-power
correction to $R^{-2}$ and about a conditional typical loop distance
of order $\log\log R$, where loop distance counts the loops in a shortest
connecting chain. The present theorem concerns the one-arm probability.
It does not establish the loop-distance assertion, a limiting logarithmic
exponent, or bounds throughout the full subcritical interval.

\subsection{Connections with Gaussian fields and random currents}

Loop soups are also connected to Gaussian fields through isomorphism
theorems. For the continuous-time soup, including loops that make no
jumps, Le Jan~\cite[Chapter~5]{LeJan} proved that the occupation field at
intensity $1/2$ has the law of one half of the square of the Gaussian
free field.
Lupu~\cite[Theorem~1 and Proposition~2.1]{Lupu} strengthened this to a
coupling with cluster geometry: the sign of the field is constant on
each discrete loop cluster, and on the associated metric graph the loop
clusters at intensity $1/2$ are exactly the sign clusters. Here the metric
graph is obtained by replacing each lattice edge by a continuous
interval. The distinction between these statements matters: connections
inside the intervals can merge clusters that remain separate in the
discrete soup. Lupu's nonpercolation theorem~\cite[Theorem~2]{Lupu}
implies $\alpha_c\ge1/2$ on $\Z^d$, $d\ge3$; Chang, Du and
Li~\cite{CDL} proved that the threshold of the metric-graph soup is
exactly $1/2$. The small-intensity discrete estimate proved here is
therefore a different regime from the Gaussian-field identification at
intensity $1/2$.

Capacity is central to this adjacent theory as well. Drewitz, Pr\'evost
and Rodriguez~\cite[Theorem~1.1]{DPRcap} determined capacity laws for
metric-graph Gaussian free field clusters under their nonpercolation
assumption, and used capacity estimates for large loops to obtain sharp
one-arm exponents on a class of transient graphs that includes
$\Z^4$~\cite[Theorem~1.1]{DPRarm}. For the critical metric-graph model
on $\Z^4$, the lower bound of Ding and Wirth~\cite{DW} and the matching
upper bound of Cai and Ding~\cite[Theorem~1.1]{CDlow} give one-arm
probability of order $R^{-1}$. This should be distinguished from the
exponent two at small discrete-soup intensity in
Theorem~\ref{thm:main}. Gao~\cite[Theorems~1.9--1.10 and
Section~7]{Gao} developed a multiscale capacity argument for discrete
and metric-graph soups in dimensions three, four and five throughout
the small-loop decay regime. The small-loop decay and macroscopic
attachment estimates used here come from that work; our additional
task is to control the expected capacity from above when all loop
durations are present.

A further connection concerns random currents. In the continuous-time
representation at soup intensity $1/2$, conditioning on its occupation
field gives the source-free Ising
random-current law for edge traversal counts, with couplings determined
by that field; Lupu and Werner~\cite[Section~2]{LW16} explain this
representation and its relation to the Gaussian free field and
FK-Ising model. Ding and Li~\cite[Theorem~2 and Section~1.2]{DL18}
used this correspondence to study chemical distances in the planar
discrete soup killed at the boundary of a box of side order $N$:
with probability bounded away from zero, there is a crossing between
two fixed macroscopic box boundaries of length at most $N^{1+o(1)}$.
Their observable and intensity differ from those of the present paper.
Our application in Section~\ref{sec:wires} instead concerns the
zero-vertex-potential specialization of the random wire family of
Benassi and Ueltschi~\cite{BU}. Its edge support has the law of the
discrete soup at intensity $\alpha=q/2$, where $q>0$ is the wire loop
weight. The theorem therefore supplies bounds for small positive $q$;
it does not reach the Gaussian case $q=1$.

\subsection{Structure of the argument}

The external loop-soup estimates are stated precisely in
Section~\ref{sec:inputs}. They imply, with constants uniform for bounded
small intensities,
\[
 p_\alpha(R)\asymp \alpha R^{-2}A_\alpha(R).
\]
We prove this comparison at the end of the paper, together with the lower
bound in~\eqref{eq:maincap}. The main work is the upper bound for
$A_\alpha(R)$.

First, we estimate the capacity of a loop range after adding a small
spatial neighborhood. Geometric killing is useful here because a loop
duration band $[M,2M)$ naturally tests capacity at killing time of order
$M$. Section~\ref{sec:potential} derives the corresponding loop-measure
estimate. Section~\ref{sec:concentration} shows that, outside a class of
bridge probability smaller than every negative power of the duration, a
loop satisfies the required capacity bound simultaneously at all killing
times. The argument uses an elementary lower-capacity estimate to control
hitting probabilities on a thin neighborhood, and then concentration for
independent translated pieces of a free walk. The latter use of block
independence and subadditivity follows the method of Popov and
Vogel~\cite[Section~2]{PV}; the capacity estimates needed here are proved
below.

Second, small-loop connections are absorbed into a neighborhood of each
retained loop. To avoid conditioning on how a loop was discovered, the
resulting estimate is proved for every deterministic parent loop.
A good loop carries weight $N/\log(N+2)$ and an exceptional loop carries
weight $N^2$. The exceptional bridge probability pays for this larger
weight. With the duration type $i=\lfloor\log_2N\rfloor$, the resulting
transition bound has the form
\[
 K_t(i,j)=\frac{t}{j}\max\{1,i/j\},\qquad t=C_K\alpha.
\]
Here $C_K$ is independent of the two durations and of the intensity in a
fixed interval. This estimate is established in
Section~\ref{sec:kernel}.

Third, we reorganize a connecting chain before applying this bound.
Whenever all interior types of a subchain are much smaller than both
endpoint types, we remove those interior loops and retain their connection
as a small-loop witness. The remaining word satisfies a condition on
every interval of indices. Factorial Mecke and the BK inequality turn
disjoint witnesses into products of fresh connection probabilities;
Section~\ref{sec:chains} gives the complete finite-volume argument.
The interval condition bounds the total variation of the logarithms of
the retained types. In Section~\ref{sec:words}, this variation estimate
controls the sum over all word lengths. Finally,
Section~\ref{sec:conclusion} removes the duration restriction, and
Section~\ref{end:arms} derives the one-arm and lower-capacity bounds,
and Section~\ref{sec:cluster-tails} proves the full-cluster tail and
moment statements. Section~\ref{sec:wires} transfers the estimates to the
zero-vertex-potential member of the random-wire family of Benassi and
Ueltschi~\cite{BU}.

\paragraph{Constants and notation.}
We fix throughout
\begin{equation}\label{eq:fixed-parameters}
 \bar\alpha=\tfrac14,\qquad \theta=\tfrac1{100},\qquad D=10.
\end{equation}
The constant $D$ is the prefactor in the neighborhood radius introduced
in Section~\ref{sec:concentration}. Unless a dependence is indicated,
$c,C\in(0,\infty)$ denote constants that may change between inequalities
and depend only on these fixed parameters and the lattice. A subscript,
such as $C_Q$ or $C_\delta$, records an additional dependence. Constants
never depend on a box radius, a loop duration, or a final duration cutoff.
Named constants such as $C_0,C_K,t_0,\alpha_*$ and $\alpha_0$ retain their
values once chosen. In particular, the good-loop threshold $C_0$ will be
chosen before the power of the desired exceptional-probability estimate.
The notation $f\asymp g$ means two inequalities with positive constants
under this convention. Every logarithm is natural unless a base is shown.

\section{Random-walk and loop-soup inputs}\label{sec:inputs}

This section states the results imported from the literature and records
the consequences of metric-graph projection that will be used later.
In particular, the duration-truncated cluster estimate below is a
deduction from Gao's metric-graph theorem, rather than its literal
statement. No estimate on the capacity of an entire loop-soup cluster
is assumed.

\subsection{Random-walk estimates}

The local central limit theorem and Green-function estimates for simple
random walk give constants $c_{\mathrm{rw}},C_{\mathrm{rw}}>0$ such that
\begin{equation}\label{in:walk}
 \begin{gathered}
 p_n(x,y)\le C_{\mathrm{rw}}(n+1)^{-2},\qquad
 p_{2n}(0,0)\ge c_{\mathrm{rw}}(n+1)^{-2},\\
 c_{\mathrm{rw}}(1+|x-y|)^{-2}
 \le G(x,y)\le C_{\mathrm{rw}}(1+|x-y|)^{-2}.
 \end{gathered}
\end{equation}
We also use the following exit and interior Green estimates. If
$\tau_{B(x,r)}=\inf\{k\ge0:S_k\notin B(x,r)\}$, then
\begin{equation}\label{in:exit}
 \sup_{y\in B(x,r)}\Pp_y(\tau_{B(x,r)}>u)
 \le C_{\mathrm{rw}}\exp[-c_{\mathrm{rw}}u/(r+1)^2],
 \qquad r\ge1,\ u\ge0.
\end{equation}
For $G_F(x,y)=\E_x[\sum_{k<\tau_F}\one_{\{S_k=y\}}]$, one has
\begin{equation}\label{in:killed-green}
 \inf_{x,y\in B_{2r}}G_{B_{4r}}(x,y)
 \ge c_{\mathrm{rw}}r^{-2},\qquad r\ge1.
\end{equation}
These are standard estimates; see Lawler and Limic~\cite[Chapters~2, 4,
and 6]{LL}. The interior bound follows from the Green estimates and
Harnack inequality in a domain with a fixed separation between the
interior and the boundary. The exponential exit bound can also be seen
directly from the heat-kernel upper bound: choose a fixed large $K$ so
that
\[
 \sup_{y\in B(x,r)}\Pp_y(S_{\lceil K(r+1)^2\rceil}\in B(x,r))
 \le |B(x,r)|C_{\mathrm{rw}}[K(r+1)^2]^{-2}\le\tfrac12,
\]
and iterate at multiples of $\lceil K(r+1)^2\rceil$ using the Markov
property. The upper bound on box capacity,
\begin{equation}\label{in:box-capacity}
 \capc(B_r)\le C_{\mathrm{rw}}(r+1)^2,
\end{equation}
follows from the probability-measure energy formula for capacity and
the Green lower bound; that formula is recalled in
Section~\ref{sec:potential}.

The range-capacity input is the mean asymptotic recorded in
Asselah, Schapira and Sousi~\cite[(1.4)]{ASS}. We use only its upper-bound
consequence: there is a fixed $C_{\mathrm{range}}<\infty$ with
\begin{equation}\label{in:range-mean}
 \E_0[\capc(U_n)]\le
 C_{\mathrm{range}}\frac{n+1}{\log(n+2)},\qquad n\ge1.
\end{equation}
That mean asymptotic is derived in the cited work from earlier
capacity and nonintersection estimates. The stronger almost-sure theorem
of Asselah, Schapira and Sousi is discussed in
Section~\ref{sec:potential}: it could also supply the fixed-confidence
lower estimate needed there. For completeness we give an elementary
energy proof of that estimate. No range-capacity estimate conditioned on
a loop-percolation event is assumed.

\subsection{Metric-graph decay and its discrete consequence}

Let $\widetilde{\Z}^4$ be the metric graph obtained by replacing each
nearest-neighbor edge of $\Z^4$ with an interval. We use the metric-graph
loop measure and normalization of Gao~\cite[Section~7.1]{Gao}.
Projection records successive visits to distinct lattice vertices and
discards loops with no nontrivial discrete trace. Under this projection,
the metric-graph soup of intensity $\alpha$ gives precisely
$\mathcal L_\alpha$ with the normalization
\eqref{eq:model-rooted}. In particular, the intensity is unchanged.
Excursions of a metric loop away from its projected vertex range stay
within the incident edge intervals. Consequently, after fixing the unit
edge convention, a projected range of diameter at most $s$ has a lift
of diameter at most $s+2$. We use boxes from the natural lattice
embedding also for the metric graph; harmless fixed changes of edge or
norm convention change only constants.

\begin{proposition}[Gao's truncated crossing estimate]\label{in:metric-tail}
Fix $\bar\alpha=1/4$. There are constants
$c_{\mathrm{met}},C_{\mathrm{met}}>0$ such that, for every $L\ge1$,
integer $q\ge1$, lattice vertex $x$, and $0\le\alpha\le\bar\alpha$,
\[
 \Pp_\alpha\!\left(
 \widetilde S(x,L)\leftrightarrow\widetilde S(x,qL)
 \text{ using metric loops of diameter at most }L\right)
 \le C_{\mathrm{met}}e^{-c_{\mathrm{met}}q}.
\]
Here $\widetilde S$ denotes the corresponding metric-graph sphere.
\end{proposition}
This is Gao~\cite[Lemma~7.1]{Gao}, applied at intensity $\bar\alpha$;
the assertion for smaller intensities follows by the monotone coupling
obtained by thinning that Poisson soup.

The exponential decay of a diameter-truncated soup has an earlier
analogue in Vogel~\cite[Lemma~3.8]{Vog26}: for $d\ge5$ and
$\alpha<\alpha_\#$, the probability of reaching distance $r$ using only
loops of diameter at most $L$ is at most $C\exp(-cr/L)$.
That argument uses finite expected cluster size. It therefore does not
provide the four-dimensional estimate, since $\alpha_\#=0$ in $d=4$.
The finite-size annular estimates of Chang and
Sapozhnikov~\cite[Lemma~4.1]{CS16}, together with a finite-range block
argument after truncation, suggest a small-intensity route to this type of decay: after a diameter
cutoff, well-separated block events depend on disjoint collections of
loops. We do not use that alternative here. Gao's result is the input used here because it gives
uniform cutoff estimates in the metric graph for every fixed intensity
ceiling below $1/2$; projection then yields the required discrete
four-dimensional bound.

For $L\ge2$, let $C_L(x)$ be the component of $x$ in the discrete soup
restricted to $N(\ell)<L$, including $x$ itself. It depends on $\alpha$,
which is suppressed in the notation. Set
$\rad C_L(x)=\sup\{|y-x|:y\in C_L(x)\}$.

\begin{lemma}\label{in:micro-tail}
There are constants $c_{\mathrm{mic}},C_{\mathrm{mic}}>0$, uniform for
$0\le\alpha\le\bar\alpha$, such that
\begin{equation}\label{eq:micro-tail}
 \Pp_\alpha(\rad C_L(x)>r)
 \le C_{\mathrm{mic}}e^{-c_{\mathrm{mic}}r/(L+1)},
 \qquad L\ge2,\ r\ge0.
\end{equation}
In particular, $C_L(x)$ is finite almost surely. For another constant
$C<\infty$ and a possibly smaller $c>0$,
\begin{align}\label{in:micro-moments}
 \E_\alpha[|C_L(x)|]&\le C(L+1)^4,\notag\\
 \E_\alpha[|C_L(x)|;\rad C_L(x)>H]
 &\le C(L+1)^4
       \left(1+\frac{H}{L+1}\right)^4
       e^{-cH/(L+1)},\qquad H\ge0.
\end{align}
\end{lemma}
\begin{proof}
Every discrete loop of duration less than $L$ has diameter at most $L$.
In the projection coupling it is therefore the trace of a metric loop
of diameter at most $L+2$. A discrete chain from $x$ to outside $B(x,r)$
lifts to a connected metric set with the same endpoints. If
$r\ge2(L+2)$, that set crosses the spheres of radii $L+2$ and
$q(L+2)$, where $q=\lfloor r/(L+2)\rfloor$.
Proposition~\ref{in:metric-tail} bounds its probability by
$C_{\mathrm{met}}e^{-c_{\mathrm{met}}q}$, which implies
\eqref{eq:micro-tail}. For smaller $r$, increasing the prefactor gives
the same inequality.

Write $X=\rad C_L(x)$ and $a=L+1$. Since a component of radius $X$
is contained in $B(x,X)$, $|C_L(x)|\le C(X+1)^4$. For $H\ge0$,
the layer-cake identity gives
\[
 \E_\alpha[(X+1)^4;X>H]
 =(H+1)^4\Pp_\alpha(X>H)
   +\int_H^\infty 4(u+1)^3\Pp_\alpha(X>u)\,du.
\]
Insert~\eqref{eq:micro-tail}, substitute $v=u/a$, and bound the
exponential integral by a constant times
$a^4(1+H/a)^4e^{-cH/a}$, decreasing $c$ if necessary. This proves the
second assertion of~\eqref{in:micro-moments}; the first follows in
the same way by integrating from zero, retaining the possible mass at
$X=0$.
\end{proof}

\subsection{Macroscopic attachments}

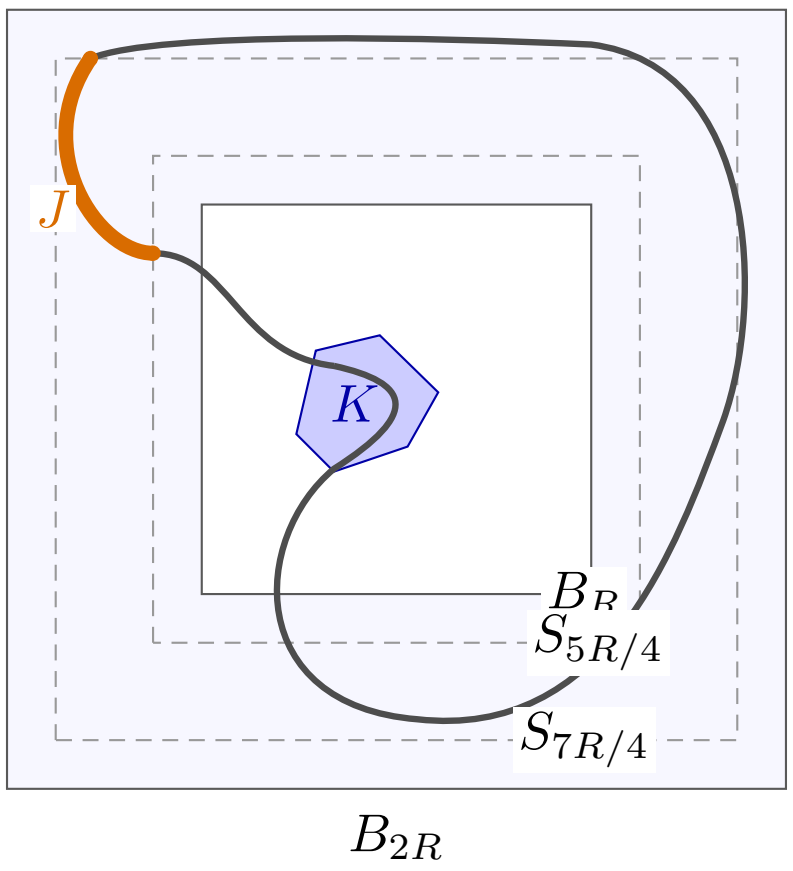
\begin{figure}[htbp]
\centering
\begin{minipage}[c]{0.47\textwidth}
\centering
\begin{tikzpicture}[x=0.92cm,y=0.92cm,font=\small,>=Stealth]
 \fill[blue!3] (-2.8,-2.8) rectangle (2.8,2.8);
 \fill[white] (-1.4,-1.4) rectangle (1.4,1.4);
 \draw[black!65] (-2.8,-2.8) rectangle (2.8,2.8);
 \draw[densely dashed,black!40] (-2.45,-2.45) rectangle (2.45,2.45);
 \draw[densely dashed,black!40] (-1.75,-1.75) rectangle (1.75,1.75);
 \draw[black!65] (-1.4,-1.4) rectangle (1.4,1.4);
 \filldraw[fill=blue!20,draw=blue!65!black]
 (-.72,-.25)--(-.45,-.52)--(.08,-.34)--(.30,.05)--(-.12,.46)--(-.58,.35)--cycle;
 \node[blue!65!black] at (-.30,-.03) {$K$};
 \draw[very thick,black!70]
 (-.45,.24) .. controls (-1.15,.32) and (-1.2,1.05) .. (-1.75,1.05)
 .. controls (-2.15,1.05) and (-2.65,1.8) .. (-2.20,2.45)
 .. controls (-1.7,2.65) and (.3,2.60) .. (1.40,2.55)
 .. controls (2.6,2.4) and (2.65,.7) .. (2.35,-.15)
 .. controls (2.0,-1.1) and (1.5,-2.45) .. (.15,-2.3)
 .. controls (-1.1,-2.2) and (-1.05,-1.0) .. (-.45,-.50)
 .. controls (.18,-.1) and (.10,.12) .. (-.45,.24);
 \draw[line width=2.8pt,orange!85!black]
 (-1.75,1.05) .. controls (-2.15,1.05) and (-2.65,1.8) .. (-2.20,2.45);
 \fill[orange!85!black] (-1.75,1.05) circle (.055);
 \fill[orange!85!black] (-2.20,2.45) circle (.055);
 \node[orange!85!black,fill=white,inner sep=1pt] at (-2.47,1.37) {$J$};
 \node[fill=white,inner sep=1pt] at (1.35,-1.40) {$B_R$};
 \node[fill=white,inner sep=1pt] at (1.45,-1.75) {$S_{5R/4}$};
 \node[fill=white,inner sep=1pt] at (1.35,-2.45) {$S_{7R/4}$};
 \node[below] at (0,-2.84) {$B_{2R}$};
\end{tikzpicture}
\end{minipage}\hfill
\begin{minipage}[c]{0.49\textwidth}
\small
The loop $\ell$ hits $K\subset B_R$ and remains in $B_{2R}$.

\medskip
The highlighted outer subset $J$ joins the two intermediate spheres and satisfies
\[
 \capc(J)\ge\kappa_0\frac{R^2}{\log R}.
\]
The loop therefore attaches to the inner cluster while carrying a controlled amount of capacity outside $B_R$.
\end{minipage}
\caption{The attachment family $T_R(K)$ from Gao's definition (7.15).
The drawing represents the nested boxes schematically in two dimensions;
the sets and their capacities in the proof are four-dimensional. The
particular radii provide fixed separation between the inner set, the
capacity-bearing outer subset, and the confining boundary.}
\label{fig:attachment}
\end{figure}

The two remaining external estimates concern the measure of loops hitting
a deterministic set. We state their spatial hypotheses explicitly.
Let $R\ge32$, put $a_4(R)=R^2/\log R$, and let $K\subset B_R$ be
finite. For $\kappa>0$, define $T_R^\kappa(K)$ to be the set of loops
$\ell$ such that
\begin{enumerate}
 \item $U(\ell)\subset B_{2R}$ and $U(\ell)\cap K\ne\varnothing$;
 \item there is a set $J\subset
 U(\ell)\cap(B_{2R}\setminus B_R)$ joining the spheres at radii
 $5R/4$ and $7R/4$, with $\capc(J)\ge\kappa a_4(R)$.
\end{enumerate}
This is the attachment family defined by Gao~\cite[(7.15)]{Gao};
Figure~\ref{fig:attachment} shows the geometry of the definition.
Boundary conventions in this definition are understood with lattice
rounding as above. Equivalently one may use Gao's annular notation
$J\subset\ell\cap A(0,R,2R)$; only the stated outer capacity and the
fact that the loop leaves $B_R$ will be used.
\begin{proposition}[Gao's attachment estimates]\label{in:attachment}
There exist fixed constants $\kappa_0,c_{\mathrm{att}},C_{\mathrm{att}}>0$
such that, for every $R\ge32$ and every finite $K\subset B_R$,
\begin{equation}\label{eq:attachment-input}
 c_{\mathrm{att}}R^{-2}\capc(K)
 \le\mu(T_R^{\kappa_0}(K))
 \le C_{\mathrm{att}}R^{-2}\capc(K).
\end{equation}
For each $\delta\in(0,1/2)$, there is a fixed
$C_{\mathrm{mac}}(\delta)<\infty$ such that
\begin{equation}\label{in:macroscopic}
 \mu\{\ell:U(\ell)\cap K\ne\varnothing,
                   \diam U(\ell)>\delta R\}
 \le C_{\mathrm{mac}}(\delta)R^{-2}\capc(K).
\end{equation}
\end{proposition}

These are Gao~\cite[Lemma~7.4, (7.15)--(7.17)]{Gao} in dimension four
and on the discrete graph. Since $T_R^{\kappa_0}(K)$ already consists
of loops contained in $B_{2R}$, its measure under $\mu$ agrees with
its measure under Gao's restricted loop measure. We henceforth abbreviate
$T_R(K)=T_R^{\kappa_0}(K)$. All capacities in these formulas have the
normalization fixed in Section~\ref{sec:introduction}.
The containment $K\subset B_R$ in~\eqref{in:macroscopic} is essential;
we will apply it only on events where it has been established.

The capacity dependence of these loop-measure estimates goes back to
Chang and Sapozhnikov~\cite[Lemma~2.7]{CS16}. Their upper bound controls
loops hitting $K\subset B_N$ and reaching $\partial B_M$, for
$M\ge\lambda N$ with fixed $\lambda>1$, by
$C\capc(K)M^{2-d}$. Their lower bound retains a comparable amount of
loop measure while requiring the \emph{whole loop} to have capacity at
least $cM^2/\log M$ in dimension four. The confinement refinement in
\cite[Section~6.2, (38)]{CS16} restricts such loops to a sufficiently
larger outer box. Gao's family specifies more spatial information:
the capacity is carried by a subset in a prescribed outer annulus,
and the loop stays in $B_{2R}$. Thus the stated attachment estimate is
Gao's localized version, rather than a literal restatement of
\cite[Lemma~2.7]{CS16}. The large-diameter upper estimate can also be
obtained from the latter lemma by covering $B_R$ with finitely many boxes
of size comparable to $\delta R$, with a constant depending on $\delta$.

\section{Capacity and geometrically killed random walk}\label{sec:potential}

The loop estimates below involve two different time scales: the duration of a
loop and the duration of the shorter loops used to connect it to another
loop.  Ordinary capacity does not keep track of either scale.  We therefore
first introduce capacity for a walk killed at an independent geometric time.
This section fixes its normalization and proves the comparison estimates
used later.  All constants in this section are universal for four-dimensional
simple random walk.

\subsection{Equilibrium measures and the two variational principles}

For a finite set $A\subset\Z^4$, define the entrance and return times by
\[
 H_A=\inf\{k\ge0:S_k\in A\},\qquad
 \widetilde H_A=\inf\{k\ge1:S_k\in A\}.
\]
The infimum of the empty set is $\infty$.  For $0<z\le1$, put
\begin{align}
 G_z(x,y)&=\sum_{k=0}^{\infty}z^k p_k(x,y),\notag\\
 h_A^z(x)&=\E_x[z^{H_A};H_A<\infty],\label{pt:potential-def}\\
 e_A^z(x)&=\one_{\{x\in A\}}
 \bigl(1-\E_x[z^{\widetilde H_A};\widetilde H_A<\infty]\bigr),
 \qquad
 \capc_z(A)=\sum_{x\in A}e_A^z(x).
 \label{pt:equilibrium}
\end{align}
At $z=1$ these are the usual Green function, equilibrium potential,
equilibrium measure, and capacity; we omit the superscript or subscript $1$.
For the empty set, the potential and capacity are zero.  The notation
$\capc_z$ always refers to killing in the whole lattice, and not to killing
on exiting a finite box.

For $z<1$, let $\xi_z$ be independent of $S$ with
$\Pp(\xi_z\ge k)=z^k$ for $k=0,1,\ldots$.  Thus $\xi_z$ is the number of
steps completed before killing.  The two probabilities in
\eqref{pt:potential-def}--\eqref{pt:equilibrium} can be written as
\[
 h_A^z(x)=\Pp_x(H_A\le\xi_z),\qquad
 e_A^z(x)=\one_{\{x\in A\}}\Pp_x(\widetilde H_A>\xi_z).
\]
Decomposing according to the last visit to $A$ before killing gives
\begin{equation}\label{pt:last-exit}
 h_A^z(x)=\sum_{y\in A}G_z(x,y)e_A^z(y).
\end{equation}
Indeed, the contribution of a last visit at time $k$ and position $y$ is
$z^k p_k(x,y)e_A^z(y)$, by the memoryless property of $\xi_z$ and the Markov
property of $S$.  Summing over $k$ and $y$ proves the formula.  At $z=1$,
transience ensures that a finite set has a last visit whenever it is hit,
so the same decomposition applies.  This is the usual last-exit identity;
see Lawler and Limic~\cite[Proposition~4.6.4]{LL} for its spatially killed
and transient versions.  In particular, $h_A^z=1$ on $A$ and
$0\le h_A^z\le1$ everywhere.

We record explicitly the variational formulas that will be used.  For
finitely supported real functions $f,g$ on $\Z^4$, write
\begin{align}
 \mathcal E(f,g)
 &=\frac12\sum_{x,y\in\Z^4}P(x,y)
       (f(x)-f(y))(g(x)-g(y)),\notag\\
 \mathcal E_z(f,g)
 &=z\mathcal E(f,g)+(1-z)\sum_x f(x)g(x).
 \label{pt:dirichlet-form}
\end{align}
Set $\mathcal E_z(f)=\mathcal E_z(f,f)$ and
$\mathcal E(f)=\mathcal E(f,f)$.  Since $P$ is symmetric with respect to
counting measure, $\mathcal E_z(f,g)=\langle f,(I-zP)g\rangle$ whenever the
expressions are defined.  For every nonempty finite $A$,
\begin{align}
 \capc_z(A)
 &=\inf\{\mathcal E_z(f): f\text{ is finitely supported and }f\ge1
                       \text{ on }A\},\label{pt:dirichlet-principle}\\
 \frac1{\capc_z(A)}
 &=\min_{\nu\in\mathcal P(A)}
       \sum_{x,y\in A}\nu(x)G_z(x,y)\nu(y),
 \label{pt:green-principle}
\end{align}
where $\mathcal P(A)$ is the set of probability measures supported on $A$.
The minimizing probability measure in \eqref{pt:green-principle} is
$e_A^z/\capc_z(A)$.

Here is a justification that also explains the boundary condition implicit
in \eqref{pt:dirichlet-principle}.  For $z<1$, the potential $h_A^z=G_z e_A^z$
belongs to $\ell^2(\Z^4)$, and
\[
 (I-zP)h_A^z=e_A^z,\qquad
 \mathcal E_z(h_A^z)=\sum_{x\in A}h_A^z(x)e_A^z(x)=\capc_z(A).
\]
If $f\ge1$ on $A$, expanding the quadratic form around $h_A^z$ gives
\[
 \mathcal E_z(f)=\capc_z(A)
 +2\sum_{x\in A}(f(x)-1)e_A^z(x)
 +\mathcal E_z(f-h_A^z)\ge\capc_z(A).
\]
Approximating $h_A^z$ by finite-support truncations gives the reverse
inequality in the infimum.  One can choose truncations that agree with
$h_A^z=1$ on $A$.  At $z=1$, first kill the walk on exiting an increasing
sequence of finite boxes containing $A$.  The corresponding equilibrium
potentials minimize the finite-box Dirichlet energy.  Their equilibrium
measures decrease to $e_A$, because the events of exiting a box before
returning to $A$ decrease to the event of never returning.  Every finitely
supported competitor is supported in one such box.  Passing to the limit
therefore proves \eqref{pt:dirichlet-principle} also for $z=1$.  The
full-space potential has finite Dirichlet energy and satisfies
\begin{equation}\label{pt:potential-energy}
 \mathcal E(h_A)=\capc(A),
\end{equation}
as can also be seen from the exhaustion.  More explicitly, let $h_j$ be the
potential killed on exiting the $j$th box, extended by zero outside it,
and let $c_j=\mathcal E(h_j)$ be its capacity.  If $j<k$, then
$h_j-h_k$ is supported in the $k$th box and vanishes on $A$.  The
Euler--Lagrange equation for the minimizer $h_k$ gives
$\mathcal E(h_k,h_j-h_k)=0$, and consequently
\[
 \mathcal E(h_j-h_k)=c_j-c_k.
\]
As $c_j\downarrow\capc(A)$, the edge gradients of $h_j$ are Cauchy in
their weighted $\ell^2$ space.  Since $h_j(x)\uparrow h_A(x)$ at every
vertex, their limit is the gradient of $h_A$.  Passing to the limit in
$\mathcal E(h_j)=c_j$ proves \eqref{pt:potential-energy}.  No square
integrability of $h_A$ is asserted.

For completeness, the probability-measure principle can be checked
directly from \eqref{pt:last-exit}.  Let
$\nu_A=e_A^z/\capc_z(A)$.  For any $\nu\in\mathcal P(A)$, set
$\rho=\nu-\nu_A$; then $\sum_x\rho(x)=0$.  Since
$G_z\nu_A=1/\capc_z(A)$ on $A$,
\[
 \langle\nu,G_z\nu\rangle
 =\frac1{\capc_z(A)}+\langle\rho,G_z\rho\rangle
 \ge\frac1{\capc_z(A)}.
\]
The last quadratic form is nonnegative: for $z<1$, $G_z$ is the inverse of
the positive operator $I-zP$; for $z=1$, take $z\uparrow1$ in the quadratic
form of the finitely supported function $\rho$.  This proves
\eqref{pt:green-principle}.  The standard theory in
\cite[Section~6.5]{LL} uses the same normalization of ordinary capacity.

\subsection{Monotonicity and comparison with ordinary capacity}

The comparison below says that killing changes capacity by at most a fixed
factor when the set lies within the walk's diffusive scale before killing.
The proof tests the killed Dirichlet form with the ordinary equilibrium
potential multiplied by a spatial cutoff. This makes the killing term
finite while keeping the boundary-energy cost proportional to ordinary
capacity.

\begin{lemma}\label{pt:massive-comparison}
For $0<z\le1$, $\capc_z$ is monotone under inclusion of finite sets,
subadditive, and at most cardinality.  If $0<z_1\le z_2\le1$, then
$\capc_{z_2}(A)\le\capc_{z_1}(A)$.  There is a constant
$C_{\mathrm{mc}}<\infty$ such that, whenever $A\subset B(x_0,r)$ is finite,
$r\ge0$, and $0<z\le1$,
\begin{equation}\label{pt:massive-comparison-eq}
 \capc_z(A)\le C_{\mathrm{mc}}
     \bigl(1+(1-z)(r+1)^2\bigr)\capc(A).
\end{equation}
Here $B(x_0,r)=x_0+B_r$ is a box in the $\ell^\infty$ norm.
\end{lemma}

\begin{proof}
Inclusion of sets reduces the admissible class in
\eqref{pt:dirichlet-principle}, giving monotonicity.  The test function
$\one_A$ gives
\[
 \capc_z(A)\le\mathcal E_z(\one_A)
 =|A|-z\sum_{x,y\in A}P(x,y)\le|A|.
\]
For subadditivity, competitors can first be truncated to take values in
$[0,1]$, which does not increase their energy.  For two such competitors
$f$ and $g$, the pointwise max--min inequality for each edge gives
\[
 \mathcal E_z(f\vee g)+\mathcal E_z(f\wedge g)
 \le\mathcal E_z(f)+\mathcal E_z(g).
\]
The killing term satisfies equality in this inequality because
$(f\vee g)^2+(f\wedge g)^2=f^2+g^2$ at each vertex.  Since $f\vee g$ is
admissible for the union, taking infima proves
$\capc_z(A\cup A')\le\capc_z(A)+\capc_z(A')$.

To prove monotonicity in $z$, use the Green, rather than the Dirichlet,
variational formula.  The entries of $G_z$ are nondecreasing in $z$.
Consequently, for every nonnegative probability measure $\nu$ on $A$,
\[
 \sum_{x,y}\nu(x)G_{z_1}(x,y)\nu(y)
 \le\sum_{x,y}\nu(x)G_{z_2}(x,y)\nu(y).
\]
Taking minima and reciprocals in \eqref{pt:green-principle} proves the
claimed direction.  This argument requires only entrywise monotonicity
against nonnegative measures; it does not assert monotonicity of $G_z$ in
the ordering of quadratic forms on signed functions.

It remains to prove \eqref{pt:massive-comparison-eq}.  The assertion is
trivial if $A$ is empty, so assume otherwise.  Write $h=h_A$.

The natural competitor in \eqref{pt:dirichlet-principle} is $h$ itself:
it equals one on $A$, and its Dirichlet energy is exactly $\capc(A)$ by
\eqref{pt:potential-energy}.  For $z<1$ it is not admissible, because the
killing term $(1-z)\sum_x h(x)^2$ need not be finite.  In four dimensions
$h$ decays only like the Green function, that is like $|x|^{-2}$, and is
not square summable.  We therefore multiply $h$ by a radial cutoff $\chi$
that equals one on the box containing $A$ and decreases linearly to zero
across a shell whose width is comparable to the box radius; see
Figure~\ref{fig:cutoff}.  The product $\chi h$ is still one on $A$. Since $0\le h\le1$,
its squared $\ell^2$ norm is bounded by the $\ell^1$ norm of $h$ on a
box of radius $O(r+1)$, hence by $C(r+1)^2\capc(A)$ using the Green
bound. The gradient of $\chi$ is bounded by $(r+1)^{-1}$, so its extra
Dirichlet-energy cost is at most $C\capc(A)$. The killing term has the
separate bound $C(1-z)(r+1)^2\capc(A)$; these two estimates explain
the factor $1+(1-z)(r+1)^2$ in
\eqref{pt:massive-comparison-eq}.  Concretely, choose
\[
 \chi(y)=\left(1-
      \frac{(\|y-x_0\|_\infty-r)_+}{r+1}\right)_+.
\]
Thus $\chi=1$ on $B(x_0,r)$, $0\le\chi\le1$, and $\chi$ vanishes outside
$B(x_0,2r+1)$.  If $P(u,v)>0$, then
$|\chi(u)-\chi(v)|\le(r+1)^{-1}$.  In particular, $f=\chi h$ is finitely
supported and equals one on $A$.

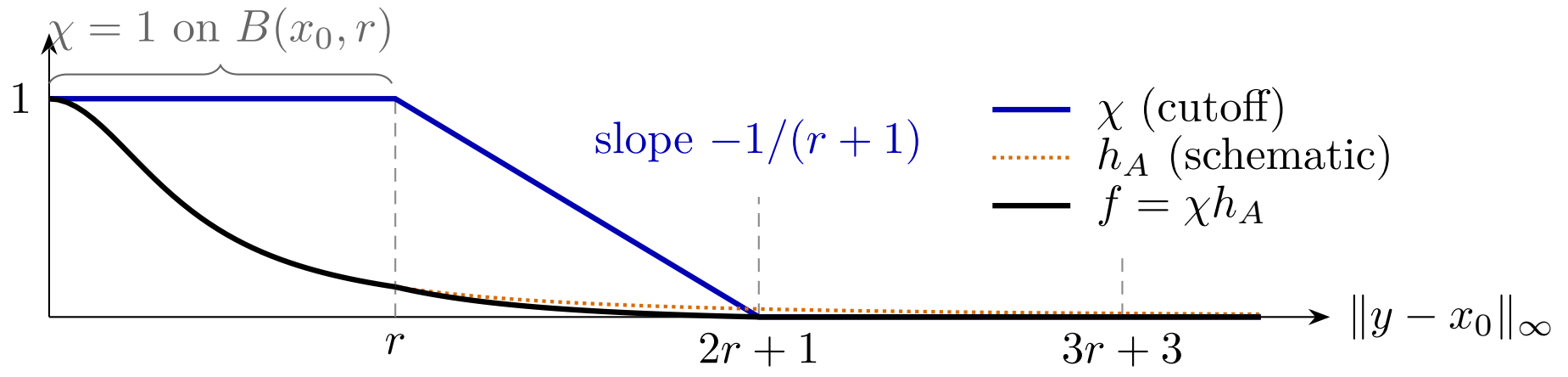
\begin{figure}[htbp]
\centering
\begin{tikzpicture}[x=1.35cm,y=1.7cm,font=\small,>=Stealth]
 \draw[->] (0,0) -- (7.4,0) node[right] {$\|y-x_0\|_\infty$};
 \draw[->] (0,0) -- (0,1.3);
 \node[left] at (0,1) {$1$};
 \draw[densely dashed,black!45] (2,0) -- (2,1);
 \draw[densely dashed,black!45] (4.1,0) -- (4.1,0.55);
 \draw[densely dashed,black!45] (6.2,0) -- (6.2,0.3);
 \node[below] at (2,0) {$r$};
 \node[below] at (4.1,0) {$2r+1$};
 \node[below] at (6.2,0) {$3r+3$};
 \draw[very thick,blue!70!black] (0,1) -- (2,1) -- (4.1,0) -- (7,0);
 \draw[thick,densely dotted,orange!85!black,domain=0:7,samples=80]
   plot (\x,{1/(1+(\x/0.8)^2)});
 \draw[very thick,black,domain=0:2,samples=40]
   plot (\x,{1/(1+(\x/0.8)^2)});
 \draw[very thick,black,domain=2:4.1,samples=40]
   plot (\x,{(1/(1+(\x/0.8)^2))*(1-(\x-2)/2.1)});
 \draw[very thick,black] (4.1,0) -- (7,0);
 \draw[decorate,decoration={brace,amplitude=4pt},black!60]
   (0,1.07) -- (1.98,1.07) node[midway,above=3pt] {$\chi=1$ on $B(x_0,r)$};
 \node[blue!70!black,anchor=south west] at (3.05,0.62) {slope $-1/(r+1)$};
 \begin{scope}[shift={(5.45,0.95)}]
  \draw[very thick,blue!70!black] (0,0) -- (0.45,0);
  \node[anchor=west] at (0.5,0) {$\chi$ (cutoff)};
  \draw[thick,densely dotted,orange!85!black] (0,-0.22) -- (0.45,-0.22);
  \node[anchor=west] at (0.5,-0.22) {$h_A$ (schematic)};
  \draw[very thick,black] (0,-0.44) -- (0.45,-0.44);
  \node[anchor=west] at (0.5,-0.44) {$f=\chi h_A$};
 \end{scope}
\end{tikzpicture}
\caption{The radial cutoff used in the proof of
Lemma~\ref{pt:massive-comparison}, drawn as a function of the
$\ell^\infty$ distance from the center $x_0$ of the box containing $A$.
The dotted curve represents only the decay of $h_A$;
the actual potential need not be radial. It equals one on $A$ but decays
only like the Green function, so in four dimensions it is not square summable and
cannot be inserted into the killed Dirichlet form directly.  The cutoff
$\chi$ (blue) equals one on $B(x_0,r)$, decreases linearly to zero on
the shell between radii $r$ and $2r+1$, and has gradient of order $1/r$
there.  The admissible test function is the product $f=\chi h_A$ (solid),
whose support lies in $D=B(x_0,3r+3)$ together with all edges on which
$\chi$ varies.}
\label{fig:cutoff}
\end{figure}

We first bound its squared $\ell^2$ norm.  Let $D=B(x_0,3r+3)$.  Since
$0\le h\le1$, last exit and the Green bound in the random-walk estimates
give
\begin{align}
 \sum_{y\in D}h(y)^2
 &\le\sum_{y\in D}h(y)
 =\sum_{a\in A}e_A(a)\sum_{y\in D}G(y,a)\notag\\
 &\le C\capc(A)
       \sum_{k=0}^{\lceil4r+4\rceil}(k+1)^3(k+1)^{-2}
 \le C(r+1)^2\capc(A).
 \label{pt:localized-potential}
\end{align}
The middle inequality uses $a\in B(x_0,r)$ and the fact that an
$\ell^\infty$ shell of radius $k$ in four dimensions contains at most
$C(k+1)^3$ vertices.  Since the support of $f$ is contained in $D$,
\eqref{pt:localized-potential} also bounds $\|f\|_2^2$.

The required product estimate for the Dirichlet energy follows from the
elementary identity
\[
 f(u)-f(v)=\chi(u)(h(u)-h(v))
                 +h(v)(\chi(u)-\chi(v))
\]
and $(a+b)^2\le2a^2+2b^2$.  Substituting this into
\eqref{pt:dirichlet-form} yields
\begin{align*}
 \mathcal E(f)
 &\le\sum_{u,v}P(u,v)\chi(u)^2(h(u)-h(v))^2
       +\sum_{u,v}P(u,v)h(v)^2(\chi(u)-\chi(v))^2\\
 &\le 2\mathcal E(h)
       +\frac{C}{(r+1)^2}\sum_{v\in D}h(v)^2
 \le C\capc(A).
\end{align*}
In the second line, a nonzero difference of $\chi$ across an edge forces
both endpoints into $D$.  The last line uses
\eqref{pt:potential-energy} and \eqref{pt:localized-potential}.
Finally, insert $f$ into \eqref{pt:dirichlet-principle}:
\[
 \capc_z(A)\le z\mathcal E(f)+(1-z)\|f\|_2^2
 \le C\bigl(1+(1-z)(r+1)^2\bigr)\capc(A).
\]
This proves the comparison.  The compact cutoff is indispensable in this
argument: in dimension four, $h_A$ itself need not belong to $\ell^2$.
\end{proof}

We will also use a spatially killed last-exit formula.  For a finite set
$D$, let $\tau_D=\inf\{k\ge0:S_k\notin D\}$ and define
\[
 G_D(x,y)=\E_x\left[\sum_{k=0}^{\tau_D-1}\one_{\{S_k=y\}}\right],
 \qquad
 e_{A,D}(y)=\one_{\{y\in A\}}\Pp_y(\widetilde H_A\ge\tau_D)
 \quad(A\subset D).
\]
The same last-visit decomposition gives
\begin{equation}\label{pt:killed-last-exit}
 \Pp_x(H_A<\tau_D)=\sum_{y\in A}G_D(x,y)e_{A,D}(y),
 \qquad e_{A,D}(y)\ge e_A(y).
\end{equation}
For example, applying the interior killed-Green estimate from the
random-walk inputs to $D=B_{4r}$ gives a constant $c_{\mathrm{hit}}>0$ such
that
\begin{equation}\label{pt:interior-hit}
 \inf_{x\in B_{r/4}}\Pp_x(H_A<\tau_{B_{4r}})
 \ge c_{\mathrm{hit}}r^{-2}\capc(A),
 \qquad A\subset B_{2r},\quad r\ge1.
\end{equation}
Indeed, each summand in \eqref{pt:killed-last-exit} is bounded below by
$c_{\mathrm{hit}}r^{-2}e_A(y)$.  The interior Green estimate is a standard
consequence of the ball Green estimates and Harnack inequality; see
\cite[Proposition~6.3.5 and Theorem~6.3.9]{LL}.  The displayed formulation
specifies exactly the geometry needed in the annular argument below.

\subsection{The mass of loops in one duration band}

The next lemma converts the measure of loops hitting a set into its
massive capacity.  It applies to every finite set, without a restriction on
its diameter.  This is why it can later be applied to the random set
obtained by attaching short loops to a deterministic parent loop.

\begin{lemma}[Duration-band estimate]\label{pt:band}
There is a constant $C_{\mathrm{band}}<\infty$ such that, for every finite
$A\subset\Z^4$ and every real $M\ge2$,
\begin{equation}\label{pt:band-eq}
 \int \one_{\{M\le N(\ell)<2M\}}
       \one_{\{U(\ell)\cap A\ne\varnothing\}}
       \frac{N(\ell)}{\log(N(\ell)+2)}\,\mu(d\ell)
 \le\frac{C_{\mathrm{band}}}{\log(M+2)}
       \capc_{e^{-1/M}}(A).
\end{equation}
\end{lemma}

\paragraph{Proof idea.}
A loop of duration comparable to $M$ which hits $A$ can be counted by its
first hit during half of its lifetime, followed by a return segment. The
heat-kernel bound controls the return segment. Summing the first-hit
probabilities over all starting vertices is then a heat-content problem,
whose geometrically killed version equals $\capc_{e^{-1/M}}(A)/(1-e^{-1/M})$.
The remaining duration sum cancels the extra factor of order $M$.

\begin{proof}
For an integer $n\ge2$, let $\mu_n$ denote the restriction of $\mu$ to
loops of duration $n$.  The rooted formula \eqref{eq:model-rooted} gives
\begin{equation}\label{pt:rooted-hit}
 \mu_n(U(\ell)\cap A\ne\varnothing)
 =\frac1n\sum_{x\in\Z^4}
       \Pp_x(S_n=x,\ H_A\le n).
\end{equation}
This includes the factor $1/n$ from the definition of the unrooted loop
measure; no additional factor is inserted for periodic loops.

For a fixed root $x$, reversal of the closed path
$(S_0,S_1,\ldots,S_n)$ preserves its probability, since the walk is
symmetric.  A visit to $A$ in the second half of a closed path becomes a
visit in its first half after reversal.  Thus, with
$q_n=\lfloor n/2\rfloor$,
\[
 \Pp_x(S_n=x,\ H_A\le n)
 \le 2\Pp_x(S_n=x,\ H_A\le q_n).
\]
The strong Markov property at $H_A$ then gives the exact formula
\[
 \Pp_x(S_n=x,\ H_A\le q_n)
 =\E_x\!\left[
    \one_{\{H_A\le q_n\}}p_{n-H_A}(S_{H_A},x)\right].
\]
On the event in this expectation, $n-H_A\ge n/2$.  The heat-kernel upper
bound therefore implies, for a constant independent of $x,n,A$,
\begin{equation}\label{pt:loop-heat-content}
 \mu_n(U(\ell)\cap A\ne\varnothing)
 \le Cn^{-3}\sum_{x\in\Z^4}\Pp_x(H_A\le q_n).
\end{equation}

Fix $z=e^{-1/M}$ for the entire duration band.  If $M\le n<2M$ and
$H_A\le q_n$, then $H_A<M$ and hence $z^{H_A}\ge e^{-1}$.  Consequently,
\[
 \Pp_x(H_A\le q_n)\le e\,h_A^z(x).
\]
The total mass of this potential is explicit.  By
\eqref{pt:last-exit}, symmetry of $P$, and Tonelli's theorem,
\begin{align}
 \sum_{x\in\Z^4}h_A^z(x)
 &=\sum_{y\in A}e_A^z(y)
       \sum_{k=0}^{\infty}z^k\sum_{x\in\Z^4}p_k(x,y)\notag\\
 &=\frac{\capc_z(A)}{1-z}.
 \label{pt:heat-content}
\end{align}
The equality $\sum_xp_k(x,y)=1$ follows from symmetry and the fact that
the rows of $P^k$ sum to one.  Since $M\ge2$,
$1-e^{-1/M}\ge c/M$ with a universal $c>0$.  Inserting
\eqref{pt:heat-content} into \eqref{pt:loop-heat-content} yields
\[
 \mu_n(U(\ell)\cap A\ne\varnothing)
 \le CMn^{-3}\capc_{e^{-1/M}}(A),\qquad M\le n<2M.
\]
Multiplying by $n/\log(n+2)$ and summing over integers in the band gives
\begin{align*}
 \sum_{M\le n<2M}\frac{n}{\log(n+2)}
       \mu_n(U(\ell)\cap A\ne\varnothing)
 &\le \frac{CM\capc_{e^{-1/M}}(A)}{\log(M+2)}
       \sum_{M\le n<2M}n^{-2}\\
 &\le\frac{C\capc_{e^{-1/M}}(A)}{\log(M+2)}.
\end{align*}
Here the last sum is at most $C/M$.  Odd values of $n$ contribute zero,
because nearest-neighbor walk on $\Z^4$ is bipartite.  This proves
\eqref{pt:band-eq} for arbitrary real $M\ge2$.
\end{proof}

\subsection{An elementary lower-capacity estimate for a walk range}

The strong law of Asselah, Schapira and Sousi~\cite[Theorem~1.1]{ASS}
already gives the sharper asymptotic statement
\begin{equation}\label{pt:capacity-slln}
 \frac{\log n}{n}\capc(U_n)\longrightarrow\frac{\pi^2}{8}
 \qquad\Pp_0\text{-almost surely}.
\end{equation}
In particular, it could supply the lower-capacity input needed below:
for any fixed confidence level, the capacity exceeds a fixed multiple of
$n/\log n$ once $n$ is sufficiently large.

We give an elementary estimate because the annular construction needs only
this fixed-confidence consequence, and the elementary proof also gives a
bound for every duration. Its threshold is allowed to decrease with the
requested error probability. After choosing that error probability once,
the strong Markov property amplifies the estimate over many annuli. Thus
using the weaker statement loses no conclusion needed for the argument;
the strong law could replace it, but neither its sharp constant nor its
almost-sure conclusion is used. The following proof uses occupation energy
and Markov's inequality.

\begin{lemma}\label{pt:lower-capacity}
There is a universal constant $c_{\mathrm{en}}>0$ such that, for
four-dimensional simple random walk, every integer $t\ge1$, and every
$\zeta\in(0,1)$,
\begin{equation}\label{pt:lower-capacity-eq}
 \Pp_0\left(\capc(S[0,t])<
       c_{\mathrm{en}}\zeta\frac{t+1}{\log(t+2)}\right)\le\zeta.
\end{equation}
\end{lemma}

\begin{proof}
Let $V_t=S[0,t]$ and consider the empirical occupation probability measure
\[
 \nu_t(x)=\frac1{t+1}\sum_{s=0}^{t}\one_{\{S_s=x\}},\qquad x\in V_t.
\]
Repeated visits are counted in $\nu_t$, so it has total mass exactly one.
Its Green energy is
\[
 \sum_{x,y\in V_t}\nu_t(x)G(x,y)\nu_t(y)
 =\frac{I_t}{(t+1)^2},\qquad
 I_t=\sum_{s,u=0}^{t}G(S_s,S_u).
\]
The variational principle \eqref{pt:green-principle} therefore implies the
deterministic inequality
\begin{equation}\label{pt:empirical-bound}
 \capc(V_t)\ge\frac{(t+1)^2}{I_t}.
\end{equation}

For $h=|u-s|$, stationary independent increments, symmetry, and the
semigroup identity imply
\begin{align*}
 \E_0[G(S_s,S_u)]
 &=\sum_{x\in\Z^4}p_h(0,x)G(0,x)\\
 &=\sum_{k=0}^{\infty}\sum_{x\in\Z^4}p_h(0,x)p_k(x,0)
 =\sum_{k=0}^{\infty}p_{h+k}(0,0)
 \le\frac{C}{h+1}.
\end{align*}
The final inequality follows by summing the heat-kernel bound
$p_j(0,0)\le C(j+1)^{-2}$.  There are $t+1$ ordered pairs with $h=0$ and
$2(t+1-h)$ ordered pairs with $h\ge1$.  Hence there is a fixed
$C_{\mathrm{en}}<\infty$ such that
\[
 \E_0[I_t]
 \le C\left((t+1)+2\sum_{h=1}^{t}\frac{t+1-h}{h+1}\right)
 \le C_{\mathrm{en}}(t+1)\log(t+2),\qquad t\ge1.
\]
Choose $c_{\mathrm{en}}=C_{\mathrm{en}}^{-1}$.  By
\eqref{pt:empirical-bound}, the event in \eqref{pt:lower-capacity-eq}
implies
$I_t>C_{\mathrm{en}}(t+1)\log(t+2)/\zeta$.  Markov's inequality bounds
the probability of the latter event by $\zeta$, as required.
\end{proof}

\section{Capacity of neighborhoods of random-walk ranges}
\label{sec:concentration}\label{conc:section}

The loop estimate used later must remain valid after a loop range has been
thickened by a small power of its duration. An ordinary mean-capacity bound
does not immediately give such an estimate: a deterministic thickening can
add many vertices. We first prove that, with arbitrarily high polynomial
probability, the original range has a uniformly positive hitting potential
throughout this neighborhood. A potential-theoretic comparison then controls
the capacity of the thickened set. We subsequently apply this estimate to
independent blocks of walk increments to obtain one good-shape event which
works at every killing scale.

Throughout this section $U_n=\{S_0,\ldots,S_n\}$ denotes the range of a
walk started at zero. We distinguish the free-walk law $\Pp_0$, with
expectation $\E_0$, from the bridge law $\Pp_n^{\mathrm{br}}$, with
expectation $\E_n^{\mathrm{br}}$. The latter is used only for even $n\ge2$.
For a finite set $A$, recall that $h_A(y)=\Pp_y(H_A<\infty)$ is its
ordinary full-space hitting potential.

\subsection{A positive potential throughout the neighborhood}

The next proposition is the main range estimate. Its mechanism is
multiscale: a walk typically supplies a substantial-capacity piece in each
of many well-separated annuli. From inside the $j$th annulus, a new test
walk has probability of order $1/j$ to hit such a piece. The annuli used
here have indices of order $\log n$, and their number is also of order
$\log n$. The accumulated hitting probability is therefore bounded away
from zero. Failure probabilities for the construction can be made smaller
than any prescribed power of $n$ by making the acceptable capacity
threshold small enough. This last freedom changes the lower bound on the
potential, but not its uniformity in $n$.

\begin{proposition}[Potential on a polynomial neighborhood]
\label{conc:sausage}
Let $0<\varepsilon<1/2$ and $Q>0$. There are constants
$a=a(\varepsilon,Q)\in(0,1)$ and
$C=C(\varepsilon,Q)<\infty$ such that
\begin{equation}
\label{conc:halo-potential}
 \begin{aligned}
 \Pp_0\left(\inf_{y\in U_n+B_{n^\varepsilon}}h_{U_n}(y)<a\right)
    &\le Cn^{-Q}, && n\ge1,\\
 \Pp_n^{\mathrm{br}}\left(
      \inf_{y\in U_n+B_{n^\varepsilon}}h_{U_n}(y)<a\right)
    &\le Cn^{-Q}, && n\ge2\text{ even}.
 \end{aligned}
\end{equation}
Consequently, the same constants satisfy
\begin{equation}
\label{conc:halo-capacity}
 \begin{aligned}
 \Pp_0\left(\capc(U_n+B_{n^\varepsilon})
                   >a^{-1}\capc(U_n)\right)
    &\le Cn^{-Q}, && n\ge1,\\
 \Pp_n^{\mathrm{br}}\left(\capc(U_n+B_{n^\varepsilon})
                   >a^{-1}\capc(U_n)\right)
    &\le Cn^{-Q}, && n\ge2\text{ even}.
 \end{aligned}
\end{equation}
There is also a constant $C_\varepsilon<\infty$ such that
\begin{equation}
\label{conc:halo-mean}
 \begin{aligned}
 \E_0\left[\capc(U_n+B_{n^\varepsilon})\right]
    &\le C_\varepsilon\frac{n}{\log(n+2)}, && n\ge1,\\
 \E_n^{\mathrm{br}}\left[\capc(U_n+B_{n^\varepsilon})\right]
    &\le C_\varepsilon\frac{n}{\log(n+2)},
       &&n\ge2\text{ even}.
 \end{aligned}
\end{equation}
\end{proposition}

We prove the proposition using two lemmas. The first constructs many
suitable pieces from the first half of the walk; the second is a
deterministic estimate on the potential of their union. Keeping the first
lemma measurable from the first half is what makes the later bridge
comparison possible.

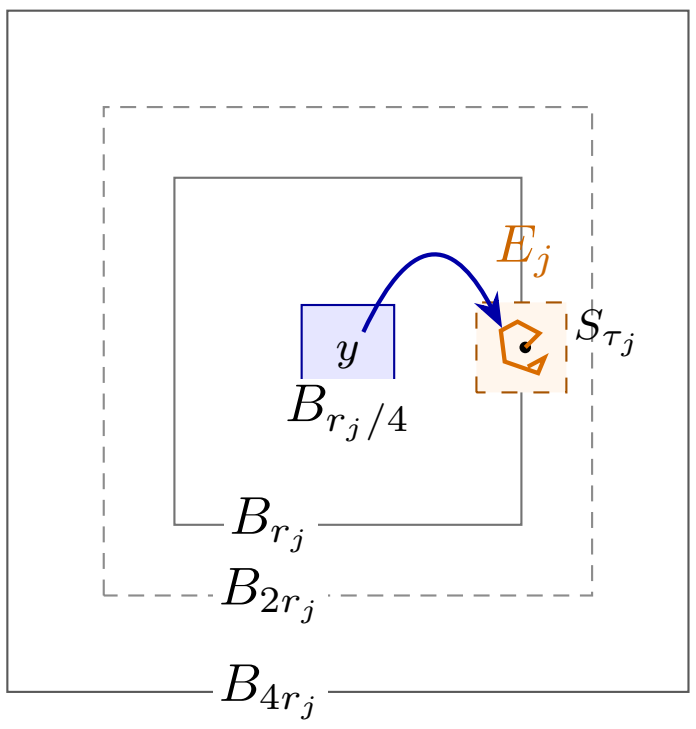
\begin{figure}[htbp]
\centering
\begin{minipage}[c]{0.43\textwidth}
\centering
\begin{tikzpicture}[x=0.85cm,y=0.85cm,font=\small,>=Stealth]
 \begin{scope}
  \draw[black!65] (-2.65,-2.65) rectangle (2.65,2.65);
  \draw[densely dashed,black!45] (-1.9,-1.9) rectangle (1.9,1.9);
  \draw[black!55] (-1.35,-1.35) rectangle (1.35,1.35);
  \filldraw[fill=blue!10,draw=blue!60!black] (-.36,-.36) rectangle (.36,.36);
  \node[font=\scriptsize] at (0,-.03) {$y$};
  \filldraw[fill=orange!7,draw=orange!65!black,dashed]
       (1.00,-.32) rectangle (1.70,.38);
  \fill[black] (1.38,.03) circle (.045);
  \draw[thick,orange!85!black]
     (1.38,.03)--(1.49,.14)--(1.32,.23)--(1.19,.16)--(1.22,-.08)
       --(1.48,-.17)--(1.53,-.05)--(1.40,-.12);
  \draw[blue!65!black,thick,->] (.12,.15)
       .. controls (.52,1.0) and (.85,.85) .. (1.20,.18);
  \node[anchor=south,orange!80!black] at (1.38,.42) {$E_j$};
  \node[anchor=south west,font=\scriptsize] at (1.60,-.20) {$S_{\tau_j}$};
  \node[fill=white,inner sep=1pt] at (0,-.48) {$B_{r_j/4}$};
  \node[fill=white,inner sep=1pt] at (-.60,-1.35) {$B_{r_j}$};
  \node[fill=white,inner sep=1pt] at (-.60,-1.90) {$B_{2r_j}$};
  \node[fill=white,inner sep=1pt] at (-.60,-2.65) {$B_{4r_j}$};
 \end{scope}
\end{tikzpicture}
\end{minipage}\hfill
\begin{minipage}[c]{0.55\textwidth}
\small
\textbf{One spatial scale.} The parent walk supplies $E_j$ inside
$B(S_{\tau_j},r_j/4)$. A test walk starting in $B_{r_j/4}$ hits it before
leaving $B_{4r_j}$ with probability at least $c\,b/j$.

\medskip
\textbf{Separation of scales.} Since $r_{j+1}=32r_j$, the next starting
radius is $r_{j+1}/4=8r_j$. An exit from $B_{4r_j}$ is still inside the
next admissible starting box.

\medskip
\textbf{Order of decisions.} We have
\[
 \tau_j\le\sigma_j<\tau_{j+1},\qquad
 \sigma_j=\tau_j+\min\{\rho_j,m_j\}.
\]
The $j$th good-piece decision is known before the next piece is started.
\end{minipage}
\caption{The annular construction and the independent test walk.
This schematic compresses the scale ratios and omits lattice overshoots.
The orange set belongs to the parent range; the blue path represents the
test walk used to estimate its hitting potential. Only pieces completed
by time $\lfloor n/2\rfloor$ are retained when passing to the bridge law.}
\label{fig:annular-scales}
\end{figure}

Fix $\varepsilon\in(0,1/2)$. For integers $j\ge1$, put $r_j=32^j$.
For all sufficiently large $n$, define
\[
 j_- = \left\lceil\log_{32}(8n^\varepsilon)\right\rceil,
 \qquad
 j_+ = \left\lfloor\log_{32}\left(\frac{\sqrt n}{\log n}\right)
                                  \right\rfloor-1,
 \qquad \mathcal J_n=\{j_-,\ldots,j_+\}.
\]
Write $k=|\mathcal J_n|$. There are constants
$\gamma_\varepsilon>0$, $\Gamma<\infty$, and a threshold
$n_\varepsilon$ such that
\begin{equation}
\label{conc:scale-count}
 j_-\ge1,\qquad
 k\ge\gamma_\varepsilon\log n,\qquad
 j_+\le\Gamma\log n,\qquad
 r_{j_++1}\le\frac{\sqrt n}{\log n},
 \qquad n\ge n_\varepsilon.
\end{equation}
Indeed, the leading terms in $j_-$ and $j_+$ are respectively
$\varepsilon\log n/\log32$ and $(\log n)/(2\log32)$;
the $\log\log n$ correction in $j_+$ is smaller. The radii and the
piece construction are illustrated in Figure~\ref{fig:annular-scales}.

\begin{lemma}[Good pieces determined by the first half of the walk]
\label{conc:annular-pieces}
For every $q>0$, there are constants $b=b(\varepsilon,q)>0$,
$C=C(\varepsilon,q)<\infty$, and $n_0=n_0(\varepsilon,q)$ such that,
for each $n\ge n_0$, there is an event
$F_n\in\sigma(S_0,\ldots,S_{\lfloor n/2\rfloor})$ satisfying
\begin{equation}
\label{conc:first-half-event}
 \Pp_0(F_n^c)\le Cn^{-q}.
\end{equation}
On $F_n$, at least $k/2$ indices $j\in\mathcal J_n$ have an associated
set $E_j$ with
\begin{equation}
\label{conc:piece-output}
 E_j\subset U_{\lfloor n/2\rfloor}
              \cap(B_{2r_j}\setminus B_{r_j/2}),
 \qquad
 \capc(E_j)\ge b\frac{r_j^2}{j}.
\end{equation}
\end{lemma}

\paragraph{Proof idea.}
The walk is observed at its successive exit times from the boxes
$B_{r_j}$. Immediately after each exit we follow it for a time $m_j$
which is a small multiple of $r_j^2$: small enough that the segment
typically stays within distance $r_j/4$ of the exit point, hence inside
the annulus $B_{2r_j}\setminus B_{r_j/2}$ and away from the next box,
but long enough that, by the elementary lower-capacity estimate of
Lemma~\ref{pt:lower-capacity}, its range has capacity of order
$m_j/\log m_j\asymp r_j^2/j$ with probability close to one. We stop
each trial at time $m_j$ or at the first failure of confinement,
whichever comes first. This stopped trial finishes before the next
box exit, even when confinement fails, so the fate of the $j$th piece
is known before the $(j+1)$st trial begins. Iterated conditioning
bounds the probability of failure at any $s$ prescribed scales by
$\zeta^s$, where $\zeta$ can be chosen arbitrarily small. A union bound
then bounds the probability that at least half of the $k\asymp\log n$
pieces fail by $2^k\zeta^{k/2}$. This is at most any prescribed negative
power of $n$ once $\zeta$ is small enough. Finally, the outermost box has radius
$r_{j_++1}\le\sqrt n/\log n$, so the walk leaves it well before time
$n/2$, which is what makes the whole construction measurable from the
first half of the walk.

\begin{proof}
Let $\mathcal F_t=\sigma(S_0,\ldots,S_t)$ and
$\tau_j=\inf\{t\ge0:S_t\notin B_{r_j}\}$. After time $\tau_j$ we
try to keep the walk near its exit point for a small multiple of $r_j^2$
steps. This supplies the capacity in \eqref{conc:piece-output}, while
leaving the next annulus unexplored.

We use the following elementary exponential maximal bound: for $r>0$ and $m\ge1$,
\begin{equation}
\label{conc:maximal}
 \Pp_0\left(\max_{0\le t\le m}|S_t|>r/4\right)
 \le 8\exp\left(-\frac{r^2}{32m}\right).
\end{equation}
It follows from Hoeffding's lemma~\cite{Hoeffding} for the bounded
symmetric increments of each coordinate, which makes
$\exp\{\lambda S_t^{(a)}-t\lambda^2/2\}$ a nonnegative supermartingale,
together with the maximal inequality for that supermartingale, the
choice $\lambda=u/m$, and a union bound over the eight coordinate
directions. For $m=0$ the probability is zero.

Fix $\zeta\in(0,1/16)$, to be chosen at the end of the proof, and choose
\[
 0<\eta\le\min\left\{\tfrac12,
                  \frac1{32\log(16/\zeta)}\right\}.
\]
For $j\in\mathcal J_n$, let $m_j=\lfloor\eta r_j^2\rfloor$ and define
\[
 \rho_j=\inf\{s\ge0:
       |S_{\tau_j+s}-S_{\tau_j}|>r_j/4\},
 \qquad
 \sigma_j=\tau_j+\min\{\rho_j,m_j\}.
\]
We call the $j$th piece good if
\begin{equation}
\label{conc:good-piece}
 \rho_j>m_j
 \quad\hbox{and}\quad
 \capc\{S_{\tau_j},\ldots,S_{\tau_j+m_j}\}
        \ge b\frac{r_j^2}{j}.
\end{equation}
When it is good, $E_j$ denotes the range in this display.

Conditionally on $\mathcal F_{\tau_j}$, the translated future is a free
walk. By \eqref{conc:maximal} and the choice of $\eta$,
\[
 \Pp_0(\rho_j\le m_j\mid\mathcal F_{\tau_j})
 \le8e^{-1/(32\eta)}\le\zeta/2.
\]
Without conditioning on confinement, the unrestricted
$m_j$-step segment has the free-walk range law. For sufficiently large
$j$, $m_j\ge\eta r_j^2/2$ and
$\log(m_j+2)\le C_{\log}j$, where $C_{\log}$ is universal.
Lemma~\ref{pt:lower-capacity}, used with error parameter $\zeta/2$,
therefore gives
\[
 \Pp_0\left(
  \capc\{S_{\tau_j},\ldots,S_{\tau_j+m_j}\}
             <b\frac{r_j^2}{j}
             \,\middle|\,\mathcal F_{\tau_j}\right)\le\zeta/2
\]
provided
$0<b\le c_{\mathrm{en}}\zeta\eta/(4C_{\log})$.
Here $c_{\mathrm{en}}$ is the fixed constant in
Lemma~\ref{pt:lower-capacity}; increasing $n_0$ ensures that all indices
being used are large enough.

Failure of \eqref{conc:good-piece} means that at least one of its two
requirements fails. The elementary inequality
$\Pp(A\cup B\mid\mathcal F)\le
\Pp(A\mid\mathcal F)+\Pp(B\mid\mathcal F)$ thus gives
\begin{equation}
\label{conc:conditional-good}
 \Pp_0(j\text{ is not good}\mid\mathcal F_{\tau_j})\le\zeta.
\end{equation}
No independence between the confinement and capacity requirements is
being assumed.

The decision whether a piece is good is known at $\sigma_j$, including
when confinement fails. Until that time the walk stays within distance
$r_j/4+1$ of $S_{\tau_j}$, and
$|S_{\tau_j}|\le r_j+1$. Since $r_{j+1}=32r_j$, this implies
$\sigma_j<\tau_{j+1}$. Consequently, for any prescribed $s$ indices,
iterated conditioning using \eqref{conc:conditional-good} bounds the
probability that all of them fail by $\zeta^s$. The pieces need not be
independent. A union bound over subsets of indices of size
$\lceil k/2\rceil$ gives
\begin{equation}
\label{conc:many-good}
 \Pp_0\left(\#\{j\in\mathcal J_n:j\text{ is not good}\}\ge k/2\right)
 \le2^k\zeta^{k/2}.
\end{equation}
By \eqref{conc:scale-count}, $\zeta$ can be fixed so small, depending
only on $\varepsilon,q$, that this is at most $n^{-q}$.
The choices are made in the order $\varepsilon,q$, then $\zeta$, then
$\eta$, then $b$; none depends on $n$.

Every good piece lies in $B_{2r_j}\setminus B_{r_j/2}$ for all indices
under consideration: its initial point has norm between $r_j$ and
$r_j+1$, and the piece remains within distance $r_j/4$ of it.
By the exit-time tail \eqref{in:exit} and
$r_{j_++1}\le\sqrt n/\log n$,
\[
 \Pp_0(\tau_{j_++1}>\lfloor n/2\rfloor)
 \le C\exp\{-c(\log n)^2\}\le C_q n^{-q}.
\]
Let $F_n$ be the event that $\tau_{j_++1}\le\lfloor n/2\rfloor$ and
at least $k/2$ pieces are good. On this event every selected piece lies
in the first half of the walk. Moreover, $F_n$ is measurable from that
half: if the outer exit has occurred, all decisions have already been
made, while otherwise $F_n$ is false. Equations
\eqref{conc:many-good} and the last display prove
\eqref{conc:first-half-event} and complete the construction.
\end{proof}

The next lemma no longer concerns the random mechanism producing the
pieces. It explains why the geometry and the capacities in
\eqref{conc:piece-output} imply a uniform positive potential.

\begin{lemma}[Hitting separated annular pieces]
\label{conc:deterministic-potential}
Let $1\le j_-\le j_+$, $k=j_+-j_-+1$, and $r_j=32^j$.
Suppose that for a subset $I\subset\{j_-,\ldots,j_+\}$ of size at least
$k/2$ there are finite sets $E_j\subset B_{2r_j}$ satisfying
$\capc(E_j)\ge br_j^2/j$, where $b>0$ is fixed. Then
\begin{equation}
\label{conc:anchored-potential}
 \inf_{y\in B_{r_{j_-}/4}}h_{\bigcup_{j\in I}E_j}(y)
 \ge1-\exp\left(-\frac{\beta_{\mathrm{hit}}k}{2j_+}\right),
 \qquad
 \beta_{\mathrm{hit}}=\min\{c_{\mathrm{hit}}b,\tfrac12\}>0,
\end{equation}
where $c_{\mathrm{hit}}$ is the universal constant in
\eqref{pt:interior-hit}.
\end{lemma}

\begin{proof}
Put $D_j=B_{4r_j}$. The killed last-exit identity
\eqref{pt:killed-last-exit}, the inequality
$e_{E_j,D_j}\ge e_{E_j}$, and the interior killed-Green bound give
\begin{align}
 \Pp_x(H_{E_j}<\tau_{D_j})
 &=\sum_{y\in E_j}G_{D_j}(x,y)e_{E_j,D_j}(y)\notag\\
 &\ge c_{\mathrm{hit}}r_j^{-2}\capc(E_j)
 \ge\frac{\beta_{\mathrm{hit}}}{j},
 \qquad x\in B_{r_j/4},\quad j\in I.
\label{conc:annular-hit}
\end{align}
Apply this estimate to a test walk started at
$y\in B_{r_{j_-}/4}$. After exiting $D_j$, its position has norm at
most $4r_j+1<r_{j+1}/4$. The strong Markov property at successive exit
times therefore gives, at each index in $I$, a conditional probability
at least $\beta_{\mathrm{hit}}/j$ to hit $E_j$ before the next outer
exit. The bound is uniform in the entrance position, including when we
condition on the entire stopped past and on avoiding all earlier pieces.
It follows that
\[
 \Pp_y\left(H_{\bigcup_{j\in I}E_j}=\infty\right)
 \le\prod_{j\in I}\left(1-\frac{\beta_{\mathrm{hit}}}{j}\right)
 \le\exp\left(-\frac{\beta_{\mathrm{hit}}|I|}{j_+}\right).
\]
Since $|I|\ge k/2$, this proves \eqref{conc:anchored-potential}.
The same deterministic argument applies to every starting point in the
inner box; no union bound over starting points is needed.
\end{proof}

\begin{proof}[Proof of Proposition~\ref{conc:sausage}]
Use Lemma~\ref{conc:annular-pieces} with $q=Q+2$.
On its event $F_n$, Lemma~\ref{conc:deterministic-potential} and
\eqref{conc:scale-count} imply
\begin{equation}
\label{conc:localized-potential}
 \inf_{y\in B_{n^\varepsilon}}
       h_{U_{\lfloor n/2\rfloor}}(y)\ge a,
 \qquad
 a=1-\exp\left(-\frac{\beta_{\mathrm{hit}}\gamma_\varepsilon}
                           {2\Gamma}\right)>0.
\end{equation}
Here $B_{n^\varepsilon}\subset B_{r_{j_-}/4}$, and enlarging the target
set can only increase its hitting potential. The constants in this
display depend only on $\varepsilon,Q$.

For even $n$, an event determined by the first $n/2$ steps has bridge
probability obtained by multiplying its free-walk law by
\[
 \frac{p_{n/2}(S_{n/2},0)}{p_n(0,0)}\le C,
\]
where the bound follows from \eqref{in:walk}. Since $F_n$ is measurable
from that half, its failure probability under $\Pp_n^{\mathrm{br}}$ is
at most $Cn^{-Q-2}$. A cyclic shift of the $n$ increments of a bridge,
followed by translation of its starting vertex to zero, has the same
bridge law. Apply \eqref{conc:localized-potential} at each of the $n$
cyclic roots. A union bound shows that, except with probability at most
$Cn\,n^{-Q-2}$, the potential of the full range is at least $a$ on the
ball of radius $n^\varepsilon$ around every visited vertex. This proves
the bridge part of \eqref{conc:halo-potential}.

For the free walk, at each time $s\in\{0,\ldots,n\}$ use the next
$\lfloor n/2\rfloor$ increments if $s\le n/2$, and that many preceding
increments in reverse order otherwise. After translation to $S_s$,
each chosen segment has the free-walk law; symmetry of the increment
law justifies reversal. There are enough increments in both cases.
The same estimate and a union bound over $n+1$ roots prove the free-walk
part of \eqref{conc:halo-potential}. Increasing $C$ handles the finite
set of $n$ below the construction threshold.

To deduce the capacity statement, let $U\subset F$ be finite sets and
suppose $h_U\ge a$ on $F$. By Green symmetry and last exit,
\[
 a\capc(F)
 \le\sum_{x\in F}e_F(x)h_U(x)
 =\sum_{y\in U}e_U(y)h_F(y)
 =\capc(U),
\]
because $h_F=1$ on $U$. Taking $F=U_n+B_{n^\varepsilon}$ proves
\eqref{conc:halo-capacity} under each of the two laws.

Finally, the mean estimate \eqref{in:range-mean}, from
Asselah, Schapira and Sousi~\cite[(1.4)]{ASS}, gives
\[
 \E_0[\capc(U_n)]\le C\frac{n}{\log(n+2)}.
\]
For even $n$, split the bridge into two halves. Subadditivity and the
bounded first-half density, reversing the second half, give
\[
 \E_n^{\mathrm{br}}[\capc(U_n)]
 \le C\E_0[\capc(U_{n/2})]
 \le C\frac{n}{\log(n+2)}.
\]
The deterministic bound
$\capc(U_n+B_{n^\varepsilon})\le C(n+1)n^{4\varepsilon}$, together
with \eqref{conc:halo-capacity} at the fixed choice $Q=5$, bounds the
exceptional contribution to either expectation by
$C(n+1)n^{4\varepsilon-5}$. On the complementary event it is bounded
by $a^{-1}$ times the ordinary-capacity expectation just estimated.
This proves both assertions of \eqref{conc:halo-mean}.
\end{proof}
\subsection{A single good-shape event for all killing scales}

Fix $\theta=1/100$ as in the rest of the paper, and fix a constant $D\ge1$
(the application below uses $D=10$).
For real $m\ge2$, define the integer radius
\begin{equation}
\label{conc:halo-radius}
 H(m)=\left\lceil Dm^\theta[\log(m+2)]^2\right\rceil.
\end{equation}
All constants in the next proposition may depend on $D$ and $\theta$.
The logarithmic factor in $H$ will later absorb the tails of components
formed by short loops.  The small power $m^\theta$ allows the block
argument below to retain a stretched-exponential error.

\begin{proposition}[Simultaneous massive-capacity profile]
\label{conc:profile}
There exists $C_0=C_0(D,\theta)<\infty$ such that the event
\begin{equation}
\label{conc:good-profile}
 \mathcal G_N=
 \left\{
 \capc_{e^{-1/M}}\bigl(U_N+B_{H(\min\{M,N\})}\bigr)
 \le\frac{C_0N}{\log(\min\{M,N\}+2)}
 \text{ for every real }M\ge2
 \right\}
\end{equation}
satisfies, for every $A>0$,
\begin{equation}
\label{conc:profile-polynomial}
 \Pp_N^{\mathrm{br}}(\mathcal G_N^c)\le C_A N^{-A},
 \qquad N\ge2\text{ even},
\end{equation}
where $C_A=C_A(A,D,\theta)<\infty$.  In particular, $C_0$ is chosen
before $A$ and does not change when a stronger polynomial error is
requested.  The property $\mathcal G_N$ depends only on the range and
duration; it is invariant under spatial translation and cyclic
rerooting of a loop.
\end{proposition}

The proof first works under the free-walk law, where translated ranges
of disjoint time blocks are independent. Blocks of length approximately
$\sqrt m$, with $m=\min\{M,N\}$, have a controlled mean thickened
capacity by Proposition~\ref{conc:sausage}; subadditivity bounds the
whole range by the sum of these block capacities. A deterministic bound
on each summand then gives concentration by Hoeffding's inequality.
Finally, dyadic bracketing handles all killing scales at once, and
conditioning on the final return costs only a polynomial factor.

\begin{proof}
We first establish the estimate at one fixed killing scale. Let
$M\ge2$, put $m=\min\{M,N\}$, and set $s=\lfloor\sqrt m\rfloor$.
For convenience we prove the estimate with the slightly larger radius
$H(2m)$; this will cover the dyadic interpolation at the end.  Choose
once and for all an exponent $\varepsilon$ with
$2\theta<\varepsilon<1/2$.
For every integer block length $\ell$ comparable to $s$, we have
\[
 H(2m)\le \ell^\varepsilon
\]
when $m$ exceeds a threshold depending only on $D,\theta,\varepsilon$.
Indeed $H(2m)=O(m^\theta\log^2 m)$, whereas
$\ell^\varepsilon$ has order $m^{\varepsilon/2}$ and
$\varepsilon/2>\theta$.
Proposition~\ref{conc:sausage}, whose mean estimate uses the fixed
confidence exponent $5$,
therefore gives
\begin{equation}
\label{conc:block-ordinary-mean}
 \E_0\left[\capc(U_\ell+B_{H(2m)})\right]
 \le C\frac{\ell}{\log(m+2)}.
\end{equation}
The comparability of $\log(\ell+2)$ and $\log(m+2)$ is used here.
For bounded $m$, the cardinality bound proves the same estimate after
increasing $C$.

The set in \eqref{conc:block-ordinary-mean} is deterministically
contained in $B_{\ell+H(2m)}$.
Since
$\ell+H(2m)+1\le C\sqrt m$ and
$1-e^{-1/M}\le1/M$, Lemma~\ref{pt:massive-comparison} implies
\begin{equation}
\label{conc:block-massive-mean}
 \E_0\left[\capc_{e^{-1/M}}(U_\ell+B_{H(2m)})\right]
 \le C\frac{\ell}{\log(m+2)}.
\end{equation}
Here the comparison factor is at most $C(1+m/M)\le2C$. The deterministic
inclusion $U_\ell\subset B_\ell$ is coarse, but it already gives a
uniformly bounded comparison factor at this block length, so no sharper
radius estimate is needed.

Partition the $N$ increments into
$b=\lceil N/s\rceil$ consecutive blocks of lengths
$\ell_1,\ldots,\ell_b$ differing by at most one and summing to $N$.
They can be chosen with $s/2\le\ell_r\le s$ when $s\ge2$; when $s=1$
all lengths equal one.  Write $U^{(r)}$ for the $r$th block range
translated so that its initial vertex is zero, and define
\[
 Y_r=\capc_{e^{-1/M}}(U^{(r)}+B_{H(2m)}).
\]
The random variables $Y_1,\ldots,Y_b$ are independent because each is a
function of a disjoint block of independent walk increments.
Translation invariance of capacity and subadditivity give
\[
 \capc_{e^{-1/M}}(U_N+B_{H(2m)})\le\sum_{r=1}^bY_r.
\]
By \eqref{conc:block-massive-mean}, there is a fixed constant $C_1$ such
that
\[
 \sum_{r=1}^b\E_0[Y_r]\le\frac{C_1N}{\log(m+2)}.
\]
The deterministic bounds $\capc_z(A)\le|A|$ and
$|U^{(r)}+B_h|\le(\ell_r+1)|B_h|$ give another fixed constant $C_2$ with
\[
 0\le Y_r\le B_m,
 \qquad B_m=C_2s[H(2m)+1]^4.
\]
The balanced partition is useful here: every block has a length
comparable to $s$, so no separate short final block needs to be charged.

Hoeffding's inequality~\cite{Hoeffding}, in the form
\[
 \Pp_0\left(\sum_{r=1}^b(Y_r-\E_0[Y_r])>u\right)
 \le\exp\left(-\frac{2u^2}{bB_m^2}\right),
\]
now applies.  Choose a fixed $C_3>C_1$ and take
$u=(C_3-C_1)N/\log(m+2)$.  Since
$b\le2N/s$, $s\le\sqrt m$, and
$H(2m)+1\le C m^\theta[\log(m+2)]^2$, it follows that
\begin{align}
 &\Pp_0\left(
 \capc_{e^{-1/M}}(U_N+B_{H(2m)})
     >\frac{C_3N}{\log(m+2)}\right)\notag\\
 &\hspace{8mm}\le
 \exp\left(-\frac{cN}
 {m^{1/2+8\theta}[\log(m+2)]^{18}}\right)
 \le
 \exp\left(-\frac{cN^{1/2-8\theta}}
 {[\log(N+2)]^{18}}\right).
\label{conc:profile-deviation}
\end{align}
The power $18$ is the sum of the two logarithmic powers from the
squared deviation threshold and the sixteen powers from
$[H(2m)+1]^8$.  Because $\theta=1/100<1/16$, the power
$1/2-8\theta$ of $N$ is strictly positive.  Constants in this estimate
are fixed independently of any later requested exponent $A$.

It remains to control all real $M$, rather than one fixed value.
For every dyadic number $a=2^j$ with
$1\le j\le\lfloor\log_2N\rfloor$, apply
\eqref{conc:profile-deviation} at killing scale $a$, with the smaller
radius $H(\min\{2a,N\})$.  Thus, outside an event of probability at most
\begin{equation}
\label{conc:dyadic-failure}
 C\log(N+2)
 \exp\left(-\frac{cN^{1/2-8\theta}}{[\log(N+2)]^{18}}\right),
\end{equation}
all inequalities
\[
 \capc_{e^{-1/a}}(U_N+B_{H(\min\{2a,N\})})
 \le\frac{C_3N}{\log(a+2)}
\]
hold simultaneously.  If $2\le M<N$, choose the dyadic $a$ satisfying
$a\le M<2a$.  Capacity decreases as the killing parameter $z$ increases,
and it increases as its set increases.  Therefore
\[
 \capc_{e^{-1/M}}(U_N+B_{H(\min\{M,N\})})
 \le
 \capc_{e^{-1/a}}(U_N+B_{H(\min\{2a,N\})}).
\]
If $M\ge N$, choose instead the largest dyadic $a\le N$; the same
inequality holds, because its neighborhood is $B_{H(N)}$ and $M\ge a$.
In both cases $\min\{M,N\}<2a$ or $\min\{M,N\}=2a$, and hence
\[
 \frac{\log(\min\{M,N\}+2)}{\log(a+2)}\le2.
\]
Consequently $C_0=2C_3$ makes all the inequalities in
\eqref{conc:good-profile} follow from the finite set of dyadic tests.
The event in \eqref{conc:good-profile} is measurable, since for fixed $N$
there are only finitely many possible walk paths.

Finally condition the entire free walk on $S_N=0$.  For any free-walk
event $E$, $\Pp_N^{\mathrm{br}}(E)\le\Pp_0(E)/p_N(0,0)$.
The return lower bound $p_N(0,0)\ge cN^{-2}$ and
\eqref{conc:dyadic-failure} thus yield
\begin{equation}
\label{conc:bridge-profile-deviation}
 \Pp_N^{\mathrm{br}}(\mathcal G_N^c)
 \le CN^2\log(N+2)
 \exp\left(-\frac{cN^{1/2-8\theta}}{[\log(N+2)]^{18}}\right).
\end{equation}
This is at most $C_A N^{-A}$ for every fixed $A>0$, because a positive
power of $N$ dominates every power of $\log N$.  Only $C_A$, not the
threshold $C_0$, changes with $A$.  Although the proof used a particular
partition of the increments, the event itself is defined by
\eqref{conc:good-profile}; its translation and cyclic invariance follow
from that definition.
\end{proof}

\section{Small-loop enlargements and a deterministic-parent kernel}
\label{ker:section}\label{sec:kernel}

The concentration estimate of Proposition~\ref{conc:profile} concerns one
loop.  A chain in the loop soup, however, can join two large loops through
several much shorter loops.  We first show that allowing such short
connections only enlarges a loop by a thin neighborhood, up to an error
which is small in expected capacity.  We then integrate over the second
loop while keeping the first loop completely deterministic.  This last
point is needed when the estimate is iterated along a chain.

\subsection{Enlargement by a fresh soup of short loops}

For $L\ge2$, let $\mathcal L_{\alpha,<L}$ be a Poisson loop soup with
intensity measure
\[
 \alpha\one_{\{N(\ell)<L\}}\mu(d\ell).
\]
The cutoff is strict: only loops with duration less than $L$ are present.
In particular, $\mathcal L_{\alpha,<2}$ is empty.  Write
$C_{\alpha,<L}(x)$ for the vertex cluster of $x$ in this soup, including
$x$ itself, and set
\begin{equation}\label{ker:enlargement}
 \Cl_{\alpha,<L}(U)
   =\bigcup_{x\in U}C_{\alpha,<L}(x)
\end{equation}
for a deterministic finite set $U\subset\Z^4$.  By a \emph{fresh} soup
we mean a sample of this Poisson law independent of the deterministic
sets or loop marks appearing in the calculation.  We do not condition
this law on any previous cluster exploration.  Its expectation will be
denoted by $\E_\alpha^{\mathrm{aux}}$ when it is useful to distinguish
it from an expectation over the original soup.

Recall the fixed choices $\theta=1/100$, $D=10$, and
$\bar\alpha=1/4$, and the radius
\[
 H(m)=\left\lceil 10m^\theta[\log(m+2)]^2\right\rceil,
 \qquad m\ge2,
\]
from the preceding section.  Constants in this section are uniform for
$0\le\alpha\le\bar\alpha$.  These fixed geometric parameters are
not changed when an error exponent is chosen.

\begin{lemma}[Expected capacity of a short-loop enlargement]
\label{ker:enlargement-bound}
Let $N,M\ge2$, $m=\min(N,M)$, $z=e^{-1/M}$, and
$2\le L\le2m^\theta$.  If $U$ is deterministic and $|U|\le N+1$, then
for every $A>0$ there is a finite constant $C_{\mathrm{err}}(A)$ such that
\begin{align}
 \E_\alpha^{\mathrm{aux}}\!\left[
       \capc_z\bigl(\Cl_{\alpha,<L}(U)\bigr)\right]
 &\le \capc_z\bigl(U+B_{H(m)}\bigr)
       +C_{\mathrm{err}}(A)Nm^{-A},
       \label{ker:enlargement-capacity}\\
 \E_\alpha^{\mathrm{aux}}\!\left[|\Cl_{\alpha,<L}(U)|\right]
 &\le C_{\mathrm{vol}}Nm^{4\theta}.
       \label{ker:enlargement-volume}
\end{align}
Here $C_{\mathrm{vol}}$ is independent of $A$.  If $U=U(\ell)$ is the
range of a loop of duration $N$ satisfying $\mathcal G_N$, then
\begin{equation}\label{ker:good-enlargement}
 \E_\alpha^{\mathrm{aux}}\!\left[
       \capc_z\bigl(\Cl_{\alpha,<L}(U(\ell))\bigr)\right]
 \le C_{\mathrm{gr}}\frac{N}{\log(m+2)},
\end{equation}
where $C_{\mathrm{gr}}$ is another fixed constant.
\end{lemma}

\paragraph{Proof idea.}
A short-loop component attached at a vertex of $U$ usually stays inside
the prescribed neighborhood. We retain that neighborhood deterministically
and charge the part outside it by its number of vertices. Summing the
expected excess over seed vertices, rather than requiring every component
to be small simultaneously, avoids a loss when the parent duration is much
larger than the child duration.

\begin{proof}
The radius of a component rooted at $x$ is measured from $x$:
\[
 \operatorname{rad}_x C_{\alpha,<L}(x)
    =\max\{|y-x|_\infty:y\in C_{\alpha,<L}(x)\}.
\]
The tail estimate \eqref{eq:micro-tail} implies that each component is
finite almost surely.  Since $U$ is finite, the same is true of its
enlargement.  Put $F=U+B_{H(m)}$.  If
$y\in\Cl_{\alpha,<L}(U)\setminus F$, there is some $x\in U$ with
$y\in C_{\alpha,<L}(x)$ and
$\operatorname{rad}_x C_{\alpha,<L}(x)>H(m)$.  Monotonicity,
subadditivity, and the cardinality upper bound for massive capacity
therefore give the pointwise estimate
\begin{align*}
 \capc_z\bigl(\Cl_{\alpha,<L}(U)\bigr)
 &\le \capc_z(F)
       +|\Cl_{\alpha,<L}(U)\setminus F|\\
 &\le \capc_z(F)
       +\sum_{x\in U}|C_{\alpha,<L}(x)|
          \one_{\{\operatorname{rad}_xC_{\alpha,<L}(x)>H(m)\}}.
\end{align*}
The component moment bound \eqref{in:micro-moments}, together with
$L+1\le3m^\theta$, bounds each expectation in the sum by
\[
 C(L+1)^4
   \left(1+\frac{H(m)}{L+1}\right)^4
       \exp\left\{-c\frac{H(m)}{L+1}\right\}.
\]
Since
\[
 \frac{H(m)}{L+1}\ge\frac D3[\log(m+2)]^2,
\]
this expression is at most $C_A m^{-A}$ for every fixed $A>0$.
For example, after absorbing the polynomial factor in the exponential,
it is bounded by
$C m^{4\theta}\exp\{-c_D[\log(m+2)]^2\}$; multiplying by $m^A$
leaves a bounded function of $m\ge2$.  Summing over
$|U|\le N+1\le3N/2$ proves \eqref{ker:enlargement-capacity}.
No independence between the components rooted at distinct vertices
of $U$ is used.

For the volume estimate, use
\[
 |\Cl_{\alpha,<L}(U)|\le\sum_{x\in U}|C_{\alpha,<L}(x)|
\]
and the first moment in \eqref{in:micro-moments}.  This gives
$\E_\alpha^{\mathrm{aux}}\!\left[|\Cl_{\alpha,<L}(U)|\right]
 \le C(N+1)(L+1)^4\le C_{\mathrm{vol}}Nm^{4\theta}$.

Finally, if the parent loop is good, its defining profile bounds
$\capc_z(U+B_{H(m)})$ by $C_0N/\log(m+2)$.
Apply \eqref{ker:enlargement-capacity} with $A=1$ and use
$\sup_{m\ge2}\log(m+2)/m<\infty$ to obtain
\eqref{ker:good-enlargement}.  In particular, the error has been absorbed
before dividing by any charge associated with the parent loop.  This
keeps the estimate uniform even when $N/M$ is arbitrarily large.
\end{proof}

\subsection{Charges and the connection probability}

Call a loop $\ell$ \emph{good} if its range satisfies
$\mathcal G_{N(\ell)}$ of Proposition~\ref{conc:profile}, and
\emph{bad} otherwise.  The definition is intrinsic to the unrooted
loop: it is unchanged by translations or cyclic changes of the root.
We use the profile estimate with exponent $10$ and fix a constant
$C_{\mathrm{bad}}$ for which
\begin{equation}\label{ker:bad-probability}
 \Pp_n^{\mathrm{br}}(\ell\text{ is bad})
       \le C_{\mathrm{bad}}n^{-10}
\end{equation}
for every even integer $n\ge2$.  Here
$\Pp_n^{\mathrm{br}}=\Pp_0(\,\cdot\mid S_n=0)$ is the ordinary
$n$-step bridge law; the corresponding expectation is
$\E_n^{\mathrm{br}}$.

Assign the loop the nonnegative charge
\begin{equation}\label{ker:charge}
 w(\ell)=
 \begin{cases}
 \displaystyle\frac{N(\ell)}{\log(N(\ell)+2)},&\ell\text{ good},\\[2mm]
 N(\ell)^2,&\ell\text{ bad}.
 \end{cases}
\end{equation}
The exponent two in the exceptional charge is dictated by the bad-parent
cases in the proof of Proposition~\ref{ker:kernel} below. For a bad
parent no capacity profile is available, and the only control on the
massive capacity of its short-loop enlargement is the volume bound
\eqref{ker:enlargement-volume}, which is of order $Nm^{4\theta}$ rather
than $N/\log(m+2)$. The charge of a bad parent must dominate this
volume, so that after normalization by the charge the kernel entry is at
most a constant times $\alpha/j$. The charge $N$ alone would leave an
uncompensated factor $m^{4\theta}$, which grows with the child type. Since
$m\le N$ and $4\theta<1$, the quadratic charge satisfies
$Nm^{4\theta}\le N^{1+4\theta}\le N^2$ and gives the needed domination.
Its exceptional contribution under the bridge law is bounded by
\[
 \E_N^{\mathrm{br}}[N^2\one_{\{\ell\text{ bad}\}}]
 \le C_{\mathrm{bad}}N^{-8}.
\]
Thus the fixed exponent ten in \eqref{ker:bad-probability} is more
than sufficient for the quadratic charge used here. Stronger polynomial
charges would require choosing a correspondingly stronger error bound
from the profile estimate first. We have
\begin{equation}\label{ker:charge-capacity}
 \capc(U(\ell))\le C_w w(\ell),\qquad C_w=\max\{C_0,1\}.
\end{equation}
Indeed, for a good loop of duration $N$, apply its profile with $M=N$
and use
$\capc(U)\le\capc_{e^{-1/N}}(U+B_{H(N)})$.
For a bad loop, use $\capc(U)\le|U|\le N\le N^2$.

For a loop $\ell$, define its duration type by
\[
 i(\ell)=\lfloor\log_2N(\ell)\rfloor\in\{1,2,\ldots\}.
\]
Thus type $i$ consists of durations $2^i\le N<2^{i+1}$.
For types $i,j\ge1$, define
\begin{equation}\label{ker:cutoff}
 L_{ij}=2^{\lceil\theta\min(i,j)\rceil}.
\end{equation}
For two deterministic loops of types $i$ and $j$, respectively, set
\begin{equation}\label{ker:connection}
 h_\alpha(\ell,\ell')
   =\Pp_\alpha^{\mathrm{aux}}\!
     \left(\Cl_{\alpha,<L_{ij}}(U(\ell))
                          \cap U(\ell')\ne\varnothing\right).
\end{equation}
Thus the two endpoint ranges may be joined through any finite chain of
auxiliary loops of duration less than $L_{ij}$.  Since the enlargement
contains its seed set, intersecting endpoint ranges have connection
probability one, even when the auxiliary soup is empty.

\subsection{Rooting identities and the uniform kernel estimate}

We first record precisely how a loop which visits a fixed vertex is
charged under the unrooted measure.  Let
$\mu_n=\one_{\{N(\ell)=n\}}\mu$ and let $F$ be a nonnegative,
translation-invariant function of an unrooted loop of duration $n$.
For even $n$, the definition of $\mu$ gives
\begin{equation}\label{ker:rooting-identity}
 \int\one_{\{0\in U(\ell)\}}F(\ell)\,\mu_n(d\ell)
   =\frac{p_n(0,0)}n
      \E_n^{\mathrm{br}}\bigl[|U_n|F(\ell_n)\bigr],
\end{equation}
where $\ell_n=[S_0,\ldots,S_n]$ denotes the bridge viewed modulo cyclic
rerooting and $U_n=\{S_0,\ldots,S_{n-1}\}$.  To see this, write the
left side as
\[
 \frac1n\sum_{x\in\Z^4}
   \E_x\!\left[
      \one_{\{0\in\{S_0,\ldots,S_{n-1}\}\}}
          F([S_0,\ldots,S_n]);\ S_n=x\right].
\]
Translate the walk by $-x$ and use translation invariance of $F$.
For a fixed translated bridge, precisely $|U_n|$ choices of $x$
satisfy $-x\in U_n$.  Tonelli's theorem now gives
\eqref{ker:rooting-identity}.  The argument uses the rooted definition
of the loop measure, and therefore includes periodic loops with their
correct multiplicities without any extra division by a period.
The same rerooting device, converting the sum over base points of a
loop into the range of a bridge, is used in the author's
thesis~\cite[(4.2.7) and Figure~4.1]{Vog20}; the based-loop
representation of Chang and Sapozhnikov~\cite[(5)--(6)]{CS16} is of the
same nature.

A useful consequence, valid for every finite deterministic
$G\subset\Z^4$, is
\begin{equation}\label{ker:bad-hit-measure}
 \mu_n\bigl(U(\ell)\cap G\ne\varnothing,\ \ell\text{ bad}\bigr)
 \le |G|p_n(0,0)\Pp_n^{\mathrm{br}}(\ell\text{ bad})
 \le C|G|n^{-12}.
\end{equation}
For the first inequality, sum the indicator of a hit over the vertices
of $G$, translate each vertex to zero, and use
\eqref{ker:rooting-identity} with $F=\one_{\{\ell\text{ bad}\}}$
and $|U_n|\le n$.  The second inequality uses
\eqref{ker:bad-probability} and the heat-kernel upper bound.  For odd
$n$ the left side is zero.

\begin{proposition}[Uniform estimate with a deterministic parent]
\label{ker:kernel}
There exist fixed finite constants $C_K,C_{\mathrm{src}}$ such that,
for every deterministic loop $\ell$ of type $i\ge1$, every $j\ge1$,
and every $0\le\alpha\le\bar\alpha$,
\begin{equation}\label{ker:kernel-bound}
 \alpha\int_{i(\ell')=j}
       h_\alpha(\ell,\ell')w(\ell')\,\mu(d\ell')
       \le w(\ell)K_t(i,j),\qquad t=C_K\alpha,
\end{equation}
where
\begin{equation}\label{ker:scalar-kernel}
 K_t(i,j)=t\frac{i}{j\min(i,j)}
      =\begin{cases}ti/j^2,&j\le i,\\t/j,&j>i.\end{cases}
\end{equation}
The source charges satisfy
\begin{equation}\label{ker:source}
 b_i(\alpha):=
   \alpha\int_{i(\ell)=i}\one_{\{0\in U(\ell)\}}
                  w(\ell)\,\mu(d\ell)
       \le\frac{C_{\mathrm{src}}\alpha}{i}.
\end{equation}
The constants are independent of the parent loop, the two types,
and any later spatial or duration cutoff.
\end{proposition}

\paragraph{Proof idea.}
For a good child, the duration-band estimate converts its integrated charge
into the massive capacity of the parent's enlargement. A good parent has
a small capacity profile; a bad parent has a larger charge that absorbs
the coarser volume estimate. Bad children are controlled directly by their
small probability under the bridge law, using a rooting identity before
any parent is selected. The four cases therefore give one bound valid for
every parent loop, which is what permits later iteration.

\begin{proof}
Fix the parent and put
\[
 N=N(\ell),\quad M=2^j,\quad m=\min(N,M),\quad
 z=e^{-1/M},\quad G=\Cl_{\alpha,<L_{ij}}(U(\ell)).
\]
Since $N\ge2^i$,
\[
 L_{ij}\le2\,2^{\theta\min(i,j)}\le2m^\theta.
\]
The enlargement estimates therefore apply.  Within the entire child
duration band $M\le N(\ell')<2M$, the auxiliary cutoff $L_{ij}$
is the same.  In particular, $G$ can be sampled independently before
the child loop is integrated.  Tonelli's theorem rewrites the integral
in \eqref{ker:kernel-bound} as
\begin{equation}\label{ker:tonelli}
 \alpha\E_\alpha^{\mathrm{aux}}\!\left[
     \int_{M\le N(\ell')<2M}
          \one_{\{U(\ell')\cap G\ne\varnothing\}}
                    w(\ell')\,\mu(d\ell')\right].
\end{equation}
We split the child integral into good and bad loops, and treat both
possible states of the parent.  The elementary comparisons used below
are
\begin{equation}\label{ker:log-comparisons}
 \begin{split}
 i\log2&\le\log(N+2)\le3i\log2,\\
 j\log2&\le\log(M+2)\le2j\log2,\\
 \min(i,j)\log2&\le\log(m+2).
 \end{split}
\end{equation}

\emph{Good child and good parent.}
For good children the charge equals
$N(\ell')/\log(N(\ell')+2)$.  Dropping their good-loop indicator,
then applying the duration-band estimate of Lemma~\ref{pt:band}
conditionally on $G$, bounds their contribution to
\eqref{ker:tonelli} by
\[
 \frac{C\alpha}{\log(M+2)}
             \E_\alpha^{\mathrm{aux}}[\capc_z(G)]
 \le\frac{C\alpha N}{\log(M+2)\log(m+2)}.
\]
Divide by the good parent's charge $N/\log(N+2)$ and apply
\eqref{ker:log-comparisons}.  The result is at most
\begin{equation}\label{ker:gg}
 C_{gg}\alpha\frac{i}{j\min(i,j)}.
\end{equation}

\emph{Good child and bad parent.}
The same duration-band estimate applies, but we now use
$\capc_z(G)\le|G|$ and \eqref{ker:enlargement-volume}.  After
normalization by the bad parent's charge $N^2$, the result is at most
\[
 \frac{C\alpha m^{4\theta}}{N\log(M+2)}
 \le\frac{C\alpha}{j}.
\]
Here $m\le N$ and $4\theta<1$ imply
$m^{4\theta}/N\le1$.  Since
$i/(j\min(i,j))\ge1/j$, this contribution is at most
\begin{equation}\label{ker:bg}
 C_{bg}\alpha\frac{i}{j\min(i,j)}.
\end{equation}

\emph{Bad child.}
For any deterministic finite $G$, multiply
\eqref{ker:bad-hit-measure} by the bad-child charge $n^2$ and sum over
$M\le n<2M$.  This yields
\begin{equation}\label{ker:bad-band}
 \int_{\substack{M\le N(\ell')<2M\\\ell'\text{ bad}}}
       \one_{\{U(\ell')\cap G\ne\varnothing\}}
                  w(\ell')\,\mu(d\ell')
       \le C|G|\sum_{M\le n<2M}n^{-10}
       \le C|G|M^{-9}.
\end{equation}
Average over the auxiliary soup and use
\eqref{ker:enlargement-volume}.  The unnormalized bad-child
contribution is at most $C\alpha Nm^{4\theta}M^{-9}$.
For a good parent, normalization and
\eqref{ker:log-comparisons} give
\[
 C\alpha\log(N+2)m^{4\theta}M^{-9}
       \le C\alpha i\,2^{-\beta_{\mathrm{bad}} j},\qquad
 \beta_{\mathrm{bad}}=9-4\theta>0.
\]
The number
$B_{\mathrm{bad}}=\sup_{j\ge1}j^2 2^{-\beta_{\mathrm{bad}} j}$ is finite.
Since $j\min(i,j)\le j^2$, this last expression is bounded by
\begin{equation}\label{ker:gb}
 C_{gb}\alpha\frac{i}{j\min(i,j)},
 \qquad C_{gb}=CB_{\mathrm{bad}}.
\end{equation}
For a bad parent, the normalized contribution is bounded instead by
\[
 C\alpha\frac{m^{4\theta}}N M^{-9}
       \le C\alpha 2^{-9j}
       \le \frac{C\alpha}{j}
       \le C_{bb}\alpha\frac{i}{j\min(i,j)}.
\]
Thus all four combinations of parent and child states satisfy the
same scalar bound.  Fix, once and for all,
\[
 C_K=\max\{1,C_{gg}+C_{gb},C_{bg}+C_{bb}\}.
\]
This proves \eqref{ker:kernel-bound}.  Every constant used in this
choice was obtained uniformly for $\alpha\le\bar\alpha$; hence
$C_K$ is fixed before the small intensity threshold is selected.

For the source, apply \eqref{ker:rooting-identity} with the charge
restricted to each of the two shape classes.  The good contribution
at duration $n$ is at most
\[
 \frac{p_n(0,0)}n
       \E_n^{\mathrm{br}}
          \left[|U_n|\frac n{\log(n+2)}
                       \one_{\{\ell_n\text{ good}\}}\right]
 \le C\frac1{n\log(n+2)}.
\]
Its sum over $2^i\le n<2^{i+1}$ is at most $C/i$.
For the bad contribution, the same identity and
\eqref{ker:bad-probability} give
\[
 \frac{p_n(0,0)}n
        \E_n^{\mathrm{br}}
          [|U_n|n^2\one_{\{\ell_n\text{ bad}\}}]
       \le Cn^{-10}.
\]
Its sum over that band is at most $C2^{-9i}\le C/i$.
Multiplying the two estimates by $\alpha$ and fixing
$C_{\mathrm{src}}$ as the sum of their constants proves
\eqref{ker:source}.
\end{proof}

The bridge law enters this proof only in the integration of intrinsic
loop classes against $\mu$.  The parent loop in
Proposition~\ref{ker:kernel} is arbitrary and deterministic.  Thus the
proposition can be applied repeatedly after a Poisson expansion has
selected endpoint loops, without any claim that an explored loop
retains an unbiased bridge distribution.

\section{From the cluster to a sum over loop chains}\label{chain:section}\label{sec:chains}
The estimate of Proposition~\ref{ker:kernel} concerns two fixed endpoint loops and an independent soup of smaller loops. To use it for a cluster, we first express connectivity through selected loop occurrences. A deterministic pruning procedure identifies the selected occurrences. The factorial Mecke formula then removes them from the soup, and the BK inequality separates the remaining connections. This order of operations is important: no assertion about the conditional law of a previously explored loop is needed.

\subsection{Finite soups and disjoint occurrence}
Let $V\subset\Z^4$ be finite with $0\in V$, and let $T\ge2$. Throughout this section, write
\[
 \Lambda=\{\ell:U(\ell)\subset V,\ N(\ell)\le T\},
 \qquad \mu_\Lambda=\one_\Lambda\mu.
\]
Here $U(\ell)$ is the vertex range of $\ell$. The set $\Lambda$ consists of finitely many loop shapes: a bounded duration permits only finitely many nearest-neighbor paths rooted in the finite set $V$. Let $\eta$ be the Poisson point measure on $\Lambda$ with intensity $\alpha\mu_\Lambda$, and denote its law and expectation by $\mathbf P_\Lambda$ and $\mathbf E_\Lambda$. Its cluster of the origin is denoted by $C_\Lambda(0)$ and includes $\{0\}$.

We distinguish a \emph{shape} $\ell\in\Lambda$ from an \emph{occurrence} of that shape in $\eta$. Thus $\eta(\{\ell\})$ is an independent Poisson random variable of mean $\alpha\mu(\{\ell\})$, and several occurrences may have the same shape. They are nevertheless different elements when choosing disjoint witnesses or ordered tuples. The total number of occurrences is finite almost surely, although it is not bounded by a deterministic constant.

A configuration is a vector in $\mathbb N_0^\Lambda$, equipped with its power-set $\sigma$-field; the inequality $\xi\le\eta$ means comparison of the counts of every shape. An event $E$ on these configurations is increasing if adding occurrences cannot destroy it. A subconfiguration $\xi\le\eta$ is a witness for $E$ if $\xi\in E$. For increasing events $E_1,\ldots,E_m$, their \emph{disjoint occurrence}, denoted by $E_1\mathbin\Box\cdots\mathbin\Box E_m$, means that there are witness subconfigurations $\xi_1,\ldots,\xi_m$ with
\[
 \xi_1+\cdots+\xi_m\le\eta,
 \qquad \xi_k\in E_k\quad(1\le k\le m).
\]
Thus the witnesses may use occurrences of the same shape, but they may not use the same occurrence twice. An event which already holds in the empty configuration may use an empty witness. For $m=0$, disjoint occurrence is the certain event.

\begin{lemma}[BK inequality for the finite Poisson soup]\label{chain:poisson-bk}
For increasing events $E_1,\ldots,E_m$ as above,
\begin{equation}\label{chain:bk}
 \mathbf P_\Lambda(E_1\mathbin\Box\cdots\mathbin\Box E_m)
 \le\prod_{k=1}^m\mathbf P_\Lambda(E_k).
\end{equation}
\end{lemma}
This is the usual Poisson version of the BK inequality, specialized to a
finite, possibly atomic loop measure. In particular, it is the same
loop-disjoint inequality used in Vogel~\cite[Lemma~3.1, equation~(3.5)]{Vog26};
the disjointness concerns loop occurrences and does not require their
vertex or edge ranges to be disjoint. Vogel discusses multiplicities by
passing to continuous-time bridges, which leaves connectivity unchanged.
Here we retain the discrete loop shapes and their multiplicities
explicitly, because the subsequent Mecke expansion selects individual
occurrences. The proof below gives a direct derivation for this convention
and for any finite number of increasing events. No dimension or
small-intensity assumption is needed for this probabilistic inequality.

The reduction to ordinary Bernoulli BK is particularly direct in this finite setting. We replace each Poisson count by many labelled Bernoulli trials, apply BK to their disjoint coordinate witnesses, and then pass to the Poisson limit.
\begin{proof}
If $\Lambda$ is empty, the configuration is deterministic and the
claim is immediate. Otherwise, enumerate the shapes as $\ell^{(1)},\ldots,\ell^{(s)}$ and put
$\lambda_r=\alpha\mu(\{\ell^{(r)}\})$. For an integer $q\ge\max\{1,\max_r\lambda_r\}$, replace the Poisson count of shape $r$ by the sum of $q$ independent Bernoulli variables with success probability $\lambda_r/q$. The variables for different shapes are also independent. Interpret each successful coordinate as one labelled occurrence.

An increasing event on the count vector lifts to an increasing event on these $sq$ Bernoulli coordinates. The count-vector definition of disjoint occurrence above agrees with the existence of disjoint sets of successful coordinates witnessing the lifted events: occurrences in prescribed subconfigurations can be assigned distinct successful coordinates, and conversely disjoint successful-coordinate witnesses give subconfigurations whose sum is bounded by the full count vector. The increasing-event BK inequality of van den Berg and Kesten~\cite{BK}, iterated over $m$ events, therefore gives \eqref{chain:bk} for this binomial approximation. To make the limiting step explicit, denote the law of the binomial
count vector by $\mathbf P_\Lambda^{(q)}$. For each fixed
$z\in\mathbb N_0^s$, its probability converges to
$\prod_{r=1}^s e^{-\lambda_r}\lambda_r^{z_r}/z_r!$ as
$q\to\infty$. These limiting masses sum to one, so convergence on this
countable state space is in total variation. Therefore, for every fixed
subset $A\subset\mathbb N_0^s$,
\[
 \lim_{q\to\infty}\mathbf P_\Lambda^{(q)}(A)
       =\mathbf P_\Lambda(A).
\]
Apply this to $A=E_k$ for each $k$ and to
$A=E_1\Box\cdots\Box E_m$. Taking $q\to\infty$ in the
binomial BK bound gives exactly \eqref{chain:bk}.
\end{proof}

We will also use the ordered factorial Mecke formula. If $F$ is a nonnegative measurable function of $k$ marked shapes and a residual configuration, then
\begin{equation}\label{chain:mecke}
\begin{split}
 &\mathbf E_\Lambda\left[
 \sum_{(a_1,\ldots,a_k)\in\eta^k_{\ne}}
 F\left(\ell(a_1),\ldots,\ell(a_k),
       \eta-\sum_{r=1}^k\delta_{\ell(a_r)}\right)\right]\\
 &\hspace{8mm}=
 \alpha^k\int_{\Lambda^k}
 \mathbf E_\Lambda[F(\ell_1,\ldots,\ell_k,\eta)]
 \prod_{r=1}^k\mu_\Lambda(d\ell_r).
\end{split}
\end{equation}
The notation $\eta^k_{\ne}$ means ordered $k$-tuples of pairwise distinct occurrences, and $\ell(a)$ is the shape of occurrence $a$. Equal shapes are permitted on either side. Formula~\eqref{chain:mecke} is the multivariate Mecke formula; see Last and Penrose~\cite[Theorem 4.5]{LP}. In the present finite setting it also follows directly from the identity
$\mathbb E[(Z)_r f(Z-r)]=\lambda^r\mathbb E[f(Z)]$ for a Poisson random variable $Z$ of mean $\lambda$, where $(Z)_r=Z(Z-1)\cdots(Z-r+1)$, applied separately to each shape. This argument also explains why no factorial divisor occurs for ordered tuples.

\subsection{A deterministic pruning procedure}
A \emph{simple occurrence chain} is a sequence
$(a_0,\ldots,a_r)$ of pairwise distinct soup occurrences such that
\[
 U(\ell(a_{k-1}))\cap U(\ell(a_k))\ne\varnothing,
 \qquad 1\le k\le r.
\]
Thus it is a path in the intersection graph of loop occurrences. The
adjective ``simple'' refers to distinct occurrences, not to disjoint
geometric ranges. Its \emph{type word} is the integer sequence
$(i(\ell(a_0)),\ldots,i(\ell(a_r)))$, where
$i(\ell)=\lfloor\log_2N(\ell)\rfloor$.

Recall the fixed separation parameter $\theta\in(0,1)$. We call a finite
word of positive integers $(i_0,\ldots,i_m)$ \emph{reduced} if
\begin{equation}\label{chain:reduced}
 \max_{a<k<b}i_k\ge\theta\min(i_a,i_b)
 \quad\hbox{for every }0\le a<b\le m\hbox{ with }b\ge a+2.
\end{equation}
In words, a reduced word has no interval whose interior types are all
much smaller than both endpoint types. The condition is imposed on every
interval, rather than only on triples of consecutive types. Words of
length one or two satisfy it automatically.

\begin{figure}[tbp]
\centering
\begin{tikzpicture}[x=1cm,y=1cm,font=\small,
  kept/.style={circle,draw=blue!65!black,fill=blue!8,thick,minimum size=8mm},
  hidden/.style={circle,draw=orange!80!black,fill=orange!12,thick,minimum size=8mm}]
 \node[anchor=west] at (-.4,1.05) {Original occurrence chain};
 \node[kept] (a0) at (0,0) {$a_0$};
 \node[hidden] (a1) at (2.2,0) {$a_1$};
 \node[hidden] (a2) at (4.4,0) {$a_2$};
 \node[kept] (a3) at (6.6,0) {$a_3$};
 \node[kept] (a4) at (8.8,0) {$a_4$};
 \draw[thick] (a0)--(a1)--(a2)--(a3)--(a4);
 \node[below=2mm of a0] {$500$};
 \node[below=2mm of a1] {$2$};
 \node[below=2mm of a2] {$3$};
 \node[below=2mm of a3] {$400$};
 \node[below=2mm of a4] {$450$};
 \node[anchor=west] at (-.4,-1.35) {Retained sequence and its hidden witnesses};
 \node[kept] (b0) at (0,-3.15) {$a_0$};
 \node[kept] (b3) at (6.6,-3.15) {$a_3$};
 \node[kept] (b4) at (8.8,-3.15) {$a_4$};
 \draw[thick,dashed,orange!80!black] (b0)--node[above=2mm,align=center]
   {hidden witness $\{a_1,a_2\}$\\durations $<L_{500,400}=16$}(b3);
 \draw[thick] (b3)--node[below=2mm,align=center] {direct\\intersection}(b4);
\end{tikzpicture}
\caption{Pruning in the intersection graph. A solid edge in the upper
row records an intersection of the two loop ranges; the numbers are
their types, not their durations. In the lower row a dashed edge records
a connection through the indicated hidden occurrences, rather than a
direct intersection of the retained ranges. Only the five displayed
occurrences constitute this witness. Additional geometric intersections
are allowed and do not affect the argument.}
\label{chain:pruning-figure}
\end{figure}
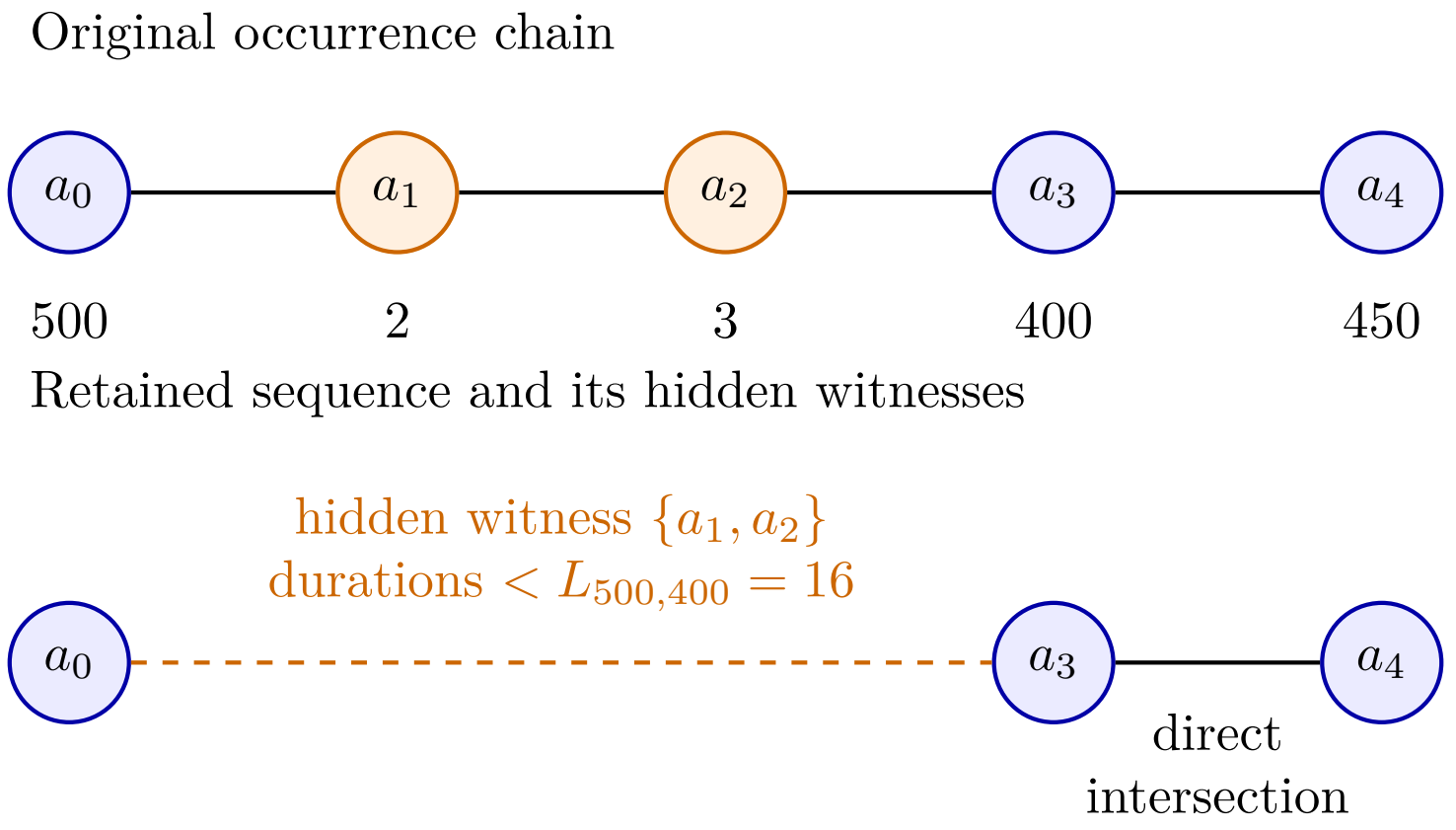

For example, take $\theta=1/100$ and a simple occurrence chain
$(a_0,\ldots,a_4)$. Suppose its type word is $(500,2,3,400,450)$.
The block $(a_1,a_2)$ may be deleted from the retained sequence because
\[
 \max\{i(\ell(a_1)),i(\ell(a_2))\}=3
       <\theta\min\{500,400\}=4.
\]
The remaining word $(500,400,450)$ is reduced. The original occurrences
$a_1,a_2$ have not been removed from the soup: they are now \emph{hidden}
from the retained sequence and are recorded as a witness connecting
$a_0$ to $a_3$. Each has duration less than $2^4=16$. The second pair
$a_3,a_4$ needs no hidden witness because these two occurrences already
intersect. Figure~\ref{chain:pruning-figure} depicts these two different
kinds of connection.

\begin{lemma}[Pruning a chain]\label{chain:pruning}
Every finite simple occurrence chain has a retained subsequence containing its first and last occurrences with the following properties.
\begin{enumerate}
\item The types of the retained occurrences form a reduced word.
\item Between any two consecutive retained occurrences of types $i,j$, the original chain supplies a connection using only hidden occurrences of durations strictly less than
\[
 L_{ij}=2^{\lceil\theta\min(i,j)\rceil}.
\]
\item The hidden witness sets for different consecutive pairs are pairwise disjoint and contain no retained occurrence.
\end{enumerate}
An empty hidden set is allowed when two consecutive retained occurrences already intersect.
\end{lemma}
\begin{proof}
Initially retain every occurrence. If a contiguous nonempty block of the currently retained interior occurrences has all its types strictly less than $\theta$ times the smaller type of its two retained endpoints, delete that block. Keep the endpoints. Each deletion reduces the number of retained occurrences, so after finitely many deletions no such block remains. This is precisely condition~\eqref{chain:reduced}.

We verify the quantitative statement about the deleted occurrences by induction over deletions. The induction invariant is that every hidden occurrence between two consecutive retained endpoints of types $i,j$ has type strictly less than $\theta\min(i,j)$. It holds initially because there are no hidden occurrences.

Consider a deletion with new endpoints of types $i,j$. A newly deleted retained occurrence has type less than $\theta\min(i,j)$ by the deletion rule. Every previously hidden occurrence absorbed into the new gap lay between a pair of consecutive old retained endpoints, at least one of which is newly deleted. If the type of that deleted endpoint is $v$, the induction hypothesis bounds the old hidden type by $\theta v$. Since $v<\theta\min(i,j)$ and $\theta<1$, this is strictly less than $\theta\min(i,j)$. The invariant therefore survives the deletion.

Let $q$ be the type of a final hidden occurrence. The invariant gives $q<\theta\min(i,j)$. As $q$ is an integer,
$q+1\le\lceil\theta\min(i,j)\rceil$, and hence
\[
 N(\ell)<2^{q+1}\le 2^{\lceil\theta\min(i,j)\rceil}=L_{ij}.
\]
All occurrences between two successive retained ones still follow their original order and intersect successively, so they supply the required connection. Finally, these hidden sets are disjoint contiguous portions of the original simple chain. They avoid all retained occurrences and cannot share an occurrence with one another.
\end{proof}

\subsection{Capacity as a test-walk hitting measure}
Let $\Omega=(\Z^4)^{\mathbb N_0}$ carry its cylinder $\sigma$-field, and let $\Pp_x$ be the law on $\Omega$ of a full-space simple random walk started at $x$. Define a finite measure on this trajectory space by
\begin{equation}\label{chain:test-measure}
 \nu_V(\mathcal A)=\sum_{x\in V}e_V(x)\Pp_x(\mathcal A),
 \qquad \mathcal A\subset\Omega\hbox{ measurable}.
\end{equation}
Its total mass is $\capc(V)$, so $\nu_V$ is generally not a probability measure. The starting vertex is weighted by the equilibrium measure $e_V$ rather than chosen uniformly in $V$. This weighting is what turns a hitting event into capacity. For a trajectory $w\in\Omega$, write
$\mathcal R(w)=\{w_n:n\ge0\}$. For every $A\subset V$, the last-exit formula and symmetry of $G$ give
\begin{align}
 \nu_V(\mathcal R(w)\cap A\ne\varnothing)
 &=\sum_{x\in V}e_V(x)h_A(x)\notag\\
 &=\sum_{x\in V}\sum_{y\in A}e_V(x)G(x,y)e_A(y)\notag\\
 &=\sum_{y\in A}e_A(y)h_V(y)
 =\capc(A).\label{chain:coverage}
\end{align}
The last equality uses $h_V(y)=1$ for $y\in V$. The walk in \eqref{chain:test-measure} is not killed on leaving $V$; this is why the identity gives full-space capacity. Identity~\eqref{chain:coverage} holds for each deterministic subset of $V$ and therefore holds, conditional on the soup, with $A=C_\Lambda(0)$.

\subsection{The capacity expansion}
For fixed endpoint loops $\ell,\ell'$ of types $i,j$, recall the connection probability from Section~\ref{ker:section},
\[
 h_\alpha(\ell,\ell')
 =\Pp_\alpha^{\mathrm{aux}}\bigl(\Cl_{\alpha,<L_{ij}}(U(\ell))\cap U(\ell')\ne\varnothing\bigr).
\]
Only the two displayed endpoint ranges are free deterministic connectors in this event. Every other connector must be an occurrence from the fresh soup with duration less than $L_{ij}$.

\begin{proposition}[Expansion over retained loop chains]\label{chain:path-series}
For every finite $V\ni0$, every $T\ge2$, and every $\alpha\ge0$,
\begin{equation}\label{chain:series}
\begin{split}
 \mathbf E_\Lambda[\capc(C_\Lambda(0))]
 \le{}&\capc\{0\}\\
 &+\sum_{m=0}^\infty\alpha^{m+1}
 \int_{\Lambda^{m+1}}
 \one_{\{0\in U(\ell_0)\}}
 \one_{\{(i(\ell_0),\ldots,i(\ell_m))\text{ reduced}\}}\\
 &\hspace{17mm}\times\capc(U(\ell_m))
 \prod_{k=1}^m h_\alpha(\ell_{k-1},\ell_k)
 \prod_{k=0}^m\mu_\Lambda(d\ell_k).
\end{split}
\end{equation}
An empty product is one. The right side is permitted to be infinite at this stage. All connection probabilities there are evaluated in unrestricted spatial volume, with the displayed duration cutoffs.
\end{proposition}
To prove this estimate, first test whether the cluster meets a fixed trajectory and prune one witnessing chain. The hidden loops then witness disjoint connection events between the retained loops. Mecke removes those retained occurrences, BK bounds the residual events by a product, and integration against $\nu_V$ replaces the terminal hitting event by its capacity. In particular, independence is introduced only after a pointwise connectivity bound has been established.
\begin{proof}
Fix a trajectory $w\in\Omega$. If $\mathcal R(w)$ hits $C_\Lambda(0)$, either it hits $0$ or there exists a finite simple occurrence chain whose first loop contains $0$ and whose last loop meets $\mathcal R(w)$. To see this, use a simple path in the finite graph whose vertices are the loop occurrences, with two vertices adjacent when their ranges intersect; the first vertex is incident to $0$ and the last to the point of intersection with $w$.

Apply Lemma~\ref{chain:pruning} to such a chain. Suppose that its retained occurrences are $(a_0,\ldots,a_m)$, with shapes $(\ell_0,\ldots,\ell_m)$. Remove all these marked occurrences and set
\[
 \eta^{\rm res}=\eta-\sum_{k=0}^m\delta_{\ell_k}.
\]
For each $k=1,\ldots,m$, let $E_k^{\boldsymbol\ell}$ be the increasing event that the two deterministic ranges $U(\ell_{k-1})$ and $U(\ell_k)$ are connected by occurrences of duration less than $L_{i(\ell_{k-1}),i(\ell_k)}$ in the configuration under consideration. Bare intersection counts as a connection. Other retained ranges are not added as deterministic connectors. By the pruning lemma, the residual configuration satisfies
\[
 \eta^{\rm res}\in
 E_1^{\boldsymbol\ell}\mathbin\Box\cdots\mathbin\Box E_m^{\boldsymbol\ell},
\]
with witnesses given by the disjoint hidden portions of the original chain.

It follows that the indicator of $\{\mathcal R(w)\cap C_\Lambda(0)\ne\varnothing\}$ is at most the indicator of $\{0\in\mathcal R(w)\}$ plus the sum, over all $m\ge0$ and all ordered distinct retained occurrence tuples, of
\begin{equation}\label{chain:tuple-indicator}
 \one_{\{0\in U(\ell_0)\}}
 \one_{\{U(\ell_m)\cap\mathcal R(w)\ne\varnothing\}}
 \one_{\{(i(\ell_0),\ldots,i(\ell_m))\text{ reduced}\}}
 \one_{\{\eta^{\rm res}\in E_1^{\boldsymbol\ell}\Box\cdots\Box E_m^{\boldsymbol\ell}\}}.
\end{equation}
This is a union bound over the retained tuples. For each tuple we ask only for the existence of disjoint hidden witnesses. We do not additionally count hidden chains or pruning histories.

Take $\mathbf E_\Lambda$ of this pointwise inequality and apply \eqref{chain:mecke} with $k=m+1$. The retained shapes become integration variables and the residual configuration becomes a fresh soup with law $\mathbf P_\Lambda$. For each such tuple, Lemma~\ref{chain:poisson-bk} gives
\begin{align*}
 \mathbf P_\Lambda\bigl(E_1^{\boldsymbol\ell}\Box\cdots\Box E_m^{\boldsymbol\ell}\bigr)
 &\le\prod_{k=1}^m\mathbf P_\Lambda(E_k^{\boldsymbol\ell})\\
 &\le\prod_{k=1}^m h_\alpha(\ell_{k-1},\ell_k).
\end{align*}
The second inequality is obtained by adding all allowed loop occurrences outside the finite restriction $\Lambda$ to each fresh connection soup. Connection events are increasing, so this enlargement can only raise their probabilities.

Finally integrate in $w$ against $\nu_V$. On the left, Tonelli's theorem
and \eqref{chain:coverage}, applied to each soup configuration, give
\begin{align*}
 &\int_\Omega\mathbf P_\Lambda
       (\mathcal R(w)\cap C_\Lambda(0)\ne\varnothing)\,\nu_V(dw)\\
 &\quad=\mathbf E_\Lambda\left[
       \int_\Omega\one_{\{\mathcal R(w)\cap C_\Lambda(0)\ne\varnothing\}}
                       \,\nu_V(dw)\right]
 =\mathbf E_\Lambda[\capc(C_\Lambda(0))].
\end{align*}
For a fixed retained tuple, put
\[
 F_m(\boldsymbol\ell)
 =\one_{\{0\in U(\ell_0)\}}
  \one_{\{(i(\ell_0),\ldots,i(\ell_m))\text{ reduced}\}}
  \prod_{k=1}^m h_\alpha(\ell_{k-1},\ell_k).
\]
This factor is independent of $w$. The contribution of that tuple in
the post-Mecke bound therefore integrates as
\begin{align*}
 &\int_\Omega F_m(\boldsymbol\ell)
       \one_{\{U(\ell_m)\cap\mathcal R(w)\ne\varnothing\}}\,\nu_V(dw)\\
 &\qquad=F_m(\boldsymbol\ell)\,
       \nu_V(\mathcal R(w)\cap U(\ell_m)\ne\varnothing)
 =F_m(\boldsymbol\ell)\capc(U(\ell_m)).
\end{align*}
Here $U(\ell_m)\subset V$, as required by
\eqref{chain:coverage}. Likewise,
$\int_\Omega\one_{\{0\in\mathcal R(w)\}}\,\nu_V(dw)=\capc\{0\}$.
All summands are nonnegative, so Tonelli also permits the integral in
$w$ to pass through the sum over $m$ and the loop integrals. Substituting
these evaluations yields \eqref{chain:series}.
\end{proof}

The independent soups appearing in the factors $h_\alpha$ are a consequence of the probability bound. They are not claimed to be independent pieces of the original cluster. Likewise, the tuple integrals may include equal geometric shapes: the factorial Mecke formula accounts for distinct occurrences of these shapes without any diagonal correction.

\section{Summation over all reduced type words}\label{word:section}\label{sec:words}
The preceding expansion still contains chains of arbitrary length. The role of pruning is to make their total contribution summable. We now prove a deterministic statement about the kernel of Proposition~\ref{ker:kernel}. For $t\ge0$ and positive integers $i,j$, write
\begin{equation}\label{word:kernel}
 K_t(i,j)=\frac{t}{j}\max\left\{1,\frac{i}{j}\right\}
 =\begin{cases}ti/j^2,&j\le i,\\t/j,&j>i.\end{cases}
\end{equation}
An unpruned sequence can make repeated large downward excursions and accumulate large factors $i/j$. Condition~\eqref{chain:reduced} limits the total logarithmic size of these downward excursions in terms of the word length and its maximal type.

\begin{lemma}[Summation of reduced words]\label{word:sum}
Fix $\theta\in(0,1)$ and define
\begin{equation}\label{word:t0}
 a_\theta=\log(1/\theta),\qquad
 b_\theta=2\bigl(a_\theta+\log2\bigr),\qquad
 t_0=\frac{\theta}{32e(a_\theta+\log2)}.
\end{equation}
There is a finite constant $C_\theta$, depending only on $\theta$, such that for every integer $J\ge1$, every initial type $i_0\in\{1,\ldots,J\}$, and every $0\le t\le t_0$,
\begin{equation}\label{word:sum-bound}
 S_{J,t}(i_0):=
 \sum_{m=0}^\infty\;
 \sum_{\substack{i_1,\ldots,i_m\in\{1,\ldots,J\}\\
                  (i_0,\ldots,i_m)\text{ reduced}}}
 \prod_{k=1}^mK_t(i_{k-1},i_k)
 \le C_\theta J^{1+16t/\theta}.
\end{equation}
For $m=0$, the inner sum and its empty product both equal one.
\end{lemma}
The proof counts oscillations of the sequence $\log i_k$ one height level at a time. Each additional excursion above a level requires an intermediate vertex lying within a fixed distance below that level; consequently, a vertex can account for such excursions only over a bounded interval of levels. This gives a bound on total variation. Logarithmic cells then convert the weighted discrete sum into a volume estimate, after which smallness of $t$ makes the sum over all lengths converge.
\begin{proof}
We separate the argument into a bound on variation and a count of words satisfying that bound. Set $H=\log J$, and for a fixed reduced word put $h_k=\log i_k$ for $0\le k\le m$. Thus $0\le h_k\le H$.

\emph{Step 1: variation of the logarithmic types.}
For $u\in(0,H)$, define
\begin{equation}\label{word:run-definition}
 c(u)=\one_{\{h_0\ge u\}}
       +\sum_{k=1}^m\one_{\{h_{k-1}<u\le h_k\}}.
\end{equation}
This counts the maximal consecutive runs in
$\{k\in\{0,\ldots,m\}:h_k\ge u\}$: the first term records a run
starting at index zero, and each summand records the start of a later
run. In particular, $c(u)=0$ if the set is empty. Consider a gap between two consecutive nonempty runs. If $r$ is the last index of the first run and $s$ the first index of the next, then $s\ge r+2$, $h_r,h_s\ge u$, and $h_k<u$ for $r<k<s$. The reduced-word condition gives
\[
 \max_{r<k<s}h_k
 \ge\min(h_r,h_s)-a_\theta\ge u-a_\theta.
\]
Thus this gap contains an index whose height lies in $[u-a_\theta,u)$. The interval condition is used here on the full gap, however many indices it contains. Different gaps contain different indices. Consequently, for every such $u$,
\begin{equation}\label{word:runs}
 c(u)\le1+\sum_{k=0}^m\one_{\{u-a_\theta\le h_k<u\}}.
\end{equation}
For a fixed $h_k$, the set of levels $u$ contributing to its indicator has length at most $a_\theta$. Integrating \eqref{word:runs} yields
\[
 \int_0^Hc(u)\,du\le H+a_\theta(m+1).
\]

Define the downward and total variations by
\[
 D(\boldsymbol h)=\sum_{k=1}^m(h_{k-1}-h_k)_+,
 \qquad
 V(\boldsymbol h)=\sum_{k=1}^m|h_k-h_{k-1}|.
\]
Except at the finitely many levels equal to some $h_k$, the number of downward crossings of a level $u$ is at most $c(u)$, since each run of heights above $u$ can end at most once. The number of crossings in either direction is at most $2c(u)$. Integrating these crossing counts, using
$(x-y)_+=\int\one_{\{y<u\le x\}}\,du$, gives
\begin{equation}\label{word:variation}
 D(\boldsymbol h)\le H+a_\theta(m+1),
 \qquad
 V(\boldsymbol h)\le2H+2a_\theta(m+1).
\end{equation}
The same bounds hold when $J=1$, when all heights are zero.

\emph{Step 2: a volume bound for the discrete sum.}
The kernel product has the exact form
\begin{equation}\label{word:product}
 \prod_{k=1}^mK_t(i_{k-1},i_k)
 =t^m\exp\{D(\boldsymbol h)\}\prod_{k=1}^m\frac1{i_k}.
\end{equation}
Indeed, for each $k$ the definition of the kernel gives
\[
 K_t(i_{k-1},i_k)
 =\frac{t}{i_k}\exp\{(\log i_{k-1}-\log i_k)_+\}
 =\frac{t}{i_k}e^{(h_{k-1}-h_k)_+}.
\]
Multiplying these identities for $1\le k\le m$ gives
\eqref{word:product}, by the definition of $D(\boldsymbol h)$.
For each integer $i\ge1$, associate the half-open interval
$I_i=[\log i,\log(i+1))$. These intervals are disjoint and satisfy
\[
 |I_i|=\log(1+1/i)\ge\frac1{2i},
 \qquad |x-\log i|\le\log2\quad(x\in I_i).
\]
The first inequality follows, for example, by integrating $1/s\ge1/(i+1)$ on $[i,i+1]$ and using $i+1\le2i$. These logarithmic cells are chosen because their lengths are comparable to the factors $1/i_k$ in \eqref{word:product}; an unweighted count of integer words would lose this useful factor.

Fix $m\ge1$. To each admissible word $(i_1,\ldots,i_m)$ associate its product cell $I_{i_1}\times\cdots\times I_{i_m}$ in $\mathbb R^m$. For a point $(x_1,\ldots,x_m)$ of this cell, put $x_0=h_0$. The triangle inequality and \eqref{word:variation} imply
\begin{align*}
 \sum_{k=1}^m|x_k-x_{k-1}|
 &\le V(\boldsymbol h)+2\sum_{k=1}^m|x_k-h_k|\\
 &\le2H+2a_\theta(m+1)+2m\log2\\
 &\le2H+b_\theta(m+1).
\end{align*}
The product cells of different words are disjoint, and their volumes dominate $2^{-m}\prod_{k=1}^m i_k^{-1}$. Their union is contained in the displayed variation region. Under the affine change of coordinates
$y_k=x_k-x_{k-1}$, this region becomes the $\ell^1$ ball of radius $2H+b_\theta(m+1)$; the Jacobian determinant is one. An $\ell^1$ ball of radius $r$ in $\mathbb R^m$ has volume $2^m r^m/m!$: in each of the $2^m$ orthants it is a simplex of volume $r^m/m!$. Therefore
\begin{equation}\label{word:cell-sum}
 \sum_{\substack{i_1,\ldots,i_m\le J\\
                 (i_0,\ldots,i_m)\text{ reduced}}}
 \prod_{k=1}^m\frac1{i_k}
 \le\frac{4^m[2H+b_\theta(m+1)]^m}{m!}.
\end{equation}
Combining \eqref{word:product}, the downward-variation bound, and \eqref{word:cell-sum}, the contribution of all words with $m$ edges is at most
\begin{equation}\label{word:length-bound}
 J e^{a_\theta}(4t e^{a_\theta})^m
 \frac{[2H+b_\theta(m+1)]^m}{m!}.
\end{equation}
The initial type $i_0$ entered only as the fixed coordinate $x_0$, so this estimate is uniform in $i_0$.

\emph{Step 3: summing the lengths.}
For $m\ge1$, use $(u+v)^m\le2^{m-1}(u^m+v^m)$ in \eqref{word:length-bound}. The resulting part containing $(2H)^m$ has sum
\[
 \frac{J e^{a_\theta}}2
 \sum_{m=1}^\infty\frac{(16t e^{a_\theta}H)^m}{m!}
 =\frac{J e^{a_\theta}}2
       \bigl(e^{16t e^{a_\theta}H}-1\bigr).
\]
For the other part, the elementary bound $m!\ge(m/e)^m$ gives
$(m+1)^m/m!\le e^{m+1}$. Its sum is at most
\[
 \frac{J e^{a_\theta}e}{2}
 \sum_{m=1}^\infty\rho^m,
 \qquad
 \rho=8t e^{a_\theta}b_\theta e.
\]
The choice of $t_0$ in \eqref{word:t0} ensures $\rho\le1/2$ for every $t\le t_0$. The geometric sum is then bounded uniformly in $t$. Since $e^{a_\theta}=1/\theta$ and $e^{16t e^{a_\theta}H}=J^{16t/\theta}$, adding the term $m=0$ proves \eqref{word:sum-bound}.
\end{proof}

For the value $\theta=1/100$ fixed in this paper, the exponent in \eqref{word:sum-bound} is uniformly bounded by
\begin{equation}\label{word:q}
 q_\theta=1+\frac{16t_0}{\theta}
 =1+\frac1{2e\log(2/\theta)}<\frac32.
\end{equation}
All chain lengths, including repeated types and arbitrarily many edges, have been included. No deterministic upper bound on the number of loop occurrences is used. The estimate deliberately keeps only a power of $J$, rather than optimizing its exponent: the largest type used for a box of radius $R$ will be of order $\log R$, so any fixed power of $J$ is already a logarithmic correction on the spatial scale. The stronger inequality $q_\theta<3/2$ leaves room for the sum over the initial type in the next section.

\section{The expected capacity of the cluster}\label{end:capacity}\label{sec:conclusion}
We now combine the estimates proved above. The purpose of first restricting
both the location and the duration of the loops is to make the marked-chain
expansion a finite-volume calculation with finitely many available loop shapes. The resulting bound will be independent
of the duration cutoff at the scale needed below.

Let $C_K$ be the constant in Proposition~\ref{ker:kernel}, and let $t_0$
be the constant in Lemma~\ref{word:sum}. Recall that $\bar\alpha=1/4$ and
$\theta=1/100$ have already been fixed. Define
\begin{equation}\label{end:intensities}
 \alpha_*:=\min\left\{\bar\alpha,\frac{t_0}{C_K}\right\},
 \qquad \alpha_0:=\frac{\alpha_*}{2},
 \qquad q_*:=q_\theta=1+\frac{16t_0}{\theta}.
\end{equation}
These are strictly positive constants independent of the radius. The
explicit choice of $t_0$ in Lemma~\ref{word:sum} gives
\[
 q_*=1+\frac{1}{2e(\log(1/\theta)+\log 2)}<\frac32.
\]
In particular, the constant $C_K$ is fixed before $\alpha_*$ is chosen.
An explicit numerical value of $t_0$ would not give a numerical value of
$\alpha_*$ without a corresponding value of $C_K$.
For the refined logarithmic bound, write $\kappa_{\log}:=16C_K/\theta$.
Then $1+\kappa_{\log}\alpha\le q_*$ whenever $\alpha\le\alpha_*$.

\begin{proposition}[Upper bound for the expected capacity]\label{end:cap-upper}
There is a finite constant $C_{\mathrm{cap}}$ such that
\begin{equation}\label{end:cap-bound}
 A_\alpha(R)\le \capc(\{0\})+
 C_{\mathrm{cap}}\alpha[\log(R+2)]^2,
 \qquad 0\le\alpha\le\alpha_*,\quad R\ge1.
\end{equation}
More precisely, there is a finite constant $C$, independent of $R$ and
$\alpha$, such that
\begin{equation}\label{end:cap-refined}
 \begin{gathered}
 A_\alpha(R)\le c_0+
 C\alpha[\log(R+2)]^{1+\kappa_{\log}\alpha}
                  [1+\log\log(R+2)],\\
 0\le\alpha\le\alpha_*,\qquad R\ge1.
 \end{gathered}
\end{equation}
\end{proposition}
The argument has two steps. For loops of duration at most $R^4$, the kernel and the reduced-word sum give a power of the largest type, which is of order $\log R$. Longer loops confined to $B_R$ have exponentially small total measure, so even the crude bound $\capc(B_R)=O(R^2)$ makes their contribution negligible. The power $R^4$ is a convenient cutoff chosen to make this final error small; it is not an assertion about the typical duration of a loop of diameter $R$.
\begin{proof}
Fix $R\ge1$ and set
\[
 T:=\max\{4,\lceil R^4\rceil\},\qquad
 J:=\lfloor\log_2 T\rfloor.
\]
Let $C_{R,T}(0)$ be the origin cluster formed by the loops whose ranges
are contained in $B_R$ and whose durations are at most $T$, together
with the vertex $0$. Every loop in this process has a type in
$\{1,\ldots,J\}$. Write $S_{J,t}(i)$ for the sum of reduced type words
starting at $i$, as in Lemma~\ref{word:sum}.

We now integrate the loop variables in the capacity expansion
\eqref{chain:series} of Proposition~\ref{chain:path-series}, starting
with the last retained loop and proceeding back towards the origin.
For each type $i$, let
\[
 \mu_{R,T,i}(d\ell)
 :=\one_{\{U(\ell)\subset B_R,\ N(\ell)\le T,\ i(\ell)=i\}}
       \mu(d\ell),\qquad t=C_K\alpha.
\]
Fix a reduced word $\boldsymbol i=(i_0,\ldots,i_m)$ with all types at
most $J$. The contribution of this word to the nonconstant part of
the upper bound \eqref{chain:series} for
$\E_\alpha[\capc(C_{R,T}(0))]$ is
\begin{equation}\label{end:fixed-word}
 I_{R,T}(\boldsymbol i)
 :=\alpha^{m+1}\int
   \one_{\{0\in U(\ell_0)\}}\capc(U(\ell_m))
   \prod_{k=1}^m h_\alpha(\ell_{k-1},\ell_k)
   \prod_{k=0}^m\mu_{R,T,i_k}(d\ell_k).
\end{equation}
The reduced-word condition is now a condition on the fixed indices,
so it does not change during the loop integrations.

Define nonnegative functions of a deterministic loop, working backwards
from the last position of the word, by
\begin{align}
 F_m(\ell_m)&:=\capc(U(\ell_m)),\notag\\
 F_{k-1}(\ell_{k-1})
 &:=\alpha\int h_\alpha(\ell_{k-1},\ell_k)
              F_k(\ell_k)\,\mu_{R,T,i_k}(d\ell_k),
       \qquad k=m,m-1,\ldots,1.\label{end:backward-functions}
\end{align}
We claim that, for every deterministic loop $\ell_k$ of type $i_k$,
\begin{equation}\label{end:backward-bound}
 F_k(\ell_k)\le C_w w(\ell_k)
             \prod_{r=k+1}^m K_t(i_{r-1},i_r),
 \qquad 0\le k\le m.
\end{equation}
At $k=m$ this is precisely the charge bound
\eqref{ker:charge-capacity}, with an empty product. If it holds at
position $k$, then enlarging the domain of the next nonnegative
integral and applying \eqref{ker:kernel-bound} gives
\begin{align*}
 F_{k-1}(\ell_{k-1})
 &\le C_w\left[\prod_{r=k+1}^mK_t(i_{r-1},i_r)\right]
     \alpha\int_{i(\ell_k)=i_k}
       h_\alpha(\ell_{k-1},\ell_k)w(\ell_k)\,\mu(d\ell_k)\\
 &\le C_w w(\ell_{k-1})
             \prod_{r=k}^mK_t(i_{r-1},i_r).
\end{align*}
This proves \eqref{end:backward-bound} by backward induction. The
parent $\ell_{k-1}$ is an integration variable held fixed at this step;
we use no assertion about its conditional distribution in an explored
cluster.

Tonelli's theorem identifies \eqref{end:fixed-word} with the remaining
source integral. Consequently,
\begin{align}
 I_{R,T}(\boldsymbol i)
 &=\alpha\int\one_{\{0\in U(\ell_0)\}}
             F_0(\ell_0)\,\mu_{R,T,i_0}(d\ell_0)\notag\\
 &\le C_w b_{i_0}(\alpha)
                    \prod_{r=1}^mK_t(i_{r-1},i_r),
                  \label{end:fixed-word-bound}
\end{align}
where $b_{i_0}(\alpha)$ is the full-space source charge in
\eqref{ker:source}. This also covers $m=0$. Summing the nonnegative
bounds \eqref{end:fixed-word-bound} over all reduced words in
Proposition~\ref{chain:path-series} gives
\begin{align}\label{end:truncated-cap}
 \E_\alpha[\capc(C_{R,T}(0))]
 &\le \capc(\{0\})+
 C_w\sum_{i=1}^J b_i(\alpha) S_{J,C_K\alpha}(i)\notag\\
 &\le \capc(\{0\})+
 C\alpha J^{1+\kappa_{\log}\alpha}\sum_{i=1}^J\frac1i\notag\\
 &\le \capc(\{0\})+
 C\alpha J^{1+\kappa_{\log}\alpha}(1+\log J).
\end{align}
The source estimate is used exactly once, at the first retained loop. Its factor $1/i$ makes the sum over initial types grow as $1+\log J$, while the uniform word bound pays for every subsequent step. The hypothesis $C_K\alpha\le t_0$ follows from
\eqref{end:intensities}. The global restriction that a word be reduced
is held fixed throughout the successive loop integrations; no Markov
property for reduced words is being asserted. The choice of $T$ gives
$J\le C\log(R+2)$. Since $1+\kappa_{\log}\alpha\le q_*$, the comparison
constants in
\[
 J^{1+\kappa_{\log}\alpha}(1+\log J)
 \le C[\log(R+2)]^{1+\kappa_{\log}\alpha}
                     [1+\log\log(R+2)]
\]
are uniform for $\alpha\le\alpha_*$. This is the logarithmic factor in
\eqref{end:cap-refined}. For the simpler bound \eqref{end:cap-bound},
use additionally $q_*<3/2$ and $1+\log J\le C\sqrt J$, which give
$J^{1+\kappa_{\log}\alpha}(1+\log J)\le C J^2$.
Thus the square logarithm is a convenient weaker formulation; retaining
the actual exponent improves the logarithmic correction itself.

It remains to restore loops of duration greater than $T$. Let
$P_{B_R}$ be the substochastic transition matrix of simple random walk
killed upon leaving $B_R$. Thus, for $x,y\in B_R$,
\[
 P_{B_R}^n(x,y)=\Pp_x(S_n=y,\ \tau_{B_R}>n),
 \qquad
 \tau_{B_R}:=\inf\{k\ge0:S_k\notin B_R\}.
\]
Here $\operatorname{Tr}(P_{B_R}^n)=\sum_{x\in B_R}P_{B_R}^n(x,x)$ is the sum of the diagonal entries of this finite matrix. The rooted-loop definition of $\mu$ gives the exact formula
\begin{equation}\label{end:long-mass}
 \lambda_{R,T}:=
 \mu\bigl(U(\ell)\subset B_R,\ N(\ell)>T\bigr)
 =\sum_{n>T}\frac1n\operatorname{Tr}(P_{B_R}^n).
\end{equation}
The exit-time estimate stated in the preliminaries yields
\[
 \sup_{x\in B_R}\Pp_x(\tau_{B_R}>n)
 \le C\exp\left\{-\frac{cn}{(R+1)^2}\right\}.
\]
Bounding each diagonal entry by this survival probability, using
$|B_R|\le C(R+1)^4$, and summing a geometric series, we obtain
\begin{align*}
 \lambda_{R,T}
 &\le \frac{C(R+1)^4}{T}
 \sum_{n>T}\exp\left\{-\frac{cn}{(R+1)^2}\right\}\\
 &\le \frac{C(R+1)^6}{T}
 \exp\left\{-\frac{cT}{(R+1)^2}\right\}
 \le C R^2 e^{-cR^2},\qquad R\ge1.
\end{align*}
Here and below the constants in successive lines may change. The last
inequality uses $T\ge R^4$, $T\ge4$, and $R+1\le2R$.

The quantity $\lambda_{R,T}$ is a loop-measure mass; the corresponding Poisson mean is $\alpha\lambda_{R,T}$. Couple $C_{R,T}(0)$ and $K_R$ using the same soup. If there is no loop
counted by $\lambda_{R,T}$, these two clusters are equal. Otherwise,
monotonicity bounds their capacity difference by
$\capc(B_R)\le C R^2$. The number of the omitted loops is Poisson with
mean $\alpha\lambda_{R,T}$. Consequently,
\begin{align*}
 0\le A_\alpha(R)-\E_\alpha[\capc(C_{R,T}(0))]
 &\le \capc(B_R)\bigl(1-e^{-\alpha\lambda_{R,T}}\bigr)\\
 &\le C\alpha R^4e^{-cR^2}\le C\alpha.
\end{align*}
The additive error $C\alpha$ is absorbed in both logarithmic bounds,
since $R\ge1$. Together with \eqref{end:truncated-cap}, this proves
\eqref{end:cap-refined} and \eqref{end:cap-bound}. This argument discards
long durations only after estimating their loop measure; it does not
replace a loop's duration by its squared diameter.
\end{proof}

\begin{remark}[The logarithmic upper exponent]\label{end:refined-log}
Combining \eqref{end:cap-refined} with the capacity--arm comparison of
Proposition~\ref{end:arm-comparison}, proved below, gives
\begin{equation}\label{end:arm-refined}
 \begin{gathered}
 p_\alpha(R)\le C\alpha R^{-2}
 \Bigl\{c_0+\alpha[\log(R+2)]^{1+\kappa_{\log}\alpha}
                       [1+\log\log(R+2)]\Bigr\},\\
 0<\alpha\le\alpha_*,\qquad R\ge1.
 \end{gathered}
\end{equation}
Consequently, for every fixed $0<\alpha\le\alpha_0$,
\begin{equation}\label{end:log-limsup}
 \limsup_{R\to\infty}
 \frac{\log(R^2p_\alpha(R))}{\log\log R}
 \le 1+\kappa_{\log}\alpha.
\end{equation}
Indeed, the extra factor $1+\log\log(R+2)$ contributes only
$O(\log\log\log R)$ after taking logarithms. The exponent of $R$
remains $-2$, while the upper exponent of $\log R$ improves from two
to $1+\kappa_{\log}\alpha$, up to this iterated logarithm. On the interval
asserted in Theorem~\ref{thm:main},
$1+\kappa_{\log}\alpha\le(1+q_*)/2<5/4$. This is an upper bound furnished
by the present summation argument, not a claimed optimal logarithmic
exponent.
\end{remark}

\section{One-arm probabilities and the lower logarithmic correction}
\label{end:arms}
The remaining argument has two parts. First, the expected capacity is
compared with the probability that a loop leaving the box attaches to
the interior cluster. Second, such attachments at different scales
supply disjoint intervals in the tail integral for capacity. Both
steps follow the framework of Gao~\cite[Section~7.3]{Gao}; we give the
details for the present four-dimensional normalization and keep the
constants uniform down to intensity zero.

\subsection{A uniform comparison between capacity and the one-arm event}

\begin{proposition}\label{end:arm-comparison}
There are constants $c_{\mathrm{arm}},C_{\mathrm{arm}}>0$ such that
\begin{equation}\label{end:comparison}
 c_{\mathrm{arm}}\alpha R^{-2}A_\alpha(R)
 \le p_\alpha(R)
 \le C_{\mathrm{arm}}\alpha R^{-2}A_\alpha(R),
 \qquad 0<\alpha\le\bar\alpha,\quad R\ge1.
\end{equation}
\end{proposition}
For the upper bound, split loops by spatial diameter. If the small-loop cluster already leaves the box, a disjoint annular crossing supplies a contracting factor; otherwise an arm requires a larger loop to hit an interior set, whose hitting measure is controlled by capacity. For the lower bound, condition on the loops contained in the box and ask for one independent good attaching loop. Both bounds must retain the factor $\alpha$, which is why the Poisson probabilities are estimated uniformly as $\alpha\downarrow0$.
\begin{proof}
We first prove the upper bound. Let $\mathcal L^{\le d}_\alpha$ denote
the soup restricted to loops of spatial diameter at most $d$. This is
a diameter restriction, in contrast with the duration restriction used
in the preceding sections. For $0<\delta<1/16$ and $R\ge1$, let
$\mathcal A_{R,\delta}$ be the event that a connected union of loops
of $\mathcal L^{\le\delta R}_\alpha$ meets both $B_{5R/8}$ and
$B_{3R/4}^{\,c}$.
See Figure~\ref{fig:small-crossing}. Equivalently, the event is
\[
 \begin{split}
 \mathcal A_{R,\delta}
 =\bigl\{&\exists\text{ a finite chain }(a_1,\ldots,a_k)
       \text{ in }\mathcal L^{\le\delta R}_\alpha:\\
       &U(\ell(a_1))\cap B_{5R/8}\ne\varnothing,
       \quad U(\ell(a_k))\cap B_{3R/4}^{\,c}\ne\varnothing\bigr\}.
 \end{split}
\]
Here a chain has successively intersecting ranges, as in
Section~\ref{sec:chains}; it need not be connected to the origin.

\begin{figure}[htbp]
\centering
\begin{tikzpicture}[x=0.82cm,y=0.82cm,font=\small,>=Stealth]
 \fill[blue!5] (-2.25,-2.25) rectangle (2.25,2.25);
 \fill[white] (-1.875,-1.875) rectangle (1.875,1.875);
 \draw[black!45] (-3,-3) rectangle (3,3);
 \draw[blue!65!black,thick] (-2.25,-2.25) rectangle (2.25,2.25);
 \draw[blue!65!black,dashed] (-1.875,-1.875) rectangle (1.875,1.875);
 \fill (0,0) circle (1.4pt) node[below left] {$0$};
 \foreach \k in {0,...,17} {
   \pgfmathsetmacro{\xx}{1.63+0.065*\k}
   \pgfmathsetmacro{\yy}{0.38+0.025*sin(70*\k)}
   \draw[orange!85!black,thick] (\xx,\yy) ellipse[x radius=0.052,y radius=0.075];
 }
 \draw[->,black!65] (3.6,1.1) -- (2.48,0.48);
 \node[anchor=west,align=left,text width=5.6cm] at (3.65,1.1)
   {A chain of small loops\\$\diam U(\ell)\le\delta R$};
 \node[anchor=west] at (3.65,-0.15) {The chain meets $B_{5R/8}$};
 \node[anchor=west] at (3.65,-0.62) {and reaches outside $B_{3R/4}$.};
 \node[anchor=south east,fill=white,inner sep=1pt] at (1.8,1.88) {$B_{5R/8}$};
 \node[anchor=south west,fill=white,inner sep=1pt] at (-2.25,2.25) {$B_{3R/4}$};
 \node[anchor=south west] at (-3,3) {$B_R$};
\end{tikzpicture}
\caption{The event $\mathcal A_{R,\delta}$, shown schematically in a
two-dimensional section. The orange chain crosses the shaded annulus
using only loops of diameter at most $\delta R$. It may lie away from
the origin. The outer box $B_R$ provides context for the later
prefix--suffix decomposition; the definition of
$\mathcal A_{R,\delta}$ does not require the chain to stay inside it.}
\label{fig:small-crossing}
\end{figure}
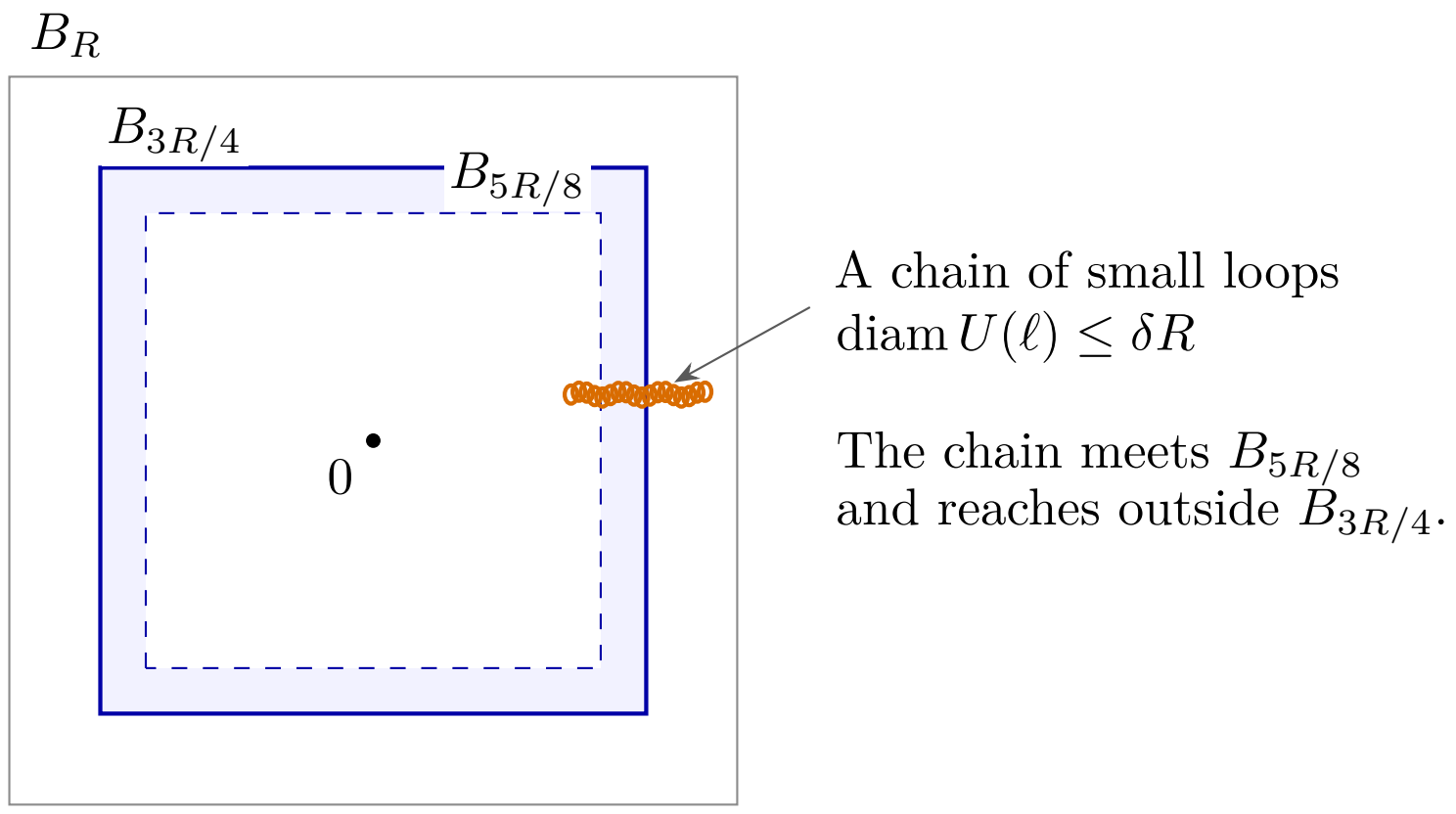

We claim that $\delta$ can be fixed sufficiently small, and then
$R_0\ge32$ sufficiently large, so that
\begin{equation}\label{end:small-crossing}
 \sup_{0\le\alpha\le\bar\alpha}
 \sup_{R\ge R_0}\Pp_\alpha(\mathcal A_{R,\delta})
 \le q,\qquad q:=\frac18.
\end{equation}
Here is the consequence of Gao's metric-graph estimate that is needed
for this claim. In the same-intensity coupling recalled in the
preliminaries, every lift of a discrete loop of diameter at most
$\delta R$ has diameter at most $\delta R+2$. For $\delta R\ge1$ this
is at most $L:=3\delta R$. Cover $B_{5R/8}$ by at most
$C\delta^{-4}$ metric boxes of radius $L$ with lattice centers. A
continuous path in the lifted connected union that starts in one of
these boxes and reaches $B_{3R/4}^{\,c}$ must cross a concentric
annulus from radius $L$ to radius $nL$, where
$n=\lfloor1/(96\delta)\rfloor$, after reducing $\delta$ so that
$n\ge2$. Indeed, its total displacement from its starting point is at
least $R/8$, while the center is within $L$ of that starting point and
$nL\le R/32$. Thus the outer sphere is reached before the path leaves
the connected union. The metric-graph estimate
of Proposition~\ref{in:metric-tail}, followed by the union bound over these boxes,
gives
\[
 \Pp_\alpha(\mathcal A_{R,\delta})
 \le C\delta^{-4}e^{-c/\delta},
 \qquad \delta R\ge1,
\]
where monotone coupling with intensity $\bar\alpha$ makes the constants
uniform in $\alpha$. The right side tends to zero as $\delta\downarrow0$.
Choose $\delta$ once so that this upper bound is at most $q=1/8$, and then take
$R_0\ge\max\{32,\delta^{-1}\}$. This proves
\eqref{end:small-crossing}. Bounded changes in the lattice boundaries
are absorbed by the strict separation of the radii in this covering
argument.

Let $K^{\mathrm{sm}}_{R,\delta}$ be the origin cluster in
$\mathcal L^{\le\delta R}_\alpha$, including $0$, and write
\[
 E_{R,\delta}:=
 \{K^{\mathrm{sm}}_{R,\delta}\not\subset B_R\}.
\]
On $E_{R,\delta}$ there is a finite simple chain of loop occurrences
from $0$ to $B_R^c$. Stop this chain at its first loop leaving
$B_{R/2}$. The preceding loops are contained in $B_{R/2}$, and the last
loop of this prefix intersects $B_{R/2}$ and has diameter at most
$\delta R$. Hence the entire prefix is contained in
$B_{R/2+\delta R}\subset B_{5R/8}$. It witnesses the full-soup arm
from $0$ to $B_{R/2}^c$. There must be a nonempty remaining suffix, since the prefix is still contained in $B_{5R/8}$. The suffix starts at a vertex of
this prefix and reaches $B_R^c$, so it witnesses
$\mathcal A_{R,\delta}$. The prefix and suffix use disjoint loop
occurrences. The Poisson BK inequality, first for finite restrictions
and then by increasing exhaustion, therefore gives
\begin{equation}\label{end:small-exit}
 \Pp_\alpha(E_{R,\delta})
 \le p_\alpha(R/2)\Pp_\alpha(\mathcal A_{R,\delta})
 \le q p_\alpha(R/2),\qquad R\ge R_0.
\end{equation}
This is the occurrence-disjoint decomposition underlying
Gao's recursion~\cite[(7.20)]{Gao}.

On $E_{R,\delta}^{\,c}$, the set $K^{\mathrm{sm}}_{R,\delta}$ is finite
and contained in $B_R$. If the full cluster reaches $B_R^c$, a chain
producing that arm must contain a first loop of diameter greater than
$\delta R$; that loop intersects $K^{\mathrm{sm}}_{R,\delta}$. The
large-loop soup is independent of the small-loop soup. Conditional on
the latter, the probability that at least one large loop hits a fixed
set $K$ is
\[
 1-\exp\{-\alpha\mu(U(\ell)\cap K\ne\varnothing,
                         \ \diam U(\ell)>\delta R)\}.
\]
Using $1-e^{-u}\le u$ and Gao's deterministic-set estimate
\eqref{in:macroscopic} on $E_{R,\delta}^{\,c}$ yields
\begin{align*}
 p_\alpha(R)
 &\le q p_\alpha(R/2)
   +C_\delta\alpha R^{-2}
       \E_\alpha\!\left[
          \one_{E_{R,\delta}^{\,c}}
          \capc(K^{\mathrm{sm}}_{R,\delta}\cap B_R)\right]\\
 &\le q p_\alpha(R/2)+C_\delta\alpha R^{-2}A_\alpha(R).
\end{align*}
For the second inequality, every small loop in the root component is
contained in $B_R$ on the indicated event; consequently
$K^{\mathrm{sm}}_{R,\delta}\subset K_R$. The containment hypothesis in
\eqref{in:macroscopic} is used only on this event.

For completeness, we iterate this recursion without restricting $R$
to a predetermined dyadic sequence. If $R\ge R_0$, choose $m\ge0$
so that $r:=2^{-m}R\in[R_0,2R_0)$. Since $A_\alpha$ is increasing,
\begin{equation}\label{end:iterate}
 p_\alpha(R)
 \le C_\delta\alpha R^{-2}A_\alpha(R)
        \sum_{k=0}^{m-1}(4q)^k
       +q^m p_\alpha(r).
\end{equation}
The finite sum is bounded by $(1-4q)^{-1}$. The factor $4$ comes from the scaling $(R/2)^{-2}=4R^{-2}$; the choice $q=1/8$ was made so that this scale factor is outweighed by the annular-crossing probability. Let
\[
 \lambda_0:=\mu(0\in U(\ell)).
\]
The rooted-loop identity and $|U_n|\le n$ imply
$\lambda_0\le\sum_{n\ge2}p_n(0,0)<\infty$. Every arm from the
origin requires a loop containing $0$, so
$p_\alpha(r)\le1-e^{-\alpha\lambda_0}\le\alpha\lambda_0$.
Moreover, $q^m\le4^{-m}=r^2/R^2\le4R_0^2/R^2$ and
$A_\alpha(R)\ge c_0:=\capc(\{0\})>0$. The last term of
\eqref{end:iterate} is therefore at most
$C\alpha R^{-2}A_\alpha(R)$. For $1\le R<R_0$, the same conclusion
follows directly from $p_\alpha(R)\le\alpha\lambda_0$ and
$R^{-2}A_\alpha(R)\ge c_0/R_0^2$. This proves the upper bound for all
real $R\ge1$.

For the lower bound, let $T_R(K)$ be the good-loop family from
Proposition~\ref{in:attachment}. Thus its loops are contained in $B_{2R}$,
intersect $K\subset B_R$, and contain a portion outside $B_R$ of
capacity at least $\kappa_0 R^2/\log R$. Denote its loop measure by
$m_R(K):=\mu(T_R(K))$. For $R\ge32$, the cited estimate is
\begin{equation}\label{end:attach-mass}
 c_{\mathrm{att}}R^{-2}\capc(K)
 \le m_R(K)\le C_{\mathrm{att}}R^{-2}\capc(K).
\end{equation}
For fixed $R$, $m_R(K)$ is a deterministic function of the finite set $K\subset B_R$; substituting $K_R$ therefore gives a well-defined random variable. Let $F_R$ be the event that a loop in the soup belongs to
$T_R(K_R)$, and set $g_\alpha(R):=\Pp_\alpha(F_R)$. The random set
$K_R$ is determined by the loops wholly contained in $B_R$, whereas
every member of $T_R(K_R)$ leaves $B_R$. Independence of the Poisson
restrictions to these disjoint classes, followed by conditioning on
$K_R$, gives the exact formula
\begin{equation}\label{end:attach-probability}
 g_\alpha(R)=\E_\alpha[1-e^{-\alpha m_R(K_R)}].
\end{equation}
Choose $C_B$ so that $\capc(B_R)\le C_B R^2$ for $R\ge1$, and put
$M_0:=C_{\mathrm{att}}C_B$. Then $m_R(K_R)\le M_0$ almost surely.
The elementary inequalities
\[
 e^{-\bar\alpha M_0}\alpha u\le1-e^{-\alpha u}\le\alpha u,
 \qquad 0\le u\le M_0,\quad 0\le\alpha\le\bar\alpha,
\]
combine with \eqref{end:attach-mass} to give
\begin{equation}\label{end:attach-comparison}
 c_g\alpha R^{-2}A_\alpha(R)
 \le g_\alpha(R)\le C_g\alpha R^{-2}A_\alpha(R),
 \qquad R\ge32,
\end{equation}
where one may take
$c_g=e^{-\bar\alpha M_0}c_{\mathrm{att}}$ and
$C_g=C_{\mathrm{att}}$. On $F_R$, the attaching loop belongs to the
full origin cluster and leaves $B_R$, so $p_\alpha(R)\ge g_\alpha(R)$.
This proves the desired lower comparison when $R\ge32$. When
$1\le R<32$, use $p_\alpha(R)\ge g_\alpha(32)\ge c\alpha$ and
$R^{-2}A_\alpha(R)\le C_B$ to adjust its constant. The constants in
\eqref{end:comparison} are thus uniform throughout the stated range.
\end{proof}

\subsection{A positive logarithmic power in the capacity}
The attachment event contains more information than the one-arm event:
it adds an outer piece whose capacity is bounded below at the natural
four-dimensional walk scale. The following tail integration exploits
that information. It is the four-dimensional counterpart of the
capacity integration in Gao~\cite[Section~7.3]{Gao}. A positive logarithmic correction to the four-dimensional one-arm lower bound and the corresponding finite-box expected-capacity lower bound were already established by Chang and Sapozhnikov~\cite[Theorem~1.6 and Lemma~6.2]{CS16}. We include the argument to make the present normalization and the uniform intensity dependence explicit; the lower-order phenomenon itself is not new.

\begin{proposition}\label{end:cap-lower}
There are constants $c_{\mathrm{low}},\beta>0$ such that
\begin{equation}\label{end:cap-lower-bound}
 A_\alpha(R)\ge
 c_{\mathrm{low}}[\log(R+2)]^{\beta\alpha},
 \qquad 0\le\alpha\le\bar\alpha,\quad R\ge1.
\end{equation}
\end{proposition}
At each dyadic spatial scale, a good attachment gives a lower bound on one tail probability of the capacity in a larger box. The corresponding capacity thresholds grow geometrically, allowing us to integrate these tail bounds over disjoint intervals of capacity values. This produces a discrete integral inequality for the expected capacities, whose equality case is an explicit product. The argument therefore combines information from many scales without assuming that their attachment events are independent.
\begin{proof}
Keep $c_0=\capc(\{0\})$. Choose an integer $j_0$ so large that
$r_0:=2^{j_0}\ge32$ and
$\kappa_0r_0^2/(2\log r_0)>c_0$. For $j\ge0$ put
\[
 r_j:=2^jr_0,\qquad
 A_j:=A_\alpha(r_j),\qquad
 b_j:=\frac{\kappa_0r_j^2}{\log r_j}.
\]
These definitions imply
\[
 \frac{b_{j+1}}{b_j}
 =4\frac{j+j_0}{j+j_0+1}\ge2.
\]
Hence the intervals $[0,c_0)$ and $[b_j/2,b_j)$, $j\ge0$, are
pairwise disjoint.

On the good-attachment event $F_{r_j}$ used above, the loops forming
$K_{r_j}$ and the attaching loop are all contained in $B_{2r_j}$.
Their union is connected to $0$ and contains a set of capacity at
least $b_j$. Thus
\[
 F_{r_j}\subset\{\capc(K_{2r_j})\ge b_j\}.
\]
If $j<n$, then $2r_j\le r_n$, and the same witnessing loops lie in
$B_{r_n}$. The coupling by restriction of one common soup therefore
gives
\[
 \Pp_\alpha(\capc(K_{r_n})\ge b_j)\ge g_\alpha(r_j),
 \qquad j<n.
\]
Apply the tail-integral identity
$\E[X]=\int_0^\infty\Pp(X>u)\,du$ to
$X=\capc(K_{r_n})$. Integrating only over the disjoint intervals
specified above, and then using \eqref{end:attach-comparison}, gives
\begin{align}\label{end:lower-recursion}
 A_n
 &\ge c_0+\frac12\sum_{j=0}^{n-1}b_jg_\alpha(r_j)\notag\\
 &\ge c_0+
       \frac{\kappa_0c_g}{2}\alpha
       \sum_{j=0}^{n-1}\frac{A_j}{\log r_j}
  =c_0+\beta\alpha\sum_{j=0}^{n-1}\frac{A_j}{j+j_0},
 \qquad
 \beta:=\frac{\kappa_0c_g}{2\log2}>0.
\end{align}
The length of the $j$th integration interval is $b_j/2$, and its tail probability is at least $g_\alpha(r_j)$ throughout; this explains each summand in the first line. Only disjoint intervals of the capacity variable were used. The events
$F_{r_j}$ need not be disjoint or independent.

Define $D_0:=c_0$ and, for $n\ge1$,
\[
 D_n:=c_0\prod_{j=0}^{n-1}
             \left(1+\frac{\beta\alpha}{j+j_0}\right).
\]
Since $D_{n+1}-D_n=\beta\alpha D_n/(n+j_0)$, telescoping shows
\[
 D_n=c_0+\beta\alpha\sum_{j=0}^{n-1}\frac{D_j}{j+j_0}.
\]
Induction in \eqref{end:lower-recursion} consequently gives $A_n\ge D_n$.
For $x\ge0$, the inequality $\log(1+x)\ge x-x^2/2$ implies
\begin{align*}
 \log(D_n/c_0)
 &\ge \beta\alpha\sum_{j=0}^{n-1}\frac1{j+j_0}
       -\frac{\beta^2\alpha^2}{2}
                       \sum_{j=0}^{n-1}\frac1{(j+j_0)^2}\\
 &\ge \beta\alpha\log\frac{n+j_0}{j_0}-C.
\end{align*}
The constant $C$ is finite and uniform for $\alpha\le\bar\alpha$.
Thus $A_n\ge c(n+j_0)^{\beta\alpha}$, with a constant independent
of $n$ and $\alpha$. Since
$\log r_n=(n+j_0)\log2$, this is
$A_\alpha(r_n)\ge c[\log r_n]^{\beta\alpha}$.
For $R\ge r_0$, choose $n$ with $r_n\le R<r_{n+1}$ and use
$A_\alpha(R)\ge A_\alpha(r_n)$ and
$\log(R+2)\le C\log r_n$. Boundedness of $\alpha$ absorbs the
fixed comparison factor into $c$. For $1\le R<r_0$, use
$A_\alpha(R)\ge c_0$ and decrease the lower constant if necessary. This
proves \eqref{end:cap-lower-bound}, including $\alpha=0$.
\end{proof}

\subsection{Proof of the main theorem and strict subcriticality}
Propositions~\ref{end:cap-upper} and \ref{end:cap-lower} prove the
capacity bounds in \eqref{eq:maincap}. Proposition~\ref{end:arm-comparison}
then gives, for constants $c,C,\beta>0$ independent of $R$ and $\alpha$,
\begin{equation}\label{end:final-arm}
 c\alpha R^{-2}[\log(R+2)]^{\beta\alpha}
 \le p_\alpha(R)
 \le C\alpha R^{-2}[\log(R+2)]^2,
 \qquad 0<\alpha\le\alpha_*,\quad R\ge1.
\end{equation}
For the constants in Theorem~\ref{thm:main}, we may take
\[
 c_A=c_{\mathrm{low}},\qquad C_A=C_{\mathrm{cap}},\qquad
 c_p=c_{\mathrm{arm}}c_{\mathrm{low}},\qquad
 C_p=C_{\mathrm{arm}}
       \left(\frac{c_0}{[\log3]^2}+C_{\mathrm{cap}}\alpha_*\right).
\]
The exponent $\beta$ is the value defined in
\eqref{end:lower-recursion}. These choices follow by inserting
\eqref{end:cap-bound} and \eqref{end:cap-lower-bound} into
\eqref{end:comparison} and using $\log(R+2)\ge\log3$.

Fix now $0<\alpha\le\alpha_0$. Taking logarithms in
\eqref{end:final-arm} shows
\[
 \beta\alpha\log\log R-O_\alpha(1)
 \le \log\bigl(R^2p_\alpha(R)\bigr)
 \le 2\log\log R+O_\alpha(1),\qquad R\longrightarrow\infty.
\]
Here $O_\alpha(1)$ denotes a quantity bounded in $R$, with a bound
allowed to depend on the fixed positive intensity. It follows that
\[
 \log\bigl(R^2p_\alpha(R)\bigr)
      =\Theta_\alpha(\log\log R),
 \qquad
 p_\alpha(R)=R^{-2+o(1)}.
\]
The latter assertion means
$\lim_{R\to\infty}\log p_\alpha(R)/\log R=-2$. The constants in
the displayed two-sided inequalities are uniform for small intensity,
but these logarithmic asymptotic statements hold with $\alpha$ fixed.

Finally, \eqref{end:final-arm} also holds at
$\alpha_*=2\alpha_0$. Its upper bound tends to zero as
$R\to\infty$. Since the event that the origin belongs to an infinite
cluster is the decreasing intersection of the events
$\{0\leftrightarrow B_R^c\}$, continuity of probability gives
\[
 \Pp_{\alpha_*}(|\mathcal C_{\alpha_*}(0)|=\infty)=0.
\]
Monotonicity in intensity implies
$\alpha_c\ge\alpha_*=2\alpha_0>\alpha_0$. Thus the interval asserted
in the theorem is strictly below the percolation threshold. This
completes the proof of Theorem~\ref{thm:main}.

The estimates determine the exponent two and the order of the
logarithmic correction for sufficiently small fixed intensity. They
do not identify a limiting coefficient multiplying $\log\log R$,
and they do not extend the conclusion to every subcritical intensity. The theorem also does not address the conditional loop distance from the origin to a distant box boundary. In particular, the loop-distance question of Chang and Sapozhnikov~\cite[Section~8]{CS16} is a different question from the logarithmic correction to the one-arm probability considered here.
\subsection{Cluster-size and capacity tails}\label{sec:cluster-tails}

\begin{proof}[Proof of Corollary~\ref{cor:cluster-tails}]
We first obtain the additional finite-box estimate
\begin{equation}\label{end:volume-mean}
 \E_\alpha[|K_R|-1]\le C\alpha[\log(R+2)]^3,
 \qquad 0\le\alpha\le\alpha_0,\quad R\ge1.
\end{equation}
This uses the same pruning, Mecke and BK argument as
Proposition~\ref{chain:path-series}, with a different final test.
For each $x\in B_R$, take the constant trajectory $w_n=x$ for all
$n\ge0$ in the pointwise bound \eqref{chain:tuple-indicator}.
Sum over $x\in B_R$ after taking the soup expectation. The left-hand
side becomes $\E_\alpha|C_{R,T}(0)|$, the singleton term becomes one,
and the terminal indicator has sum
\[
 \sum_{x\in B_R}\one_{\{x\in U(\ell_m)\}}=|U(\ell_m)|.
\]
Thus \eqref{chain:series} holds with the two capacity terms replaced
by $|C_{R,T}(0)|$ and $|U(\ell_m)|$, and its constant term replaced
by one. This step is made in the same finite soup as before; all sums
and integrals are nonnegative. Constant trajectories need not have
positive random-walk probability, because the underlying connectivity
bound was proved pointwise for every trajectory.

Take $T=\max\{4,\lceil R^4\rceil\}$ and
$J=\lfloor\log_2T\rfloor$, as in Section~\ref{sec:conclusion}.
For every marked loop of duration $N\le T$, its charge satisfies
\[
 |U(\ell)|\le N\le\log(T+2)w(\ell).
\]
For a good loop this follows from $w=N/\log(N+2)$; for a bad loop,
it follows from $w=N^2$ and $\log(T+2)>1$. The backward integrations
in \eqref{end:backward-functions}--\eqref{end:fixed-word-bound}
therefore apply with terminal bound $\log(T+2)w$ in place of $C_ww$.
Using the numerical sum estimated in \eqref{end:truncated-cap}, we get
\begin{align*}
 \E_\alpha[|C_{R,T}(0)|-1]
 &\le\log(T+2)\sum_{i=1}^J
                   b_i(\alpha)S_{J,C_K\alpha}(i)\\
 &\le C\alpha\log(T+2)J^2
 \le C\alpha[\log(R+2)]^3.
\end{align*}
Restoring the omitted long loops changes the cardinality by at most
$|B_R|$ and only if at least one such loop occurs. Their measure
$\lambda_{R,T}$ was bounded after \eqref{end:long-mass}, so
\[
 0\le\E_\alpha[|K_R|-|C_{R,T}(0)|]
 \le |B_R|\alpha\lambda_{R,T}
 \le C\alpha R^6e^{-cR^2}\le C\alpha.
\]
This proves \eqref{end:volume-mean}.

On the event $\rho_\alpha\le R$, the full cluster equals $K_R$:
every loop belonging to that cluster has its entire range in $B_R$.
For $t>1$, Markov's inequality applied to $|K_R|-1$ consequently gives
\[
 \Pp_\alpha(|\mathcal C_\alpha(0)|>t)
 \le p_\alpha(R)
       +\frac{\E_\alpha[|K_R|-1]}{t-1}.
\]
Choose $R=\sqrt t$ and apply \eqref{eq:mainarm} and
\eqref{end:volume-mean}. This proves the upper bound in
\eqref{eq:full-volume-tail}; subtracting the singleton preserves the
factor $\alpha$ uniformly as $\alpha\downarrow0$.
For capacity, \eqref{in:box-capacity} implies that
$\capc(\mathcal C_\alpha(0))>t$ requires
$\rho_\alpha>c\sqrt t$ for a sufficiently small fixed $c>0$ and
sufficiently large $t$. The upper one-arm bound proves the upper
bound in \eqref{eq:full-capacity-tail}.

For the lower bounds, use the event $F_R$ from
\eqref{end:attach-probability}, based on Gao's attachment estimate
(Proposition~\ref{in:attachment}). On this event the full cluster
contains a set of capacity at least $\kappa_0R^2/\log R$. Thus
\eqref{end:attach-comparison} gives
\[
 \Pp_\alpha\left(
  \capc(\mathcal C_\alpha(0))\ge\frac{\kappa_0R^2}{\log R}\right)
 \ge g_\alpha(R)
 \ge c_g\alpha R^{-2}A_\alpha(R),\qquad R\ge32.
\]
Set $R=\lceil M\sqrt{t\log(t+2)}\rceil$, with a fixed sufficiently
large $M$. For all sufficiently large $t$, the displayed capacity
threshold is strictly greater than $t$, while
$R^2\asymp t\log(t+2)$ and $\log R\asymp\log(t+2)$.
The lower bound in \eqref{eq:maincap} proves the lower capacity-tail
bound. Its constants are uniform because
$0<\beta\alpha\le\beta\alpha_0$. The same lower bound holds for
cardinality, since $\capc(A)\le|A|$ for every finite set $A$.

Finally, for every nonnegative random variable $X$ and $s>0$,
Tonelli's theorem gives
\[
 \E[X^s]=s\int_0^\infty t^{s-1}\Pp(X>t)\,dt.
\]
The two tail bounds give finite cardinality and capacity moments for
$s<1$, and infinite ones for $s>1$. At $s=1$, the lower bound gives
the divergent integral
$\int^\infty t^{-1}(\log t)^{\beta\alpha-1}\,dt$ because
$\beta\alpha>0$. The radius statement follows from
$\Pp_\alpha(\rho_\alpha>R)=p_\alpha(R)$ and
\eqref{eq:mainarm}; the lower bound also makes its second moment
infinite. This proves all the assertions.
\end{proof}

\section{An application to random wires}\label{sec:wires}

The theorem also describes the support clusters of a member of the
random-wire family of Benassi and Ueltschi~\cite[Section~3.2,
equations~(3.3)--(3.4)]{BU}. We specify the member and the parameter
conversion, since the wire loop weight differs from the soup intensity.

Let $D\subset\Z^4$ be finite, let $E_D$ be its set of unoriented
nearest-neighbor edges, and let $P_D=(P(x,y))_{x,y\in D}$. Thus the walk
is killed on leaving $D$: boundary transition probabilities are not
renormalized, and $P_D$ has spectral radius less than one.
A wire configuration consists of nonnegative integers
$n=(n_e)_{e\in E_D}$ and pairings $\pi=(\pi_x)_{x\in D}$.
The $n_e$ links on edge $e$ have labels $1,\ldots,n_e$, and $\pi_x$
pairs all link endpoints incident to $x$. In particular,
$\sum_{e\ni x}n_e$ must be even. Alternating links and vertex pairings
partitions the links into unoriented closed circuits; let
$k(n,\pi)$ be their number. The empty configuration has one empty
pairing and no circuits.
For any real $q>0$, put
\begin{equation}\label{wire:law}
 \mathbb W_{q,D}(n,\pi)
 =\frac{q^{k(n,\pi)}}{Z_q(P_D)}
       \prod_{e\in E_D}\frac{(1/8)^{n_e}}{n_e!},
\end{equation}
where $Z_q(P_D)$ is the sum of the displayed weights over all
configurations. In the notation of~\cite{BU}, this is the model with
vertex potential $U\equiv0$, edge constants $J_e=1/8$, and loop
parameter $q$. Its finiteness follows from the next proposition.

\begin{proposition}[Wire--soup correspondence]\label{wire:identity}
For every $q>0$,
\begin{equation}\label{wire:partition}
 Z_q(P_D)=\det(I-P_D)^{-q/2}.
\end{equation}
The edge counts $n$ under $\mathbb W_{q,D}$ have the same law as the
total unoriented edge-traversal counts of $\mathcal L_{q/2}^{D}$.
Consequently, the connected components of the graph
$(D,\{e:n_e>0\})$ have the same law as the soup clusters in $D$.
\end{proposition}

\begin{proof}
For a positive integer $q$, let $(X_x)_{x\in D}$ be independent standard
Gaussian vectors in $\R^q$. Gaussian integration gives
\[
 \E\exp\left\{\sum_{\{x,y\}\in E_D}P_D(x,y)X_x\cdot X_y\right\}
 =\det(I-P_D)^{-q/2}.
\]
Expand first near zero in the edge weights. Wick's formula pairs the
factors at each vertex; a coordinate index is constant along each
paired circuit and has $q$ choices. The resulting coefficients are
exactly the weights in \eqref{wire:law}. For each fixed edge monomial,
both this pairing coefficient and the corresponding coefficient of
\[
 \exp\left\{\frac q2\sum_{j\ge1}\frac{\tr(P_D^j)}{j}\right\}
\]
are polynomials in $q$. Equality at all positive integers proves their
equality for every real $q>0$. Replacing $P_D$ by $sP_D$, the last
series is analytic for $|s|<\rho(P_D)^{-1}>1$ (with infinite radius
if $P_D=0$). Its nonnegative coefficients therefore sum at $s=1$
to the finite partition function in \eqref{wire:partition}.

For $z=(z_e)_{e\in E_D}\in[0,1]^{E_D}$, define $P_D(z)(x,y)$ to be
$P_D(x,y)z_{\{x,y\}}$ on adjacent pairs and zero otherwise. The wire
generating function is
\begin{equation}\label{wire:pgf}
 \E_{\mathbb W_{q,D}}\prod_{e\in E_D}z_e^{n_e}
 =\left(\frac{\det(I-P_D)}{\det(I-P_D(z))}\right)^{q/2}.
\end{equation}
By \eqref{eq:model-rooted}, the loop measure has total mass
$\sum_{j\ge1}\tr(P_D^j)/j$. Inserting the edge weights $z$ in this
identity and applying the Poisson exponential formula gives exactly
\eqref{wire:pgf} for the traversal counts of $\mathcal L_{q/2}^{D}$.
This identifies their joint law, including backtracking and periodic
loops. Finally, each loop has connected edge support, and loops meeting
at a vertex belong to the same support component. This proves the
cluster assertion.
\end{proof}

Write $\mathbb W_q$ for the infinite-volume law of the edge counts
obtained by increasing finite domains to $\Z^4$. This limit exists:
the preceding identity couples these finite-domain laws as counts
of complete soup loops contained in the domains, and these counts
increase to the full-soup counts. For $e=\{x,y\}$ their limiting
expectation, obtained by differentiating \eqref{wire:pgf} and passing
to the limit, is $qP(x,y)G(x,y)<\infty$, so every edge count is finite
almost surely. This definition requires only the current marginal,
not a choice of pairings at infinity. Let $C_q(0)$ denote its origin
support component, including $0$ when it is isolated, and let
$C_{q,D}(0)$ be the corresponding component under $\mathbb W_{q,D}$
when $0\in D$. The notation $0\leftrightarrow B_R^c$ below refers to
connectivity in this occupied-edge graph.

\begin{corollary}[Small loop-weight wire clusters]\label{wire:corollary}
Let $\alpha_0,\beta,c_A,C_A,c_p,C_p$ be the constants of
Theorem~\ref{thm:main}. For $0<q\le2\alpha_0$ and $R\ge1$,
\begin{align}
 \frac{c_pq}{2}R^{-2}[\log(R+2)]^{\beta q/2}
 &\le\mathbb W_q(0\leftrightarrow B_R^c)
 \le\frac{C_pq}{2}R^{-2}[\log(R+2)]^2,
 \label{wire:arm}\\
 c_A[\log(R+2)]^{\beta q/2}
 &\le\E_{\mathbb W_{q,B_R}}\capc(C_{q,B_R}(0))
 \le c_0+\frac{C_Aq}{2}[\log(R+2)]^2.
 \label{wire:capacity}
\end{align}
In particular, for every fixed $q$ in this interval, the wire support
has one-arm probability $R^{-2+o(1)}$ and no infinite component.
\end{corollary}

\begin{proof}
Apply Theorem~\ref{thm:main} with $\alpha=q/2$ and use
Proposition~\ref{wire:identity} and the increasing-domain construction.
Absence of an infinite origin component follows by letting $R$ tend
to infinity in \eqref{wire:arm}; translation invariance and a countable
union give absence of all infinite components.
\end{proof}

Capacity in \eqref{wire:capacity} is full-space capacity of a cluster
sampled in the killed domain $B_R$. Restricting the edges of the
infinite-volume current to $B_R$ would give a different cluster.
The corollary concerns the zero-vertex-potential model
\eqref{wire:law}; nonzero vertex potentials in the general wire family,
including those used for fixed-length spin models, are not covered by
this identification.

\section{AI Usage}
I have worked on this problem over the course of many years. I gave my notes and ideas to ChatGPT and asked it for help, resulting eventually in this preprint. Some of the ideas presented here are mine, others are supplied by ChatGPT, and others are joint. I do not think I could have finished this project without the use of AI tools. 
\bigskip
\begingroup
\raggedright

\endgroup


\begin{thebibliography}{99}

\bibitem{ASS}
A.~Asselah, B.~Schapira and P.~Sousi,
\emph{Capacity of the range of random walk on $\mathbb Z^4$},
Annals of Probability \textbf{47} (2019), no.~3, 1447--1497.
\href{https://doi.org/10.1214/18-AOP1288}{doi:10.1214/18-AOP1288}.
See equation~(1.4) for the mean-capacity asymptotic.

\bibitem{BBS}
R.~Bauerschmidt, D.~C.~Brydges and G.~Slade,
\emph{Logarithmic correction for the susceptibility of the 4-dimensional
weakly self-avoiding walk: a renormalisation group analysis},
Communications in Mathematical Physics \textbf{337} (2015), 817--877.
\href{https://doi.org/10.1007/s00220-015-2352-6}{doi:10.1007/s00220-015-2352-6}.

\bibitem{BU}
C.~Benassi and D.~Ueltschi,
\emph{Loop correlations in random wire models},
Communications in Mathematical Physics \textbf{374} (2020), 525--547.
\href{https://doi.org/10.1007/s00220-019-03474-9}{doi:10.1007/s00220-019-03474-9}.

\bibitem{BDKM26}
V.~Betz, A.~Drewitz, A.~Klippel and C.~M\"onch,
\emph{(Non-)coincidence of critical parameters for Poisson Zoos and Loop
Soup Percolation on $\mathbb Z^d$, $d>4$, and $T_d$, $d\ge3$},
preprint, arXiv:2608.11198, 2026.
\url{https://arxiv.org/abs/2608.11198}.

\bibitem{CDlow}
Z.~Cai and J.~Ding,
\emph{One-arm probabilities for metric graph Gaussian free fields below
and at the critical dimension},
preprint, arXiv:2406.02397v2, 2024.
\url{https://arxiv.org/abs/2406.02397}.

\bibitem{Cha17}
Y.~Chang,
\emph{Supercritical loop percolation on $\mathbb Z^d$ for $d\ge3$},
Stochastic Processes and their Applications \textbf{127} (2017),
no.~10, 3159--3186.
\href{https://doi.org/10.1016/j.spa.2017.02.003}{doi:10.1016/j.spa.2017.02.003}.

\bibitem{CDL}
Y.~Chang, H.~Du and X.~Li,
\emph{Percolation threshold for metric graph loop soup},
Bernoulli \textbf{30} (2024), no.~4, 3324--3333.
\href{https://doi.org/10.3150/23-BEJ1716}{doi:10.3150/23-BEJ1716}.

\bibitem{CS16}
Y.~Chang and A.~Sapozhnikov,
\emph{Phase transition in loop percolation},
Probability Theory and Related Fields \textbf{164} (2016),
no.~3--4, 979--1025.
\href{https://doi.org/10.1007/s00440-015-0624-x}{doi:10.1007/s00440-015-0624-x}.

\bibitem{DL18}
J.~Ding and L.~Li,
\emph{Chemical distances for percolation of planar Gaussian free fields
and critical random walk loop soups},
Communications in Mathematical Physics \textbf{360} (2018),
no.~2, 523--553.
\href{https://doi.org/10.1007/s00220-018-3140-x}{doi:10.1007/s00220-018-3140-x}.

\bibitem{DW}
J.~Ding and M.~Wirth,
\emph{Percolation for level-sets of Gaussian free fields on metric graphs},
Annals of Probability \textbf{48} (2020), no.~3, 1411--1435.
\href{https://doi.org/10.1214/19-AOP1397}{doi:10.1214/19-AOP1397}.

\bibitem{DPRcap}
A.~Drewitz, A.~Pr\'evost and P.-F.~Rodriguez,
\emph{Cluster capacity functionals and isomorphism theorems for Gaussian
free fields},
Probability Theory and Related Fields \textbf{183} (2022), 255--313.
\href{https://doi.org/10.1007/s00440-021-01090-0}{doi:10.1007/s00440-021-01090-0}.

\bibitem{DPRarm}
A.~Drewitz, A.~Pr\'evost and P.-F.~Rodriguez,
\emph{Arm exponent for the Gaussian free field on metric graphs in
intermediate dimensions},
Annals of Probability \textbf{54} (2026), no.~2, 644--668.
\href{https://doi.org/10.1214/25-AOP1775}{doi:10.1214/25-AOP1775}.

\bibitem{Gao}
Y.~Gao,
\emph{On loop-soup clusters in three dimensions and higher},
preprint, arXiv:2610.38510v1, 28 October 2026.
\url{https://arxiv.org/abs/2610.38510v1}.

\bibitem{Hoeffding}
W.~Hoeffding,
\emph{Probability inequalities for sums of bounded random variables},
Journal of the American Statistical Association \textbf{58} (1963),
no.~301, 13--30.
\href{https://doi.org/10.1080/01621459.1963.10500830}{doi:10.1080/01621459.1963.10500830}.

\bibitem{LP}
G.~Last and M.~Penrose,
\emph{Lectures on the Poisson Process},
Institute of Mathematical Statistics Textbooks, vol.~7,
Cambridge University Press, Cambridge, 2018.
\href{https://doi.org/10.1017/9781316104477}{doi:10.1017/9781316104477}.

\bibitem{LL}
G.~F.~Lawler and V.~Limic,
\emph{Random Walk: A Modern Introduction},
Cambridge Studies in Advanced Mathematics, vol.~123,
Cambridge University Press, Cambridge, 2010.
\href{https://doi.org/10.1017/CBO9780511750854}{doi:10.1017/CBO9780511750854}.

\bibitem{LTF07}
G.~F.~Lawler and J.~A.~Trujillo Ferreras,
\emph{Random walk loop soup},
Transactions of the American Mathematical Society \textbf{359} (2007),
no.~2, 767--787.
\href{https://doi.org/10.1090/S0002-9947-06-03916-X}{doi:10.1090/S0002-9947-06-03916-X}.

\bibitem{LW04}
G.~F.~Lawler and W.~Werner,
\emph{The Brownian loop soup},
Probability Theory and Related Fields \textbf{128} (2004),
no.~4, 565--588.
\href{https://doi.org/10.1007/s00440-003-0319-6}{doi:10.1007/s00440-003-0319-6}.

\bibitem{LeJan}
Y.~Le Jan,
\emph{Markov Paths, Loops and Fields},
Lecture Notes in Mathematics, vol.~2026,
Springer, Heidelberg, 2011.
\href{https://doi.org/10.1007/978-3-642-21216-1}{doi:10.1007/978-3-642-21216-1}.

\bibitem{LJL13}
Y.~Le Jan and S.~Lemaire,
\emph{Markovian loop clusters on graphs},
Illinois Journal of Mathematics \textbf{57} (2013),
no.~2, 525--558.
\href{https://doi.org/10.1215/ijm/1408453593}{doi:10.1215/ijm/1408453593}.

\bibitem{Lupu}
T.~Lupu,
\emph{From loop clusters and random interlacements to the free field},
Annals of Probability \textbf{44} (2016), no.~3, 2117--2146.
\href{https://doi.org/10.1214/15-AOP1019}{doi:10.1214/15-AOP1019}.

\bibitem{LW16}
T.~Lupu and W.~Werner,
\emph{A note on Ising random currents, Ising-FK, loop-soups and the
Gaussian free field},
Electronic Communications in Probability \textbf{21} (2016),
paper no.~13, 1--7.
\href{https://doi.org/10.1214/16-ECP4733}{doi:10.1214/16-ECP4733}.

\bibitem{PV}
S.~Popov and Q.~Vogel,
\emph{On large deviations for the range of a two-dimensional random walk},
Electronic Communications in Probability \textbf{31} (2026), 1--6.
\href{https://doi.org/10.1214/26-ECP817}{doi:10.1214/26-ECP817}.

\bibitem{SW12}
S.~Sheffield and W.~Werner,
\emph{Conformal loop ensembles: the Markovian characterization and the
loop-soup construction},
Annals of Mathematics (2) \textbf{176} (2012),
no.~3, 1827--1917.
\href{https://doi.org/10.4007/annals.2012.176.3.8}{doi:10.4007/annals.2012.176.3.8}.

\bibitem{BK}
J.~van den Berg and H.~Kesten,
\emph{Inequalities with applications to percolation and reliability},
Journal of Applied Probability \textbf{22} (1985), no.~3, 556--569.
\href{https://doi.org/10.2307/3213860}{doi:10.2307/3213860}.

\bibitem{Vog20}
Q.~Vogel,
\emph{Geometric properties of random walk loop soups},
PhD thesis, University of Warwick, 2020.
\url{http://wrap.warwick.ac.uk/154434}.

\bibitem{Vog26}
Q.~Vogel,
\emph{$1+O(1)$ asymptotics for loop percolation in five and higher dimensions},
International Mathematics Research Notices \textbf{2026} (2026),
no.~9, rnag084, 18~pp.
\href{https://doi.org/10.1093/imrn/rnag084}{doi:10.1093/imrn/rnag084}.

\end{thebibliography}
\end{document}